\documentclass{amsart}
\usepackage{{amsmath,amsfonts,latexsym, amstext, amssymb}}

\newtheorem{definition}{Definition}[section]
\newtheorem{prop}[definition]{Proposition}
\newtheorem{thm}[definition]{Theorem}
\newtheorem{corollary}[definition]{Corollary}

\newtheorem{lemma}[definition]{Lemma}

\newcommand{\E}{\mathcal{E}}
\newcommand{\F}{\mathcal{F}}

\newcommand{\la}{\langle}
\newcommand{\ra}{\rangle}

\newcommand{\n}{\nabla}
\newcommand{\ep}{\epsilon}

\newcommand{\CE}{\mathcal{E}}

\newcommand{\CM}{\mathbb{C} \mathbb{M}}
\newcommand{\CP}{\mathbb{C} \mathbb{P}}
\newcommand{\CC}{\mathbb{C}}
\newcommand{\CO}{\mathcal{O}}
\newcommand{\Z}{\mathbb{Z}}

\newcommand{\R}{\mathbb{R}}
\newcommand{\D}{\textup{d}}
\newcommand{\un}{\underline \n}
\newcommand{\tz}{\tilde{z}}
\newcommand{\tw}{\tilde{w}}
\newcommand{\ts}{\tilde{s}}

\newcommand{\uD}{\underline \D}
\newcommand{\TCM}{\widetilde{{\CM}^\#}}
\newcommand{\hCM}{{\CM}^\#}

\newcommand{\pA}{A^\prime}
\newcommand{\pB}{B^\prime}

\newcommand{\extd}{\textup{d}}

\newcommand{\pc}{c}
\newcommand{\pd}{d}
\newcommand{\mup}{\mu}
\newcommand{\nup}{\nu}
\newcommand{\alphap}{\alpha}
\newcommand{\betap}{\beta}
\renewcommand{\S}{\textup{S}}
\newcommand{\GL}{\textup{G}\textup{L}}
\newcommand{\PGL}{\textup{P}\textup{G}\textup{L}}
\newcommand{\SL}{\textup{S}\textup{L}}
\newcommand{\SU}{\textup{S}\textup{U}}
\newcommand{\U}{\textup{U}}
\newcommand{\p}{\partial}
\newcommand{\bu}{\bullet}

\renewcommand{\o}{{}_{\scriptscriptstyle(1)}}
\renewcommand{\t}{{}_{\scriptscriptstyle(2)}}
\newcommand{\thr}{{}_{\scriptscriptstyle(3)}}
\newcommand{\fo}{{}_{\scriptscriptstyle(4)}}

\newcommand{\tens}{{\otimes}}
\newcommand{\eps}{{\epsilon}}
\newcommand{\del}{{\partial}}
\newcommand{\C}{{\mathbb C}}
\newcommand{\Tr}{{\mathrm{Tr}}}
\renewcommand{\Im}{\textup{Im}}

\begin{document}

\title{Quantisation of twistor theory by cocycle twist}
\author{S.J. Brain}
\address{Mathematical Institute\\ 24-29 St Giles', Oxford OX1 3LB, UK}
\author{S. Majid}
\address{School of Mathematical Sciences\\
Queen Mary, University of London\\ 327 Mile End Rd,  London E1
4NS, UK}
\thanks{The work was mainly completed while S.M. was visiting July-December 2006 at the Isaac Newton Institute, Cambridge,  which both authors thank for support. }


\date{January, 2007}%

\maketitle

\begin{abstract} We present the main ingredients of twistor theory leading up to and including the
Penrose-Ward transform in a coordinate algebra form which we can
then `quantise' by means of a functorial cocycle twist. The quantum
algebras for the conformal group, twistor space $\CP^3$,
compactified Minkowski space $\hCM$ and the twistor correspondence
space are obtained along with their canonical quantum differential
calculi, both in a local form and in a global $*$-algebra formulation which even in the classical commutative case provides a useful alternative to the formulation in terms of projective varieties. We outline how the Penrose-Ward transform then quantises. As an example, we show that the pull-back of the tautological bundle on $\hCM$ pulls back to the basic instanton on $S^4\subset\hCM$ and that this observation quantises to obtain the Connes-Landi instanton on $\theta$-deformed $S^4$ as the pull-back of the tautological bundle on our $\theta$-deformed $\hCM$. We likewise quantise the fibration $\CP^3\to S^4$ and use it to construct the bundle on $\theta$-deformed $\CP^3$ that maps over under the transform to the $\theta$-deformed instanton.
 \end{abstract}

\bibliographystyle{plain}

\section{Introduction and preliminaries}

There has been a lot of interest in recent years in the
`quantisation' of space-time (in which the algebra of coordinates
$x_\mu$ is noncommutative), among them one class of examples of the
Heisenberg form
\[ [x_\mu,x_\nu]=\imath\theta_{\mu\nu}\]
where the deformation parameter is an antisymmetric tensor or (when
placed in canonical form) a single parameter $\theta$. One of the
motivations here is from the effective theory of the ends of open
strings in a fixed D-brane\cite{sw:stng} and in this context a lot
of attention has been drawn to the existence of noncommutative
instantons and other nontrivial noncommutative geometry that
emerges, see \cite{ns:ins} and references therein to a large
literature. One also has $\theta$-versions of $S^4$ coming out of
considerations of cyclic cohomology in noncommutative geometry (used
to characterise what a noncommutative 4-sphere should be),  see
notably \cite{cl:id,CS}.

In the present paper we show that underlying and bringing together
these constructions is in fact a systematic theory of what could be
called $\theta$-deformed or `quantum' twistor theory. Thus we
introduce noncommutative versions of conformal complexified
space-time $\hCM$, of twistor space $\CP^3$ as well as of the
twistor correspondence space $\F_{12}$ of 1-2-flags in $\C^4$ used
in the Penrose-Ward transform \cite{pr:sst,rw:sdgf}. Our approach is
a general one but we do make contact for specific parameter values
with some previous ideas on what should be noncommutative twistor
space, notably with \cite{kko:nitt,kch:ncts} even though these works
approach the problem entirely differently. In our approach we
canonically find not just the noncommutative coordinate algebras but
their algebras of differential forms, indeed because our
quantisation takes the form of a `quantisation functor' we find in
principle the noncommutative versions of all suitably covariant
constructions. Likewise, inside our $\theta$-deformed $\hCM$ we find
(again for certain parameter values) exactly the $\theta$-deformed
$S^4$ of \cite{cl:id} as well as its differential calculus.

While the quantisation of twistor theory is our main motivation,
most of the present paper is in fact concerned with properly setting
up the classical theory from the `right' point of view after which
quantisation follows functorially. We provide in this paper {\em
two} classical points of view, both of interest. The first is purely
local and corresponds in physics to ordinary (complex) Minkowski
space as the flat `affine' part of $\hCM$. Quantisation at this
level gives the kind of noncommutative space-time mentioned above
which can therefore be viewed as a local `patch' of the actual
noncommutative geometry. The actual varieties $\hCM$ and $\CP^3$ are
however projective varieties and cannot therefore be simply
described by generators and relations in algebraic geometry, rather
one should pass to the `homogeneous coordinate algebras'
corresponding to the affine spaces $\widetilde{\hCM},
\widetilde{\CP^3}=\C^4$ that project on removing zero and
quotienting by an action of $\C^*$ to the projective varieties of
interest. Let us call this the `conventional approach'. We explain
the classical situation in this approach in Sections~1.1, 2 below,
and quantise it (including the relevant quantum group of conformal
transformations and the algebra of differential forms) in
Sections~4,5. The classical Sections~1.1,~2 here are not intended to
be anything new but to provide a lightning introduction to the
classical theory and an immediate coordinate algebra reformulation
for those unfamiliar either with twistors or with algebraic groups.
The quantum Sections~4,5 contain the new results in this stream of
the paper and provide a more or less complete solution to the basic
noncommutative differential geometry at the level of the quantum
homogeneous coordinate algebras $\C_F[\widetilde{\hCM}]$,
$\C_F[\widetilde{\CP^3}]$ etc. Here $F$ is a 2-cocycle which is the
general quantisation data in the cocycle twisting method
\cite{Ma:book,MaOec:twi} that we use.

Our second approach even to classical twistor theory is a novel one
suggested in fact from quantum theory. We call this the unitary or
$*$-algebraic formulation of our projective varieties $\hCM,\CP^3$
as real manifolds, setting aside that they are projective varieties.
The idea is that mathematically $\hCM$ is the Grassmannian of
2-planes in $\C^4$ and every point in it can therefore be viewed not
as a 2-plane but as a self-adjoint rank two projector $P$ that picks
out the two-plane as the eigenspace of eigenvalue 1. Working
directly with such projectors as a coordinatisation of $\hCM$, its
commutative coordinate $*$-algebra is therefore given by 16
generators $P^\mu{}_\nu$ with relations that $P.P=P$ as an
algebra-valued matrix, $\Tr~{P}=2$ and the $*$-operation
$P^\mu{}_\nu{}^*=P^\nu{}_\mu$. Similarly $\CP^3$ is the commutative
$*$-algebra with a matrix of generators $Q^\mu{}_\nu$, the relations
$Q.Q=Q$, $\Tr~{Q}=1$ and the $*$-operation
$Q^\mu{}_\nu{}^*=Q^\nu{}_\mu$. One may proceed similarly for all
classical flag varieties. The merit of this approach is that if one
forgets the $*$-structure one has  affine varieties defined simply
by generators and relations (they are  the complexifications of our
original projective varieties viewed as real manifolds), while the
$*$-structure picks out the real forms that are $\hCM,\CP^3$ as real
manifolds in our approach (these cannot themselves be described
simply by generators and relations). Finally, the complex structure
of our projective varieties appears now in real terms  as a
structure on the cotangent bundle. This amounts to a new approach to
projective geometry suggested by our theory for classical flag
varieties and provides a second stream in the paper starting in
Section~3. Note that there is no simple algebraic formula for change
of coordinates from describing a 2-plane as a 2-form and as a rank 2
projector, so the projector coordinates have a very different
flavour from those usually used for $\hCM,\CP^3$. For example the
tautological vector bundles in these coordinates are now immediate
to write down and we find that the pull-back of the tautological one on $\hCM$ to
a natural $S^4$ contained in it is exactly the instanton bundle given by the
known projector for $S^4$ (it is the analogue of the Bott projector that gives
the basic monopole bundle on $S^2$). We explain this calculation in
detail in Section~3.1. The Lorentzian version is also mentioned and
we find that Penrose's diamond compactification of Minkowski space
arises very naturally in these coordinates. In Section~3.2 we explain the known fibration $\CP^3\to S^4$ in our new approach, used to construct an auxiliary bundle that maps over under the
Penrose-Ward transform to the basic instanton.

The second merit of our approach is that just as commutative
$C^*$-algebras correspond to (locally compact) topological spaces,
quantisation has a precise meaning as a noncommutative $*$-algebra
with (in principle) $C^*$-algebra completion. Moreover, one does not
need to consider completions but may work at the $*$-algebra level,
as has been shown amply in the last two decades in the theory of
quantum groups \cite{Ma:book}. The quantisation of all flag
varieties, indeed of all varieties defined by `matrix' type
relations on a matrix of generators is given in Section~6, with the
quantum tautological bundle looked at explicitly in Section~6.1. Our
quantum algebra $\C_F[\hCM]$ actually has three independent real
parameters in the unitary case and takes a `Weyl form' with phase
factor commutation relations (see Proposition~6.3). We also show
that only a 1-parameter subfamily gives a natural quantum $S^4$ and
in this case we recover exactly the $\theta$-deformed $S^4$ and its
instanton as in \cite{cl:id,CS}, now from a different point of view
as `pull back' from our $\theta$-deformed $\hCM$.

Finally, while our main results are about the coordinate algebras
and differential geometry behind twistor theory in the classical and
quantum cases, we look in Section~7,8 at enough of the deeper theory
to see that our methods are compatible also with the Penrose-Ward
transform and ADHM construction respectively. In these sections we
concentrate on the classical theory but formulated in a manner that
is then `quantised' by our functorial method. Since their
formulation in noncommutative geometry is not fully developed we
avoid for example the necessity of the implicit complex structures.
We also expect our results to be compatible with another approach to
the quantum version based on groupoid $C^*$-algebras
\cite{sb:thesis}. Although we only sketch the quantum version, we do
show that our formulation includes for example the quantum basic
instanton as would be expected. A full account of the quantum
Penrose-Ward transform including an explicit treatment of the
noncommutative complex structure is deferred to a sequel.

\subsection{Conformal space-time}
\label{conformal space time}

Classically, complex Minkowski space $\CM$ is the four-dimensional
affine vector space $\CC^4$ equipped with the metric
$$ \D s^2 = 2(\D z \D \tz - \D w \D \tw)$$
written in double null coordinates \cite{mw:isdtt}. Certain
conformal transformations, such as isometries and dilations, are
defined globally on $\CM$, whereas others, such as inversions and
reflections, may map a light cone to infinity and \textit{vice
versa}.  In order to obtain a group of globally defined conformal
transformations, we adjoin a light cone at infinity to obtain
\textit{compactified} Minkowski space, usually denoted $\CM^\#$.

This compactification is achieved geometrically as follows (and is
just the Pl\"{u}cker embedding, see for example
\cite{mw:isdtt,ww:tgft,be:book}). One observes that the exterior
algebra $\Lambda^2 \CC^4$ can be identified with the set of $4
\times 4$ matrices as
$$x = \left( \begin{array}{cccc}
0 & s & -w & \tz \\
-s & 0 & -z & \tw \\
w & z & 0 & t \\
-\tz & -\tw & -t & 0 \end{array} \right),$$
the points of $\Lambda^2 \CC^4$ being identified with the six
entries $x^{\mu \nu}, \ \mu < \nu$.  Then $\GL_4 = \GL(4,\CC)$ acts
from the left on $\Lambda^2 \CC^4$ by conjugation,
$$x \mapsto a x a^t, \quad a \in \GL_4.$$
We note that multiples of the identity act trivially, and that this
action preserves the quadratic relation $\textup{det} \, x \equiv (
st - z \tz + w \tw)^2 = 0$.  From the point of view of $\Lambda^2
\CC^4$ this quadric, which we shall denote $\TCM$, is the subset of
the form $\{a\wedge b~:~a,b\in \CC^4\}\subset\Lambda^2\CC^4$, (the
antisymmetric projections of rank-one matrices, i.e. of decomposable
elements of the tensor product). We exclude $x=0$. Note that $x$ of
the form
\[\mbox{ $x = \left( \begin{array}{cccc}
0 & a_{11}a_{22} - a_{21}a_{12} & -(a_{31}a_{12} - a_{11}a_{32}) &  a_{11}a_{42}- a_{41}a_{12}  \\
-(a_{11}a_{22} - a_{21}a_{12}) & 0 &  -(a_{31}a_{22}-a_{21}a_{32})   & a_{21}a_{42}-a_{41}a_{22} \\
a_{31}a_{12} - a_{11}a_{32} &  a_{31}a_{22}-a_{21}a_{32}  & 0 & a_{31}a_{42}-a_{41}a_{32} \\
-( a_{11}a_{42}- a_{41}a_{12})  & -(a_{21}a_{42}-a_{41}a_{22}) &
-(a_{31}a_{42}-a_{41}a_{32}) & 0
\end{array} \right)$, }\]
or $x^{\mu\nu}=a^{[\mu}b^{\nu]}$ where $a=a_{\cdot 1}, b=a_{\cdot
2}$, automatically has determinant zero.  Conversely, if the
determinant vanishes then an antisymmetric matrix has this form over
$\CC$. To see this, we provide a short proof as follows. Thus, we
have to solve
\[ a_1b_2-a_2b_1=s,\quad a_1b_3-a_3b_1=-w,\quad a_1b_4-a_4b_1=\tz\]
\[ a_2 b_3-a_3 b_2=-z,\quad a_2 b_4-a_4b_2=\tw,\quad a_3b_4-a_4b_3=t.\]
We refer to the first relation as the (12)-relation, the second as
the (13)-relation and so forth. Now if a solution for $a_i,b_i$
exists, we make use of a `cycle' consisting of the $(12)b_3,
(23)b_1, (13)b_2$ relations (multiplied as shown) to deduce that
\[ a_1b_2b_3=a_2b_1b_3+sb_3=a_3b_1b_2+sb_3-zb_1=a_1b_3b_2+sb_3-zb_1+wb_2\]
hence a linear equation for $b$. The cycles consisting of the
$(12)b_4, (24)b_1, (14)b_2$ relations, the $(13)b_4, (34)b_1,
(14)b_3$ relations, and the $(23)b_4, (34)b_2, (24)b_3$ relations
give altogether the necessary conditions
\[ \left( \begin{array}{cccc}
0 & -s & -w & z \\
s & 0 & -\tz & \tw \\
w & \tz & 0 &- t \\
-z & -\tw & t & 0 \end{array} \right)\begin{pmatrix}b_4\\ b_3\\ b_2\\ b_1\end{pmatrix}=0.\]
The matrix here is not the matrix $x$ above but it has the same
determinant. Hence if $\det x=0$ we know that a nonzero vector $b$
obeying these necessary conditions must exist. We now fix such a
vector $b$, and we know that at least one of its entries must be
non-zero. We treat each case in turn. For example, if $b_2\ne 0$
then from the above analysis, the (12),(23) relations imply the (13)
relations. Likewise  $(12),(24) \Rightarrow (14)$,
$(23),(24)\Rightarrow (34)$. Hence the six original equations to be
solved become the three linear equations in four unknowns $a_i$:
\[ a_1b_2-a_2b_1=s,\quad a_2 b_3-a_3 b_2=-z,\quad a_2 b_4-a_4b_2=\tw\]
with general solution
\[ a=\lambda b+b_2^{-1}\begin{pmatrix} s\\ 0 \\ z\\ -\tw\end{pmatrix},\quad \lambda\in\CC.\]
One proceeds similarly in each of the other cases where a single
$b_i\ne 0$. Clearly, adding any multiple of $b$ will not change
$a\wedge b$, but we see that apart from this $a$ is uniquely fixed
by a choice of zero mode $b$ of a matrix with the same but permuted
entries as $x$. It follows that every $x$ defines a two-plane in
$\CC^4$ spanned by the obtained linearly independent vectors $a,b$.

Such matrices $x$ with $\det x=0$ are the orbit under $GL_4$ of the point where $s=1,
t=z=\tz=w=\tw=0$.  It is easily verified that this point has
isotropy subgroup $\tilde H$ consisting of elements of $GL_4$ such
that $a_{3 \mu} = a_{4 \mu} = 0$ for $\mu = 1,2$ and $a_{11}a_{22} -
a_{21}a_{12} = 1$. Thus $\TCM=GL_4/\tilde H$ where we quotient from the right.

Finally, we may identify conformal space-time $\CM^\#$ with the rays
of the above quadric cone $st = z \tz - w \tw$ in $\Lambda^2 \CC^4$,
identifying the finite points of space-time with the rays for which
$t \neq 0$ (which have coordinates $z, \tz, w, \tw$ up to scale):
the rays for which $t=0$ give the light cone at infinity. It follows
that the group $\textup{P} \GL(4,\CC) = \GL(4,\CC) / \CC$ acts
globally on $\CM^\#$ by conformal transformations and that every
conformal transformation arises in this way.  Observe that $\CM^\#$
is, in particular, the orbit of the point $s=1, z=\tz=w=\tw=t=0$
under the action of the conformal group $\textup{P} \GL(4,\CC)$.
Moreover, by the above result we have that $\hCM=\F_2(\CC^4)$, the
Grassmannian of two-planes in $\CC^4$.

We may equally identify $\CM^\#$ with the resulting quadric in the
projective space $\CP^5$ by choosing homogeneous coordinates
$s,z,\tz,w,\tw$ and projective representatives with $t=0$ and $t=1$.
In doing so, there is no loss of generality in identifying the
conformal group $\textup{P} \GL(4,\CC)$ with $\SL(4,\CC)$ by
representing each equivalence class with a transformation of unit
determinant.  Observing that $\Lambda^2 \CC^4$ has a natural metric
\[ \tilde \upsilon = 2(-\D s \D t + \D z \D \tz - \D w \D \tw), \]
we see that $\TCM$ is the null cone through the origin in $\Lambda^2
\CC^4$.  This metric may be restricted to this cone and moreover it
descends to give a metric $\upsilon$ on $\hCM$ \cite{pr:sst}.
Indeed, choosing a projective representative $t=1$ of the coordinate
patch corresponding to the affine piece of space-time, we have
$$\upsilon = 2(\D z \D \tz - \D w \D \tw),$$
thus recovering the original metric.  Similarly, we find the metric
on other coordinate patches of $\hCM$ by in turn choosing projective
representatives $s=1, z=1, \tz=1, w=1, \tw=1$.

Passing to the level of coordinates algebras let us denote by
$a^\mu_\nu$ the coordinate functions in $\CC[GL_4]$ (where we have
now rationalised indices so that they are raised and lowered by the
metric $\tilde \upsilon$) and by $s,t,z,\tz,w,\tw$ the coordinates
in $\CC[\Lambda^2\CC^4]$. The algebra $\CC[\Lambda^2\CC^4]$ is the
commutative polynomial algebra on the these six generators with no
further relations, whereas the algebra $\CC[\TCM]$ is the quotient
by the further relation $st-z\tz+w\tw=0$.  (Although we are
ultimately interested in the projective geometry of the space
described by this algebra, we shall put this point aside for the
moment). In the coordinate algebra (as an affine algebraic variety)
we do not see the deletion of the zero point in $\TCM$.

As explained, $\CC[\TCM]$ is essentially the algebra of functions on the orbit
of the point $s=1, t=z=\tz=w=\tw=0$ in $\Lambda^2 \CC^4$ under the
action of $\GL_4$. The specification of a $GL_4/\tilde H$ element that moves
the base point to a point of $\TCM$ becomes at the level of coordinate algebras
the map
\[\phi: \CC[\TCM]\cong \CC[GL_4]^{\CC[\tilde H]},\quad \phi(x^{\mu\nu})=a^\mu_1 a^\nu_2 - a^\nu_1 a^\mu_2. \]
As shown, the relation $st = z \tz - w \tw$ in $\CC[\Lambda^2
\CC^4]$ automatically holds for the image of the generators, so this
map is well-defined. Also in these dual terms there is a left coaction
\[ \Delta_L (x^{\mu\nu}) = a^\mu_\alpha
a^\nu_\beta \tens x^{\alpha \beta}   \]
of $\CC[GL_4]$ on $\CC[\Lambda^2\CC^4]$.  One should view the orbit
base point above as a linear function on $\CC[\Lambda^2\CC^4]$ that sends
$s=1$ and the rest to zero. Then applying this to $\Delta_L$ defines
the above map $\phi$. By
construction, and one may easily check if in doubt, the image of
$\phi$ lies in the fixed subalgebra under the right coaction
$\Delta_R=({\rm id}\otimes\pi)\Delta$ of $\CC[\tilde H]$ on
$\CC[GL_4]$, where $\pi$ is the canonical surjection to
\[ \CC[\tilde H]=\CC[GL_4]/ \la a^3_1=a^3_2=a^4_1=a^4_2=0,\quad a^1_1a^2_2-a^2_1a^1_2=1 \ra \]
and $\Delta$ is the matrix coproduct of $\CC[GL_4]$.

Ultimately we want the same picture for the projective variety
$\hCM$.  In order to do this the usual route in algebraic geometry
is to work with rational functions instead of polynomials in the
homogeneous coordinate algebra and make the quotient by $\C^*$ as
the subalgebra of total degree zero. Rational functions here may
have poles so to be more precise, for any open set $U\subset X$ in a
projective variety, we take the algebra
\[ \CO_X(U)=\{a/b\ |\ a,b\in \CC[\tilde X],\quad a,b\ {\rm same\ degree},\quad b(x)\ne 0\quad \forall x\in U\}\]
where $\CC[\tilde X]$ denotes the homogeneous coordinate algebra
(the coordinate algebra functions of the affine (i.e.
non-projective) version $\tilde X$) and $a,b$ are homogeneous. Doing
this for any open set gives a sheaf of algebras. Of particular
interest are principal open sets of the form $U_f=\{x\in X\ |\
f(x)\ne 0\}$ for any nonzero homogeneous $f$. Then
\[ \CO_X(U_f)=\CC[\tilde X][f^{-1}]_0\]
where we adjoin $f^{-1}$ to the homogeneous coordinate algebra and
$0$ denotes the degree zero part. In the case of $\PGL_4$ we in fact
have a coordinate algebra of regular functions
\[ \CC[\PGL_4] :=\CC[\GL_4]^{\CC[\C^*]}=\CC[\GL_4]_0\]
constructed as the affine algebra analogue of $\GL_4/\C^*$. It is an
affine variety and not projective (yet one could view its coordinate
algebra as $\CO_{\CP^{15}}(U_D)$ where $D$ is the determinant). In
contrast, $\hCM$ is projective and we have to work with sheaves. For
example, $U_t$ is the open set where $t\ne 0$. Then
\[ \CO_{\hCM}(U_t):= \CC[\TCM][t^{-1}]_0. \]
There is a natural inclusion
\[ \CC[\CM] \rightarrow \CO_{\hCM}(U_t)\]
of the coordinate algebra of affine Minkowski space $\CM$
(polynomials in the four coordinate functions $x_1,x_2,x_3,x_4$ on
$\CC^4$ with no further relations) given by
\[ x_1 \mapsto t^{-1}z, \quad x_2 \mapsto t^{-1} \tz, \quad x_3
\mapsto t^{-1}w, \quad x_4 \mapsto t^{-1} \tw. \]
This is the coordinate algebra version of identifying the affine
piece of space-time $\CM$ with the patch of $\hCM$ for which $t \ne
0$.

\section{Twistor Space and the Correspondence Space}
\label{Twistor Space and the Correspondence Space}

Next we give the coordinate picture for twistor space
$T=\CP^3=\F_{1}(\CC^4)$ of lines in $\CC^4$. As a partial flag
variety this is also known to be a homogeneous space. At the
non-projective level we just mean $\tilde T=\CC^4$ with coordinates
$Z=(Z^\mu)$ and the origin deleted. This is of course a homogeneous
space for $\GL_4$ and may be identified as the orbit of the point
$Z^1=1, \, Z^2=Z^3=Z^4=0$: the isotropy subgroup $\tilde K$ consists
of elements such that $a^1_1=1, a^2_1=a^3_1=a^4_1=0$, giving the
identification $\tilde T = \GL_4 /\tilde K$ (again we quotient from
the right).

Again we pass to the coordinate algebra level. At this level we do
not see the deletion of the origin, so we define $\CC[\tilde
T]=\CC[\CC^4]$. We have an isomorphism
\[ \phi : \CC [\tilde T] \rightarrow \CC [\GL_4]^{\CC[\tilde K]},\quad  \phi (Z^\mu) = a^\mu_1 \]
according to a left coaction
\[ \Delta_L (Z^\mu) =  a^\mu_\alpha \tens Z^\alpha. \]
One should view the principal orbit base point as a linear function on $\CC^4$ that
sends $Z^1 = 1$ and the rest to zero:  as before, applying this to
the coaction $\Delta_L$ defines $\phi$ as the dual of the orbit
construction. It is easily verified that the image of this isomorphism is exactly
the subalgebra of $\CC [\GL_4]$ fixed under
\[ \CC[\tilde K]=\CC[\GL_4]/ \la a^1_1=1, a^2_1=a^3_1=a^4_1=0 \ra \]
by the right coaction on $\CC[\GL_4]$ given by projection from the
coproduct $\Delta$ of $\CC[GL_4]$.

Finally we introduce a new space $\F$, the `correspondence space',
as follows. For each point $Z\in T$ we define the associated
`$\alpha$-plane'
\[ \hat Z=\{x\in \hCM\ |\  x \wedge Z = x^{[\mu\nu}Z^{\rho]}=0\}\subset \hCM.\]
The condition on $x$ is independent both of the scale of $x$ and of
$Z$, so constructions may be done `upstairs' in terms of matrices,
but we also have a well-defined  map at the projective level. The
$\alpha$-plane $\hat Z$ contains for example all points in the
quadric of the form $W\wedge Z$ as $W\in \CC^4$ varies. Any multiple
of $Z$ does not contribute, so $\hat Z$ is a 3-dimensional space in
$\TCM$ and hence a $\CP^2$ contained in $\hCM$ (the image of a
two-dimensional subspace of $\CM$ under the conformal
compactification, hence the term `plane').

Explicitly, the condition $x\wedge Z=0$ in our coordinates is:
\begin{eqnarray}
\label{alpha planes}
\tz Z^3 + w Z^4 -t Z^1=0, \quad \tw Z^3 + z Z^4 - t Z^2=0,\\
\nonumber s Z^3+w Z^2-z Z^1=0, \quad s Z^4-\tz Z^2+\tw Z^1=0.
\end{eqnarray}
If $t\ne 0$ one can check that the second pair of equations is
implied by the first (given the quadric relation $\det \, x=0$), so
generically we have two equations for four unknowns as expected.
Moreover, at each point of a plane $\hat Z$ we have in the
Lorentzian case the property that $\nu(A,B)=0$ for any two tangent
vectors to the plane (where $\nu$ is the aforementioned metric on
$\hCM$). One may check that the plane $\hat Z$ defined by $x \wedge
Z=0$ is null if and only if the bivector $\pi=A\wedge B$ defined at
each point of the plane (determined up to scale) is self-dual with
respect to the Hodge $*$-operator. We note that one may also
construct `$\beta$-planes', for which the tangent bivector is
anti-self-dual:  these are instead parameterized by 3-forms in the
role of $Z$.

Conversely, given any point $x\in\hCM$ we define the `line'
\[ \hat x= \{ Z\in T\ |\ x\in \hat Z\}=\{Z\in T\ | \ x^{[\mu\nu}Z^{\rho]}=0\}\subset T.\]
We have seen that we may write $x=a\wedge b$ and indeed $Z=\lambda
a+\mu b$ solves this equation for all $\lambda,\mu\in\CC$. This is a
plane in $\tilde T = \CC^4$ which projects to a $\CP^1$ contained in
$\CP^3$, thus each $\hat x$ is a projective line in twistor space
$T=\CP^3$.

We then define $\F$ to be the set of pairs $(Z,x)$, where $x\in\TCM$
and $Z\in \CP^3=T$ are such that $x\in\hat Z$ (or equivalently
$Z\in\hat x$), i.e. such that $x\wedge Z=0$.  This space fibres
naturally over both space-time and twistor space \textit{via} the
obvious projections
\begin{eqnarray}
\label{double fibration}
\setlength{\unitlength}{1mm}
\begin{picture}(50,15)(0,0)
\put(25,14){\makebox(0,0){$\F$}} \put(10,0){\makebox(0,0){$\CP^3$}}
\put(42,0){\makebox(0,0){$\CM^\#.$}} \put(23,11){\vector(-1,-1){9}}
\put(27,11){\vector(1,-1){9}} \put(15,7){\makebox(0,0){$p$}}
\put(35,7){\makebox(0,0){$q$}}
\end{picture}
\end{eqnarray}
Clearly we have
\[ \hat Z=q (p^{-1}(Z)),\quad \hat x=p(q^{-1}(x)).\]
It is also clear that the defining relation of $\F$ is preserved
under the action of $\GL_4$.

From the Grassmannian point of view, $Z \in \CP^3 = \F_1(\CC^4)$ is
a line in $\CC^4$ and $\hat Z \subset \F_2(\CC^4)$ is the set of
two-planes in $\CC^4$ containing this line. Moreover, $\hat{x}
\subset \F_1(\CC^4)$ consists of all one-dimensional subspaces of
$\CC^4$ contained in $x$ viewed as a two-plane in $\CC^4$.  Then
$\F$ is the partial flag variety
\[ \F_{1,2}(\CC^4)\]
of subspaces $\CC\subset\CC^2\subset\CC^4$.  Here
$x\in\hCM=\F_2(\CC^4)$ is a plane in $\CC^4$ and $Z\in T=\CP^3$ is a
line in $\CC^4$ contained in this plane.  From this point of view it
is known that the homology $H_4(\hCM)$ is two-dimensional and indeed
one of the generators is given by any $\hat Z$ (they are all
homologous and parameterized by $\CP^3$). The other generator is
given by a similar construction of `$\beta$-planes'
\cite{pr:sst,mw:isdtt} with correspondence space $\F_{2,3}(\CC^4)$
and with $\F_{3}(\CC^4)=(\CP^3)^*$. Likewise, the homology
$H_2(\CP^3)$ is one-dimensional and indeed any flag $\hat x$ is a
generator (they are all homologous and parameterized by $\hCM$). For
more details on the geometry of this construction, see
\cite{mw:isdtt}.  More on the algebraic description can be found in
\cite{be:book}.

This $\F$ is known to be a homogeneous space. Moreover, $\tilde\F$
(the non-projective version of $\F$) can be viewed as a quadric in
$(\Lambda^2 \CC^4) \otimes \CC^4$ and hence the orbit under the
action of $\GL_4$ of the point in $\tilde\F$ where $s=1, Z^1 = 1$
and all other coordinates are zero. The isotropy subgroup $\tilde
{R}$ of this point consists of those $a \in \GL_4$ such that $a^2_1
=0$, $a^3_\mu = a^4_\mu = 0$ for $\mu = 1,2$ and $a^1_1=a^2_2=1$. As
one should expect, $\tilde {R}= \tilde {H} \cap \tilde {K}$.

At the level of the coordinate rings, the identification of
$\tilde\F$ with the quadric in $(\Lambda^2 \CC^4) \otimes \CC^4$
gives the definition of $\CC[\tilde \F]$ as the polynomials in the
coordinate functions $x^{\mu \nu},Z^\alpha$, modulo the quadric
relations and the relations (\ref{alpha planes}). That it is an
affine homogeneous space is the isomorphism
\[ \phi: \CC [\tilde \F ] \rightarrow \CC [\GL_4]^{\CC [\tilde {R}]}, \quad \phi (x^{\mu \nu} \otimes Z^{\beta}) = (a^\mu_1 a^\nu_2 - a^\nu_1 a^\mu_2)a^\beta_1, \]
according to the left coaction
\[ \Delta_L (x^{\mu \nu} \otimes Z^{\sigma}) =  a^\mu_\alpha a^\nu_\beta a^\sigma_\gamma \tens (x^{\alpha \beta} \otimes Z^{\gamma}). \]
The image of $\phi$ is the invariant subalgebra under the right coaction of
\[ \CC[\tilde R]=\CC[GL_4]/\la a^2_1 =0, \quad a^1_1a^2_2=1, \quad a^3_1 = a^3_2=a^4_1=a^4_2 = 0 \ra.\]
on $\CC[\GL_4]$ given by projection from the coproduct $\Delta$ of
$\CC[\GL_4]$.

\section{$\SL_4$ and Unitary Versions}
\label{Unitary Versions}

As discussed, the group $\GL_4$ acts on $\TCM = \{ x \in \Lambda^2
\CC^4 ~|~ \det x = 0 \}$ by conjugation, $x \mapsto axa^t$, and
since multiples of the identity act trivially on $\TCM$, this
picture descends to an action of the projective group $\PGL_4$ on
the quotient space $\hCM$.  Our approach accordingly was to
work at the non-projective level in order for the algebraic structure to have an affine form and pass
at the end to the projective spaces $\hCM, T$ and $\F$ as rational functions
of total degree zero.

If one wants to work with these spaces directly as homogeneous
spaces one may do this as well, so that $\hCM=\PGL_4/P\tilde H$, and
so on. From a mathematician's point of view one may equally well
define
\[ \hCM=\F_2(\CC^4)=\GL_4/H,\quad T=\F_1(\CC^4)=\GL_4/K,\quad\F=\F_{1,2}(\CC^4)=\GL_4/R,\]
\[ H=\begin{pmatrix}*&*&*&*\\ * & * & * &*\\ 0 & 0 & * &*\\ 0 & 0 & * &*\end{pmatrix},\quad
K=\begin{pmatrix}*&*&*&*\\ 0 & * & * &*\\ 0 & * & * &*\\ 0 & * & *
&*\end{pmatrix},\quad R=H\cap K=\begin{pmatrix}*&*&*&*\\ 0 & * & *
&*\\ 0 & 0 & * &*\\ 0 & 0 & * &*\end{pmatrix},\]

\medskip \noindent
where the overall $\GL_4$ determinants are non-zero. Here
$H$ is slightly bigger than the subgroup $\tilde H$ we had before.
As homogeneous spaces, $\hCM, T$ and $\F$ carry left actions of
$\GL_4$ which are essentially identical to those given above at both
the $\PGL_4$ and at the non-projective level.

One equally well has
\[ \hCM=\F_2(\CC^4)=\SL_4/H,\quad T=\F_1(\CC^4)=\SL_4/K,\quad\F=\F_{1,2}(\CC^4)=\SL_4/R\]
where $H,K,R$ are as above but now viewed in $\SL_4$, and now
$\hCM,T$ and $\F$ carry canonical left actions of $\SL_4$ similar to
those previously described.

These versions would be the more usual in algebraic geometry but at
the coordinate level one does need to then work with an appropriate
construction to obtain these projective or quasi-projective
varieties. For example, if one simply computes the invariant
functions $\CC[\SL_4]^{\CC[K]}$ etc. as affine varieties, one will
not find enough functions.

As an alternative, we mention a version where we consider all our
spaces in the double fibration as real manifolds, and express this
algebraically in terms of $*$-structures on our algebras. Thus for
example $\CP^3$ is a real 6-dimensional manifold which we construct
by complexifying it to an affine 6-dimensional variety over $\C$,
but we remember its real form by means of a $*$-involution on the
complex algebra. The $*$-algebras in this approach can then in
principle be completed to an operator-algebra setting and the
required quotients made sense of in this context, though we shall
not carry out this last step here.

In this case the most natural choice is
\[ \hCM=\SU_4/H,\quad T=\SU_4/K,\quad \F=\SU_4/R\]
\[H=\rm \S(\U(2)\times \U(2)),\quad K=\rm \S(\U(1)\times \U(3)),\quad R=H\cap K=\rm \S(\U(1)\times \U(1)\times \U(2)),\]
embedded in the obvious diagonal way into $\SU_4$.   As homogeneous
spaces one has canonical actions now of $\SU_4$ from the left on
$\hCM,T$ and $\F$.

For the coordinate algebraic version one expresses $\SU_4$ by
generators $a^\mu_\nu$, the determinant relation and in addition the
$*$-structure
\[ a^\dagger=Sa\]
where $\dagger$ denotes transpose and $*$ on each matrix generator
entry, $(a^\mu_\nu)^\dagger = (a^\nu_\mu)^*$, and $S$ is the Hopf
algebra antipode characterised by $aS(a)=(Sa)a={\rm id}$. This is as
for any compact group or quantum group coordinate algebra. The coordinate algebras
of the subgroups are similarly defined as $*$-Hopf algebras.

There is also a natural $*$-structure on twistor space. To see this
let us write it in the form
\[ T=\CP^3=\{ Q \in M_4(\CC),\quad Q=Q^\dagger, \quad Q^2=Q,\quad \textup{Tr}~Q=1\} \]
in terms of Hermitian-conjugation $\dagger$. Thus $\CP^3$ is the
space of Hermitian rank one projectors on $\C^4$. Such projectors
can be written explicitly in the form $Q^\mu_\nu=Z^\mu \bar Z^\nu$
for some complex vector $Z$ of modulus 1 and determined only up to a
$\U(1)$ normalisation. Thus $\CP^3=S^7/\U(1)$ as a real
6-dimensional manifold. In this description the left action of
$\SU_4$ is given by conjugation in $M_4(\C)$, i.e. by unitary
transformation of the $Z$ and its inverse on $\bar Z$. One can
exhibit the identification with the homogeneous space picture, as
the orbit of the projector ${\rm diag}(1,0,0,0)$ (i.e.
$Z=(1,0,0,0)=Z^*$). The isotropy group of this is the intersection
of $\SU_4$ with $\U(1)\times \U(3)$ as stated.

The coordinate $*$-algebra version is
\[\CC[\CP^3] = \CC[Q^\mu_\nu]/\la   Q^2=Q,\quad
\textup{Tr}~Q=1\ra,\quad Q=Q^\dagger\]
with the last equation now as a definition of the $*$-algebra
structure via $\dagger=(\ )^{*t}$. We can also realise this as the
degree zero subalgebra,
\[ \CC[\CP_3]=\CC[S^7]_0,\quad \CC[S^7]=\C[Z^\mu, Z^\nu{}^*]/\la \sum_\mu Z^\mu{}^* Z^\mu =1\ra,\]
where $Z, Z^*$ are two sets of generators related by the
$*$-involution. The grading is given by $\deg(Z)=1$ and
$\deg(Z^*)=-1$, corresponding to the $\U(1)$ action on $Z$ and its
inverse on $Z^*$. Finally, the left coaction is
\[ \Delta_L   (Z^\mu) =a^\mu_\alpha \tens Z^\alpha,\quad
\Delta_L (Z^\mu{}^*) =Sa^\alpha_\mu \tens Z^\alpha{}^*,\]
as required for a unitary coaction of a Hopf $*$-algebra on a
$*$-algebra, as well as to preserve the relation.

We have similar $*$-algebra versions of $\F$ and $\hCM$ as well. Thus
\[ \CC[\hCM] = \CC[P^\mu_\nu]/\la   P^2=P,\quad
\textup{Tr}~P=2\ra,\quad P=P^\dagger\]
in terms of a rank two projector matrix of generators, while
\[ \CC[\F]=\CC[Q^\mu_\nu,P^\mu_\nu]/\la Q^2=Q,\ P^2=P,\ \Tr~ Q=1,\ \Tr~ P=2,\ PQ=Q=QP\ra\]
\[ Q=Q^\dagger,\ P=P^\dagger.\]
We see that our $*$-algebra approach to flag varieties has a
`quantum logic' form. In physical terms, the fact that $P,Q$ commute
as matrices (or matrices of generators in the coordinate algebras)
means that they may be jointly diagonalised, while $QP=Q$ says that
the 1-eigenvectors of $Q$ are a subset of the 1-eigenvectors of $P$
(equivalently, $PQ=Q$ says that the $0$-eigenvectors of $P$ are a
subset of the $0$-eigenvectors of $Q$). Thus the line which is the
image of $Q$ is contained in the plane which is the image of $P$:
this is of course the defining property of pairs of projectors
$(Q,P) \in \F=\F_{1,2}(\C^4)$. Clearly this approach works for all
flag varieties $\F_{k_1,\cdots k_r}(\C^n)$ of
$k_1<\cdots<k_r$-dimensional planes in $\C^n$ as the $*$-algebra
with $n\times n$ matrices $P_i$ of generators
\[ \C[\F_{k_1,\cdots k_r}]=\C[P_i^\mu{}_\nu;\ i=1,\cdots,r]/\la P_i^2=P_i,\  \Tr ~P_i=k_i,\  P_iP_{i+1}=P_i=P_{i+1}P_i,\ra \]
\[ P_i=P_i^\dagger.\]

In this setting all our algebras are now complex affine varieties
(with $*$-structure) and we can expect to be able to work
algebraically. Thus one may expect for example that
$\CC[\CP^3]=\CC[\SU_4]^{\CC[\textup{S}(\U(1)\times \U(3))]}$
(similarly for other flag varieties) and indeed we may identify the
above generators and relations in the relevant invariant subalgebra
of $\C[\SU_4]$. This is the same approach as was successfully used
for the Hopf fibration construction of $\CP^1=S^2=\SU(2)/\U(1)$ as
$\CC[\SU_2]^{\CC[\U(1)]}$, namely as $\CC[\SL_2]^{\CC[\C^*]}$ with
suitable $*$-algebra structures \cite{Ma:sphere}. Note that one
should not confuse such Hopf algebra (`GIT') quotients with complex
algebraic geometry quotients, which are more complicated to define
and typically quasi-projective. Finally, we observe that in this
approach the tautological bundle of rank $k$ over a flag variety
$\F_k(\C^n)$ appears tautologically as a matrix generator viewed as
a projection $P\in M_n(\C[\F_k])$. The classical picture is that the
flag variety with this tautological bundle is universal for rank $k$
vector bundles.

\subsection{Tautological bundle on $\hCM$ and the instanton as its Grassmann connection}
\label{section taut bundle on cmhash}

Here we conclude with a result that is surely known to some, but
apparently not well-known even at the classical level, and yet drops
out very naturally in our $*$-algebra approach. We show that the
tautological bundle on $\hCM$ restricts in a natural way to
$S^4\subset \hCM$, where it becomes the 1-instanton bundle, and for
which the Grassmann connection associated to the projector is the
1-instanton.

We first explain the Grassmann connection for a projective module
$\CE$ over an algebra $A$. We suppose that $\CE=A^ne$ where $e\in
M_n(A)$ is a projection matrix acting on an $A$-valued row vector.
Thus every element $v\in \CE$ takes the form
\[ v=\tilde v\cdot e=\tilde v_je_j=(\tilde v_k)e_{kj}e_j,\]
where $e_j=e_{j\cdot}\in A^n$ span $\CE$ over $A$ and $\tilde v_i\in
A$. The action of the Grassmann connection is the exterior
derivative on components followed by projection back down to $\CE$:
\[ \nabla v=\nabla(\tilde v.e)= (\extd (\tilde v.e))e=((\extd \tilde v_j)e_{jk}+\tilde v_j \extd e_{jk})e_k=(\extd \tilde v+\tilde v\extd e)e.\]
One readily checks that this is both well-defined and a connection
in the sense that
\[ \nabla(av)=\extd a.v+a\nabla v,\quad \forall a\in A,\quad v\in\CE,\]
and that its curvature operator $F=\nabla^2$ on sections is
\[ F(v)=F(\tilde v.e)=(\tilde v.\extd e.\extd e).e\]

As a warm-up example we compute the Grassmann connection for the
tautological bundle on $A=\C[\CP^1]$. Here $e=Q$, the projection
matrix of coordinates in our $*$-algebraic set-up:
\[ e=Q=\begin{pmatrix} a &  z\\  z^* & 1-a\end{pmatrix};\quad a(1-a)= z z^*,\quad a\in \R,\quad  z\in \C,\]
where $s=a-{1\over 2}$ and $ z=x+\imath y$ describes a usual sphere of radius $1/2$ in Cartesian coordinates $(x,y,s)$. We note that
\[ (1-2a)\extd a=\extd  z. z^*+ z\extd z^*\]
allows to eliminate $\extd a$ in the open patch where $a\ne {1\over 2}$ (i.e. if we delete the north pole of $S^2$). Then
\[ \extd e.\extd e=\begin{pmatrix}\extd a&\extd z \\ \extd  z^* & -\extd a\end{pmatrix}\begin{pmatrix}\extd a&\extd z \\ \extd  z^* & -\extd a\end{pmatrix}=\extd z\extd z^*\begin{pmatrix} 1&-{2 z\over 1-2a} \\ -{2 z^*\over 1-2a} & -1\end{pmatrix}={\extd z\extd z^*\over 1-2a}(1-2e)\]
and hence
\[ F(\tilde v.e)=-{\extd z\extd z^*\over 1-2a}\tilde v.e.\]
In other words, $F$ acts as a multiple of the identity operator on
$\CE=A^2.e$ and this multiple has the standard form for the charge 1
monopole connection if one converts to usual Cartesian coordinates.
We conclude that the Grassmann connection for the tautological
bundle on $\CP^1$ is the standard 1-monopole. This is surely
well-known. The $q$-deformed version of this statement can be found
in \cite{HajMa:pro} provided one identifies the projector introduced
there as the defining projection matrix of generators for
$\C_q[\CP^1]=\C_q[\SL_2]^{\C[t,t^{-1}]}$ as a $*$-algebra in the
$q$-version of the above picture (the projector there obeys
$\Tr_q(e)=1$ where we use the $q$-trace).  Note that if one looks
for any algebra $A$ containing potentially non-commuting elements
$a,z,z^*$ and a projection $e$ of the form above with $\Tr(e)=1$,
one immediately finds that these elements commute and obey the
sphere relation as above. If one performs the same exercise with the
$q$-trace, one finds exactly the four relations of the standard
$q$-sphere as a $*$-algebra.

Next, we look in detail at $A=\C[\hCM]$ in our projector $*$-algebra
picture.  This has a $4\times 4$ matrix of generators which we write
in block form
\[ P=\begin{pmatrix} A & B \\ B^\dagger & D\end{pmatrix},\quad \Tr A+\Tr D=2,\quad A^\dagger=A,\quad D^\dagger=D\]
\begin{equation}\label{hcma}A(1-A)=BB^\dagger,\quad D(1-D)=B^\dagger B\end{equation}
\begin{equation}\label{hcmb} (A-{1\over 2})B+B(D-{1\over 2})=0,\end{equation}
where we have written out the requirement that $P$ be a Hermitian
$A$-valued projection without making any assumptions on the
$*$-algebra $A$ (so that these formulae also apply to any
noncommutative version of $\C[\hCM]$ in our approach).

To proceed further, it is useful to write
\[ A=a+  \alpha\cdot\sigma,\quad B=t+\imath x\cdot \sigma,\quad B^\dagger=t^*-\imath x^*\cdot\sigma,\quad D=1-a+ \delta\cdot\sigma\]
in terms of usual Pauli matrices $\sigma_1,\sigma_2,\sigma_3$. We
recall that these are traceless and Hermitian, so $a,\alpha,\delta$
are self-adjoint, whilst $t,x_i$, $i=1,2,3$ are not necessarily so
and are subject to (\ref{hcma})-(\ref{hcmb}).

\begin{prop} \label{hcmproj} The commutative $*$-algebra $\C[\hCM]$ is defined by the above generators $a=a^*, \alpha=\alpha^*, \delta=\delta^*,t,t^*,  x,  x^*$ and the relations
\[ tt^*+xx^*=a(1-a)-\alpha.\alpha,\quad (1-2a)(\alpha-\delta)=2\imath(t^*x-tx^*),\quad (1-2a)(\alpha+\delta)=2\imath x\times x^*\]
\[  \alpha.\alpha=\delta.\delta,\quad (\alpha+\delta).x=0,\quad (\alpha+\delta)t=(\alpha-\delta)\times x\]
\end{prop}
\begin{proof}
This is a direct computation of (\ref{hcma})-(\ref{hcmb}) under the
assumption that the generators commute. Writing our matrices in the
form above, equations (\ref{hcma}) become
\[ a(1-a)-\alpha\cdot\alpha+(1-2a)\alpha\cdot\sigma=tt^*+x\cdot\sigma x^*\cdot\sigma+\imath (xt^*-tx^*)\cdot\sigma\]
\[ (1-a)a-\delta\cdot\delta-(1-2a)\delta\cdot\sigma=t^*t+x^*\cdot \sigma x\cdot\sigma+\imath(t^*x-x^*t)\cdot\sigma.\]
Taking the sum and difference of these equations and in each case
the parts proportional to $1$ (which is the same on both right hand
sides) and the parts proportional to $\sigma$ (where the difference
of the right hand sides is proportional to $x\times x^*$) gives four
of the stated equations (all except those involving terms
$(\alpha+\delta)\cdot x$ and $(\alpha+\delta)t$). We employ the key
identity
\[ \sigma_i\sigma_j=\delta_{ij}+\imath\eps_{ijk}\sigma_k,\]
where $\eps$ with $\eps_{123}=1$ is the totally antisymmetric tensor used in the definition of the vector cross product. Meanwhile, (\ref{hcmb}) becomes
\[ \imath \alpha.\sigma x\cdot\sigma+\alpha t\cdot\sigma+\imath x\cdot\sigma \delta\cdot\sigma+t\delta\cdot\sigma=0\]
after cancellations, and this supplies the remaining two relations
using our key identity.   \end{proof}

We see that in the open set where $a\ne {1\over 2}$ we have
$\alpha,\delta$ fully determined by the second and third relations,
so the only free variables are the complex generators $t,\vec x$,
with $a$ determined from the first equation. The complex affine
variety generated by the independent variables $x$, $x^*$, $t$,
$t^*$ modulo the first three equations is reducible; the second
`auxiliary' line of equations makes $\C_F[\hCM]$ reducible (we
conjecture this).

\begin{prop}
There is a natural $*$-algebra quotient $\C[S^4]$ of $\C[\hCM]$
defined by the additional relations $x^*=x$, $t^*=t$ and
$\alpha=\delta=0$.  The tautological projector of $\C[\hCM]$ becomes
\[ e=\begin{pmatrix} a & t+\imath x\cdot \sigma\\ t-\imath x\cdot\sigma & 1-a\end{pmatrix}\in M_2(\C[S^4]).\]
The Grassmann connection on the projective module $\CE=\C[S^4]^4e$ is the 1-instanton with local form
\[ (F\wedge F)(\tilde v.e)=-4!{\extd t\extd^3 x\over 1-2a}\tilde v.e\]
\end{prop}
\begin{proof}
All relations in Proposition~\ref{hcmproj} are trivially satisfied
in the quotient except $tt^*+xx^*=a(1-a)$, which is that of a
4-sphere of radius $1\over 2$ in usual Cartesian coordinates
$(t,x,s)$ if we set $s=a-{1\over 2}$.  The image $e$ of the
projector exhibits $S^4\subset\hCM$ as a projective variety in our
$*$-algebra projector approach.  We interpret this as providing a
projective module over $\C[S^4]$, the pull-back of the tautological
bundle on $\hCM$ from a geometrical point of view. To compute the
curvature of its Grassmann connection we first note that
\[ \extd x\cdot \sigma \extd x\cdot\sigma=\imath (\extd x\times \extd x)\cdot\sigma,\quad (\extd x\times\extd x)\cdot(\extd x\times\extd x)=0,\]
\[ (\extd x\times\extd x)\times(\extd x\times\extd x)=\extd x\times (\extd x\times\extd x)=0,\quad  \extd x\cdot (\extd x\times\extd x)=3!\extd^3x,\]
since 1-forms anticommute and since any four products of the $\extd
x_i$ vanish.  Now we have
\[ \extd e\extd e=\begin{pmatrix}\extd a& \extd t+\imath\extd x\cdot \sigma\\
\extd t-\imath\extd x\cdot\sigma& -\extd a\end{pmatrix}\begin{pmatrix}\extd a& \extd t+\imath\extd x\cdot \sigma\\
\extd t-\imath\extd x\cdot\sigma& -\extd a\end{pmatrix}\]
\[=\begin{pmatrix}-2\imath\extd t\extd x\cdot\sigma+\imath\extd x\times\extd x\cdot\sigma& 2\extd a\extd t+2\imath\extd a \extd x\cdot \sigma\\
-2\extd a\extd t+2\imath\extd a \extd x\cdot \sigma & 2\imath\extd
t\extd x\cdot\sigma+\imath\extd x\times\extd
x\cdot\sigma\end{pmatrix}\]
and we square this matrix to find that
\[ (\extd e)^4=\begin{pmatrix}1 & -{2\over 1-2a}(t+\imath x\cdot\sigma)\\ -{2\over 1-2a}(t-\imath x\cdot\sigma) & -1\end{pmatrix}4!\extd t\extd^3x={1-2e\over 1-2a}4!\extd t\extd^3x\]
after substantial computation. For example, since $(\extd a)^2=0$,
the 1-1 entry is
\[ -(\extd x\times\extd x-2\extd t\extd x).\sigma (\extd x\times\extd x-2\extd t\extd x)\cdot\sigma=2\extd t\extd x\cdot(\extd x\times\extd x)+2(\extd x\times\extd x)\cdot\extd t\extd x=4!\extd t\extd^3x\]
where only the cross-terms contribute on account of the second
observation above and the fact that $(\extd t)^2=0$. For the 1-2
entry we have similarly that
\begin{eqnarray*}&& 2\imath\extd a(\extd x\times\extd x-2\extd t\extd x)\cdot\sigma(\extd t+\imath\extd x\cdot\sigma)
+2\imath\extd a(\extd t+\imath\extd x\cdot\sigma)(\extd x\times\extd x+2\extd t\extd x)\cdot\sigma\\
 &=&12\imath\extd a\extd t(\extd x\times\extd x)\cdot\sigma-4\extd a(\extd x\times\extd x)\cdot\sigma\extd x\cdot\sigma\\
 &=&-{24\over 1-2a}2\imath x\cdot \sigma\extd t\extd^3x-{8\over 1-2a}t 3!\extd t\extd^3x=-{2\over 1-2a}(t+\imath x\cdot \sigma)4!\extd t\extd^3 x,\end{eqnarray*}
where at the end we substitute
\[ \extd a={2( t\extd t+ x\cdot\extd x) \over 1-2a}\]
and note that
\[ x\cdot\extd x(\extd x\times\extd x)=x_i\extd x_i\eps_{jkm}\extd x_j\extd x_k=2x_m\extd^3x,\]
since in the sum over $i$ only $i=m$ can contribute for a nonzero
3-form. The $2$-$2$ and $2$-$1$ entries are analogous and left to
the reader. We conclude that $(\extd e)^4e$ acts on $\C[S^4]^4$ from
the right as a multiple of the identity as stated. \end{proof}

One may check that $F=\extd e \extd e.e$ is anti-self-dual with
respect to the usual Euclidean Hodge $*$-operator. Note also that
the off-diagonal corners of $e$ are precisely a general quaternion
$q=t+\imath x\cdot \sigma$ and its conjugate, which relates our
approach to the more conventional point of view on the 1-instanton.
However, that is not our starting point as we come from $\hCM$,
where the top right corner is a general $2\times 2$ matrix $B$ and
the bottom left corner its adjoint.

If instead we let $B$ be an arbitrary Hermitian matrix in the form
\[ B=t+x\cdot\sigma \]
(i.e. replace $\imath x$ above by $x$ and let $t^*=t,x^*=x$) then the quotient $t^*=t, x^*=x, \alpha=-\delta=2tx/(1-2a)$ gives us
\[ s^2+t^2+(1+{t^2\over s^2})x^2={1\over 4}\]
when $s=a-{1\over 2}\ne 0$ and $tx_i=0, t^2+x^2+\alpha^2={1\over 4}$
when $s=0$. We can also approach this case directly from
(\ref{hcma})-(\ref{hcmb}). We have to find Hermitian $A,D$, or
equivalently $S=A-{1\over 2}$, $T=D-{1\over 2}$ with $\Tr(S+T)=0$
and $S^2=T^2={1\over 4}-B^2$. Since $B$ is Hermitian it has real
eigenvalues and indeed after conjugation we can rotate $x$ to $|x|$
times a vector in the 3-direction, i.e. $B$ has eigenvalues $t\pm
|x|$. It follows that square roots $S,T$ exist precisely when
\begin{equation}\label{diamond} |t|+|x|\le {1\over 2}\end{equation}
and are diagonal in the same basis as was $B$, hence they
necessarily commute with $B$.  In this case (\ref{hcmb}) becomes
that $(S+T).B=0$. If $B$ has two nonzero eigenvalues then $S+T=0$.
If $B$ has one nonzero eigenvalue then $S+T$ has a zero eigenvalue,
but the trace condition then again implies $S+T=0$. If $B=0$ our
equations reduce to those for two self-adjoint $2\times 2$
projectors $A,D$ with traces summing to two. In summary, if $B\ne 0$
there exists a projector of the form required if and only if $(t,x)$
lies in the diamond region (\ref{diamond}), with $S=T$ and a
fourfold choice (the choice of root for each eigenvalue) of $S$ in
the interior. These observations about the moduli of projectors with
$B$ Hermitian means that the corresponding quotient of $\C[\hCM]$ is
a fourfold cover of a diamond region in affine Minkowski space-time
(viewed as the space of $2\times 2$ Hermitian matrices). The diamond
is conformally equivalent to a compactification of all of usual
Minkowski space (the Penrose diagram for Minkowski space), while its
fourfold covering reminds us of the Penrose diagram for a
black-white hole pair. It is the analogue of the disk that one
obtains by projecting $S^4$ onto its first two coordinates. One may
in principle compute the connection associated to the pull-back of
the tautological bundle to this region as well as the 4-dimensional
object of which it is a projection. The most natural version of this
is to slightly change the problem to two $2\times 2$ Hermitian
matrices $S,B$ with $S^2+B^2={1\over 4}$ (a `matrix circle'),  a
variety which will be described elsewhere.

Note also that in both cases  $D=1-A$ and if we suppose this at the
outset our equations including (\ref{hcmb}) and (\ref{hcma})
simplify to
\begin{equation}\label{hcmd}[A,B]=[B,B^\dagger]=0,\quad A(1-A)=BB^\dagger,\quad A=A^\dagger.\end{equation}
In fact this is the same calculation as for any potentially
noncommutative $\CP^1$ which (if we use the usual trace) is forced
to be commutative as mentioned above.

Finally, returning to the general case of $\C[\hCM]$, we have
emphasised `Cartesian coordinates' with different signatures. From a
twistor point of view it is more natural to work with the four
matrix entries of $B$ as the natural twistor coordinates. This will
also be key when we quantise. Thus equivalently to Proposition~3.1
we write
\begin{equation}
\label{twistorB} B=\begin{pmatrix} z & \tilde w \\ w &  \tilde z
\end{pmatrix},\ A=\begin{pmatrix} a+\alpha_3 & \alpha \\ \alpha^* &
a-\alpha_3 \end{pmatrix},\ D=\begin{pmatrix} 1-a +\delta_3& \delta
\\ \delta^* & 1-a-\delta_3 \end{pmatrix},\end{equation}
where $a=a^*,\alpha_3=\alpha_3^*,\delta_3=\delta_3^*$ as before but
all our other notations are different. In particular,
$\alpha,\alpha^*,\delta,\delta^*,z,z^*,w,w^*,\tilde z,\tilde
z^*,\tilde w,\tilde w^*$ are now complex generators.

\begin{corollary} The relations of $\C[\hCM]$ in these new notations appear as
\[ zz^*+ww^*+\tilde z\tilde z^*+\tilde w\tilde w^*=2(a(1-a)-\alpha\alpha^*-\alpha_3^2)\]
\[ (1-2a)\alpha=zw^*+\tilde w\tilde z^*,\quad (1-2a)\delta=-z^*\tilde w-\tilde z w^*\]
\[ (1-2a)(\alpha_3+\delta_3)=\tilde w\tilde w^*-w w^*,\quad (1-2a)(\alpha_3-\delta_3)=zz^*-\tilde z\tilde z^*\]
and the auxiliary relations
\[ \alpha\alpha^*+\alpha_3^2=\delta\delta^*+\delta_3^2\]
\[  (\alpha_3+\delta_3)\begin{pmatrix}z\\ \tilde z\end{pmatrix}=\begin{pmatrix}-\alpha & - \delta^*\\ \delta & \alpha^*\end{pmatrix}\begin{pmatrix}w\\ \tilde w\end{pmatrix},\quad (\alpha_3-\delta_3)\begin{pmatrix}w\\ \tilde w\end{pmatrix}=\begin{pmatrix}-\alpha^* & -\delta^*\\ \delta & \alpha\end{pmatrix}\begin{pmatrix}z\\ \tilde z\end{pmatrix}.\]
Moreover, $S^4\subset\hCM$ appears as the $*$-algebra quotient $\C[S^4]$ defined by $w^*=-\tilde w$, $z^*=\tilde z$ and $\alpha=\delta=0$.
\end{corollary}
\begin{proof} It is actually easier to recompute these, but of course this is just a change of
generators from the equations in Proposition~3.1. \end{proof}

Note that these affine $*$-algebra coordinates are more similar in
spirit but not the same as those for $\hCM$ as a projective quadric
in Section 1.

\subsection{Twistor space $\CP^3$ in the $*$-algebra approach}

For completeness, we also describe $\CP^3$ more explicitly in our
affine $*$-algebra approach. As a warm up we start with $\CP^2$
since $\CP^1$ is already covered above. Thus $\C[\CP^2]$  has a
trace 1 matrix of generators
\[ Q=\begin{pmatrix} a& x & y\\ x^* & b & z\\ y^* & z^* & c\end{pmatrix},\quad a+b+c=1\]
with $a,b,c$ self-adjoint.

 \begin{prop} $\C[\CP^2]$ is the algebra with the above matrix of generators with $a+b+c=1$ and the projector relations
 \[ x^*x=ab,\quad y^*y=ac,\quad z^*z=bc\]
 \[ cx=yz^*,\quad by=xz,\quad az=x^*y.\]
\end{prop}
\begin{proof} First of all, the `projector relations' $Q^2=Q$  come out as
the second line of relations stated and the relations
\[ a(1-a)=X+Y,\quad b(1-b)=X+Z,\quad c(1-c)=Y+Z\]
where we use the shorthand $X=x^*x, Y=y^*y, Z=z^*z$.  We subtract these
from each other to obtain
\[ X-Z=(a-c)b,\quad Y-Z=(a-b)c,\quad X-Y=(b-c)a\]
(in fact there are only two independent ones here).  Combining with
the original relations allows to solve for $X,Y,Z$ as stated.
\end{proof}

Clearly, if $a\ne 0$ (say), i.e. if we look at $\C[\CP^2][a^{-1}]$, we can regard $x,y$ (and their adjoints) and $a,a^{-1},b,c$ as generators with the relations
\begin{equation}\label{radCP2} x^*x=ab,\quad y^*y=ac,\quad a+b+c=1\end{equation}
and all the other relations become empty. Thus $az=x^*y$ is simply
viewed as a definition of $z$ and one may check for example that
$z^*za^2=y^*xx^*y=XY=bca^2$, as needed. Likewise, for example,
$ayz^*=yy^*x=Yx=acx$ as required. We can further regard
(\ref{radCP2}) as defining $b,c$, so the localisation viewed in this
way is a punctured $S^4$ with complex generators $x,y$, real
invertible generator $a$ and the relations
\[ x^*x+y^*y=a(1-a),\]
conforming to our expectations for $\CP^2$ as a complex 2-manifold.

We can also consider setting $a,b,c$ to be real numbers with
$a+b+c=1$ and $b,c>0, b+c<1$. The inequalities here are equivalent
to $ab,ac>0$ and $a\ne 0$ (with $a>0$ necessarily following since if
$a<0$ we would need $b,c<0$ and hence $a+b+c<0$, which is not
allowed). We then have
\[ \C[\CP^2]|_{{b,c>0\atop b+c<1}}=\C[S^1\times S^1],\]
so the passage to this quotient algebra is geometrically an
inclusion $S^1\times S^1\subset \CP^2$ with (\ref{radCP2}) defining
the two circles (recall that $x,y$ are complex generators). As the
parameters vary the circles vary in size so the general case with
$a$ inverted can be viewed in that sense as an inclusion
$\C^*\times\C^*\subset \CP^2$. This holds classically as an open
dense subset (since $\CP^2$ is a toric variety). We have the same
situation for $\C[\CP^1]$ where there is only one relation
$x^*x=a(1-a)$, i.e. circles $S^1\subset \CP^1$ of different size as
$0<a<1$. They are the circles of constant latitude and as $a$ varies
in this range they map out $\C^*$ (viewed as $S^2$ with the north
and south pole removed).

We now find similar results for $\CP^3$  (the
general $\CP^n$ case is analogous). We now have a matrix of generators
\[ Q=\begin{pmatrix} a & x & y & z\\ x^* & b & w & v \\ y^* & w^* & c & u \\ z^* & v^* & u^* & d\end{pmatrix},\quad a^*=a,\ b^*=b,\ c^*=c,\ d^*=d,\ a+b+c+d=1\]
and make free use of the shorthand notation
\[ X=x^*x,\quad Y=y^*y,\quad Z=z^*z,\quad U=u^*u,\quad V=v^*u\quad, W=w^*w.\]

\begin{prop} $\C[\CP^3]$ is the commutative $*$-algebra with generators $Q$ of the form above with $a+b+c+d=1$ and projector relations
\[  a(1-a)=X+Y+Z,\quad X-U=ab-cd,\quad Y-V=ac-bd,\quad Z-W=ad-bc,\]
\[ au=y^*z,\quad av=x^*z,\quad  aw=x^*y,\quad  bu=w^*v,\quad cv=wu,\quad dw=vu^*,\]
\[ cx=yw^*,\quad by=xw,\quad bz=xv,\quad dx=zv^*,\quad dy=zu^*,\quad cz=yu.\]
\end{prop}

\begin{proof} We first write out the relations $P^2=P$ as
\begin{equation}\label{pra} a(1-a)=X+Y+Z,\quad b(1-b)=X+V+W,\end{equation}
\begin{equation} \label{prb} c(1-c)=Y+U+W,\quad d(1-d)=Z+U+V,\end{equation}
\begin{equation}\label{prc} yw^*+zv^*=x(c+d),\quad xw+zu^*=y(b+d),\quad xv+yu=z(b+c),\end{equation}
\begin{equation} \label{prd} y^*z+w^*v=u(a+b),\quad x^*z+wu=v(a+c),\quad x^*y+vu^*=w(a+d).\end{equation}
We add and subtract several combinations of (\ref{pra})-(\ref{prb})
to obtain the equivalent four equations stated in the proposition.
For example, subtracting (\ref{pra}) gives $(Y-V)+(Z-W)=(c+d)(a-b)$
while subtracting (\ref{prb}) gives $(Y-V)-(Z-W)=(c-d)(a+b)$ and
combining these gives the $Y-V$ and $Z-W$ relations stated.
Similarly, for the $X-U$ relation.  We can also write our three
equations as
\begin{equation}\label{prab}  (a+c)(a+d)+X-U=(a+b)(a+d)+Y-V=(a+b)(a+c)+Z-W=a\end{equation}
using $a+b+c+d=1$.

Next, we compute (\ref{prd}) assuming (\ref{prc}), for example
\begin{eqnarray*} (a+b)(a+c)u&=&(a+c)(y^*z+w^*v)=(a+c)y^*z+w^*(x^*z+wu)\\
&=&(a+c)y^*z+Wu+(y^*(b+d)-uz^*)z=y^*z+(W-Z)u\end{eqnarray*}
which, using (\ref{prab}), becomes $au=y^*z$. We similarly obtain
$av=x^*z, aw=x^*y$. Given these relations, clearly (\ref{prd}) is
equivalent to the next three stated equations, which completes the
first six equations of this type. Similarly for the remaining six.
\end{proof}

\begin{lemma} In $\C[\CP^3]$ we have
\[ (X-ab)(Y-(ac-bd))=0,\quad (X-ab)(Z-(ad-bc))=0\]
\[(Y-ac)(X-(ab-cd))=0,\quad (Y-ac)(Z-(ad-bc))=0\]
\[(Z-ad)(X-(ab-cd))=0,\quad (Z-ad)(Y-(ac-bd))=0\]
\[ (X-ab)(X-b(1-a))=0,\quad (Y-ac)(Y-c(1-a))=0,\quad (Z-ad)(Z-d(1-a))=0\]
\[ (X-ab)(X-a(1-b))=0,\quad (Y-ac)(Y-a(1-c))=0,\quad (Z-ad)(Z-a(1-d))=0\]
\end{lemma}
\begin{proof}
For example, $adu=dy^*z=uz^*z=uZ$. In this way one has
\[  (X-ab)v=(X-ab)w=(Y-ac)u=(Y-ac)w=(Z-ad)u=(Z-ad)v=0.\]
Multiplying by $u^*,v^*,w^*$ and replacing $U,V,W$ using
Proposition~3.5 give the first two lines of relations.  Next,
$by=xw$ in Proposition~3.5 implies $b^2Y=XW$ and similarly for $bz$
gives $b^2(Y+Z)=X(V+W)$. We then use (\ref{pra}) to obtain
$X^2-bX+b^2a(1-a)$ which factorises to one of the quadratic
equations stated. Similarly, the equations $a^2V=XZ$, $a^2W=XY$
imply $a^2(V+W)=X(Y+Z)$, which yields the other quadratic equation
for $X$. Similarly for the other quadratic equations. \end{proof}

\begin{lemma} If we consider the trace one projection $Q$ as a numerical Hermitian matrix of the form above, then
\[ X=ab,\quad Y=ac,\quad Z=ad\]
necessarily holds.
\end{lemma}
\begin{proof}
We use the preceding lemma but regarded as for real numbers
(equivalently one can assume that our algebra has no zero divisors).
Suppose without loss of generality that $X\ne ab$.  Then by the
lemma, $V=W=0$ or $Y=ac-bd$, $Z=ad-bc$. We also have $v=w=0$ and
hence from Proposition~3.5 that $x^*z=x^*y=0$. We can also deduce
from the quadratic equations of $X$ that $a=b$ and $X=a(1-a)$ or
$U=a(1-2a)+cd$.  We distinguish two cases: (i) $x=0$ in which case
$a=b=0$ (since $X=a(1-a)\ne ab=a^2$) and (ii) $x\ne 0$, $y=z=0$, in
which case $a=b\ne 0,1$, $c+d=0$. In either case since $b\ne 1$ we
have $Y\ne ac$ and $Z\ne ad$, hence $X=ab-cd$ or $U=0$, $u=0$ and
hence $y^*z=0$ while $c,d\ne 0$. This means that at most one of
$x,y,z$ is non-zero. We can now go through all of the subcases and
find a contradiction in every case.  Similar arguments prove that
$Y=ac, Z=ad$. \end{proof}

Let us denote by $\C^-[\CP^3]$ the quotient of $\C[\CP^3]$ by the
relations in the lemma. We call it the `regular form' of the
coordinate algebra for $\CP^3$ in our $*$-algebraic approach and
will work with it henceforth.  The lemma means that there is no
discernible difference (if any) between the $*$-algebras
$\C^-[\CP^3]$ and $\C[\CP^3]$ in the sense that if we were looking
at $\CP^3$ as a set of projector matrices and the above variables as
real or complex numbers, we would not see any distinction.  (As long
as the relevant intersections are transverse the same would then be
true in the algebras also, but it is beyond our scope to prove this
here.)  Moreover, if either $x,y,z$ or $u,v,w$ are made invertible
then one can show that the solution in the lemma indeed holds and
does not need to be imposed, i.e. $\C[\CP^3]$ and $\C^-[\CP^3]$ have
the same localisations in this respect.

\begin{prop} $\C^-[\CP^3]$ can be viewed as having generators $a,b,c,x,y,z$ with $a+b+c+d=1$ and the relations
\[ X=ab,\quad Y=ac,\quad Z=ad\]
as well as auxiliary generators $u,v,w$ and auxiliary relations
\[U=cd,\quad V=bd,\quad W=bc,\]
\[ au=y^*z,\quad av=x^*z\quad  aw=x^*y,\quad bu=w^*v,\quad cv=wu,\quad dw=vu^*,\]
\[ cx=yw^*,\quad by=xw,\quad bz=xv,\quad dx=zv^*,\quad dy=zu^*,\quad cz=yu.\]
If $a\ne 0$ these auxiliary variables and equations are redundant.
\end{prop}
\begin{proof} If $a\ne 0$ (i.e. if we work in the algebra with $a^{-1}$ adjoined) we regard three of
the auxiliary equations as a definition of $u,v,w$. We then verify
that the other equations hold automatically. The first line is clear
since these equations times $a^2$ were solved in the lemma above.
For example, from the next line we have
$a^2w^*v=y^*xx^*z=Xy^*z=aby^*z=a^2bu$, as required. Similarly
$a(yw^*+zv^*)=yy^*x+zz^*x=(Y+Z)x=a(c+d)x$ as required.
\end{proof}

{}From this we see that the `patch' given by inverting $a$ is
described by just three independent complex variables $x,y,z$ and
one invertible real variable $a$ with the single relation
\[ x^*x+y^*y+z^*z=a(1-a)\]
(the relations stated can be viewed as a definition of $b,c,d$ but
we still need $a+b+c+d=1$), in other words a punctured $S^6$ where
the point $x=y=z=a=0$ is deleted. This conforms to our expectations
for $\CP^3$ as a complex 3-manifold or real 6-manifold. Of course,
our original projector system was symmetric and we could have
equally well analysed and presented our algebra in a form adapted to
one of $b,c,d\ne 0$.

Finally, we also see that if we set $a,b,c,d$ to actual real values
with $a+b+c+d=1$ and $b,c,d>0, b+c+d<1$  (the inequalities here are
equivalent to $a\ne 0$ and $ab,ac,ad>0$) then
\[ \C^-[\CP^3]|_{{b,c,d>0\atop b+c+d<1}}  =\C[S^1\times S^1\times S^1]\]
for three circles $x^*x=ab$, $y^*y=ac$, $z^*z=ad$ and no further
relations.  This is the analogue in our $*$-algebra approach of
inclusions  $S^1\times S^1\times S^1\subset \CP^3$ and as the
circles vary in radius we have part of the fact in the usual picture
that $\CP^3$ is a toric variety (namely that
$\C^*\times\C^*\times\C^*\subset \CP^3$ is open dense).

We now relate this description of twistor space to the space-time
algebras in the previous section.  In particular, we note that
classically there is a fibration of twistor space over the Euclidean
four-sphere, $\CP^3 \rightarrow S^4$, whose fibre is $\CP^1$ (see
for example \cite{ww:tgft}).  This fibration arises through the
observation that each $\alpha$-plane in $\hCM$ intersects $S^4$ at a
unique point (essentially because there are no null lines in
Euclidean signature). The double fibration (\ref{double fibration})
thus collapses to a single fibration $\CP^3 \rightarrow S^4$ (making
the twistor theory of the real space-time $S^4$ much easier to study
than that of its complex counterpart). To see this we make use of
the following nondegenerate antilinear involution on $\C^4$,
\[J(Z)=J(Z^1,Z^2,Z^3,Z^4) := (-\bar Z^2, \bar Z^1, -\bar
Z^4, \bar Z^3).\]
Once again we recall that points of twistor space are
one-dimensional subspaces of $\C^4$, whereas points of $\hCM$ are
two-dimensional subspaces.  Of course, given a $1$d subspace
(spanned by $Z \in \C^4$), there are many $2$d subspaces in which it
lies, and these constitute exactly the set $\hat Z = \CP^2$.
However, the involution $J$ serves to pick out a unique such $2$d
subspace, the one spanned by $Z$ and $J(Z)$.

Now recall our `quantum logic' interpretation of the correspondence
space $\F$, as pairs of projectors $(Q,P)$ on $\C^4$ with $Q$ of
rank one and $P$ of rank two such that $QP=Q=PQ$.  Then since we
have
\[ \CP^3 = S^7/\U(1), \qquad  S^7 = \{ Z^\mu, \bar Z^\nu ~|~ \sum Z^\mu \bar Z^\mu = 1\}, \]
the involution $J$ extends to one on $\CP^3$, given by
\[J(Q^{\mu}{}_\nu) = J(Z^\mu \bar Z^\nu) = J(\bar
Z^\nu)J(Z^\mu).\]
At the level of the coordinate algebra $\C^-[\CP^3]$ we have the
following interpretation.

\begin{lemma} There is an antilinear involution $J: \C^-[\CP^3]
\rightarrow \C^-[\CP^3]$, given in the notation of
Proposition~3.8 by
\[ J(a) = b, \quad J(b)=a, \quad J(c) = 1-(a+b+c),\]
\[ J(x) = -x, \quad J(y)=v^*, \quad J(z) = -w^*,\quad J(u)=-u, \quad J(v) = y^*, \quad J(w) = -z^*,\]
\end{lemma}

\proof  This is by direct computation, noting that if we write
\[ Q = \begin{pmatrix} \bar Z^1 Z^1 & \bar Z^1 Z^2 & \bar Z^1 Z^3 & \bar Z^1
Z^4 \\ \bar Z^2 Z^1 & \bar Z^2 Z^2 & \bar Z^2 Z^3 & \bar Z^2 Z^4 \\
\bar Z^3 Z^1 & \bar Z^3 Z^2 & \bar Z^3 Z^3 & \bar Z^3 Z^4 \\ \bar
Z^4 Z^1 & \bar Z^4 Z^2 & \bar Z^4 Z^3 & \bar Z^4 Z^4 \end{pmatrix},
\qquad \Tr \, Q = 1,\]
we see that
\[ J(Q) = \begin{pmatrix} \bar Z^2 Z^2 & -\bar Z^1 Z^2 & \bar Z^4 Z^2 & -\bar
Z^3 Z^2 \\ -\bar Z^2 Z^1 & \bar Z^1 Z^1 & -\bar Z^4 Z^1 & \bar Z^3 Z^1 \\
\bar Z^2 Z^4 & -\bar Z^1 Z^4 & \bar Z^4 Z^4 & -\bar Z^3 Z^4 \\ -\bar
Z^2 Z^3 & \bar Z^1 Z^3 & -\bar Z^4 Z^3 & \bar Z^3 Z^3 \end{pmatrix},
\qquad \Tr \, J(Q) = \Tr \, Q = 1,\]
and the result follows by comparing with the notation of
Proposition~3.8. In particular we see that $J(X) = X$ and $J(U)=U$.
The relations of Proposition~3.5 indicate that $J$ extends to the
full algebra as an antialgebra map (indeed this needs to be the case
for $J$ to be well-defined), since then we have
\[ J(au) = J(u)J(a) = -ub = -bu = -w^*v = J(z) J(y^*)
= J(y^*z),\] similarly for the remaining relations.
\endproof

We remark that since the algebra is commutative here, we may treat
$J$ as an algebra map rather than as an antialgebra map as required
in the notion of an antilinear involution.  This will no longer be
the case when we come to quantise, when it is the notion of
antilinear involution that will survive.

Now $J$ extends further to an involution on $\hCM$: given $P \in
\hCM$ we write $P=Q+Q'$ for $Q,Q' \in \CP^3$ and define
\[ J(P) := J(Q) + J(Q').\]
Indeed, we note that the $1$d subspaces defined by a pair of rank
one projectors $Q_1$, $Q_2$ are distinct if and only if
$Q_1Q_2=Q_2Q_1=0$, and it is easily checked that this is equivalent
to the condition that $J(Q_1)J(Q_2)=J(Q_2)J(Q_1)=0$ (computed either
by direct calculation with matrices or by working with vectors
$Z_1$, $Z_2 \in \C^4$ which define $Q_1$, $Q_2$ up to scale and
using the inner product on $\C^4$ induced by $J$). Thus if $P$ is a
rank two projector with $P = Q_1+Q_1' = Q_2+Q_2'$, then if $Q_1$,
$Q_2$ are distinct, so are $J(Q_1)$, $J(Q_2)$. Elementary linear
algebra then tells us that $J(Q_1)+J(Q_1')$ and $J(Q_2)+J(Q_2')$
must define the same projector $J(P)$, i.e. the map $J$ is
well-defined on $\hCM$.

\begin{prop} $P \in \hCM$ is invariant under $J$ if and only if
$P \in S^4$.
\end{prop}

\proof Writing $Q' = (\bar W^\mu W^\nu)$, we have that
\[ P = Q+Q' = (\bar Z^\mu Z^\nu + \bar W^\mu W^\nu),\]
and that this supposed to be identified with the $2 \times 2$ block
decomposition
\[ P = \begin{pmatrix} A & B \\ B^\dagger & D \end{pmatrix}, \qquad
A=A^\dagger, D=D^\dagger,\; \Tr \, A + \Tr \, D = 2.\]
Here we shall use the notation of Proposition~3.1.  Examining $A$,
we have
\[A = a + \alpha \cdot \sigma = \begin{pmatrix} \bar Z^1 Z^1 + \bar W^1 W^1 & \bar Z^1 Z^2 + \bar Z^2 Z^1 \\ \bar Z^2 Z^1 + \bar Z^1 Z^2 & \bar Z^2
Z^2 + \bar W^2 W^2
\end{pmatrix},\]
and hence an identification
\[ a = \frac{1}{2}(\bar Z^1 Z^1 + \bar Z^2 Z^2 + \bar W^1 W^1 + \bar W^2
W^2),\]
\[ \alpha_3 = \frac{1}{2}(\bar Z^1 Z^1 - \bar Z^2 Z^2 + \bar W^1 W^1 - \bar W^2
W^2),\]
as well as the off-diagonal entries
\[ \alpha_1 = \frac{1}{2}(\bar Z^1 Z^2 - \bar Z^2 Z^1 + \bar W^1 W^2 - \bar
W^2 W^1),\]
\[ \alpha_2 = \frac{1}{2\imath}(\bar Z^1 Z^2 + \bar Z^2 Z^1 + \bar W^1 W^2 + \bar
W^2 W^1).\]
Clearly the relations $a=a^*$, $\alpha = \alpha^*$ hold under this
identification.  Under the involution $J$ we calculate that
\[J (a)=a, \qquad J(\alpha) = - \alpha.\]

Similarly we look at the block $D$,
\[D = d + \delta \cdot \sigma = \begin{pmatrix} \bar Z^3 Z^3 + \bar W^3 W^3 & \bar Z^3 Z^4 + \bar W^3 W^4 \\ \bar Z^4 Z^3 + \bar W^4 W^3 & \bar Z^4
Z^4 + \bar W^4 W^4
\end{pmatrix}.\]
The same computation as above shows that the relations $d=d^*$,
$\delta=\delta^*$ hold here, and moreover the trace relation implies
that $d=1-a$, in agreement with Section 3.1.  Under the involution
$J$ we also see that
\[J(d)=d, \qquad J(\delta) = - \delta.\]
Finally we look at the matrix $B$,
\[ B = t + \imath x \cdot \sigma = \begin{pmatrix} \bar Z^1 Z^3 + \bar W^1 W^3 & \bar Z^1 Z^4 + \bar W^1 W^4 \\ \bar Z^2 Z^3 + \bar W^2 W^3 & \bar
Z^2 Z^4 + \bar W^2 W^4
\end{pmatrix}.\]
Solving, we have the identification of generators
\[t = \frac{1}{2}(\bar Z^1 Z^3 + \bar Z^2 Z^4 + \bar W^1 W^3 + \bar W^2
W^4),\]
\[x_3 = \frac{1}{2\imath}(\bar Z^1 Z^3 - \bar Z^2 Z^4 + \bar W^1 W^3 - \bar W^2
W^4)\]
for the diagonal entries, as well as
\[x_1 = \frac{1}{2\imath}(\bar Z^2 Z^3 + \bar Z^1 Z^4 + \bar W^2 W^3 + \bar
W^1 W^4),\]
\[x_2 = \frac{1}{2}(\bar Z^2 Z^3 - \bar Z^1 Z^4 + \bar W^2 W^3 - \bar
W^1 W^4)\]
on the off-diagonal.  This is in agreement with the fact as in
Proposition~3.1 that the generators $t$, $x$ are not necessarily
Hermitian.  Moreover, it is a simple matter to compute that under
the involution $J$ we have
\[ J(t)=t^*, \qquad J(x) = x^*.\]
Overall, we see that $J$ has fixed points in $\hCM$ consisting of
those with coordinates subject to the additional constraints
$\alpha=\delta=0$, $t=t^*$, $x=x^*$.  Thus (upon verification of the
extra relations) the fixed points of $\hCM$ under $J$ are exactly
those lying in $S^4$, in accordance with proposition~3.2.  \endproof

In the notation of Proposition~3.3, the action of $J$ on
$\C[\hCM]$ is to map
\[ J(a)=a, \quad J(\alpha_3)=\alpha_3, \quad
J(\delta_3)=\delta_3, \quad J(\alpha)=-\alpha, \quad
J(\delta)=-\delta,\]
\[ J(w)=-\tilde w^*, \quad J(z)=\tilde z^*, \quad J(\tilde
w)=-w^*, \quad J(\tilde z)=z^*.\]
This may either be recomputed, or obtained simply by making the same
change of variables as was made in going from Proposition~3.1 to
Proposition~3.3.  The fixed points in these coordinates are those
with $\alpha=\delta=0$, $w^*=-\tilde w$, $z^* = \tilde z$, in
agreement with Proposition~3.3.

\begin{prop}  For each $P \in \hCM$ we have $P \in S^4$ if and only
if there exists $Q \in \CP^3$ such that $P = Q + J(Q)$.
\end{prop}

\proof By the previous proposition, $P \in S^4$ if and only if
$J(P)=P$.  Of course, the reverse direction of the claim is easy,
since if $P = Q + J(Q)$, we have $J(P) = J(Q) + J^2 (Q)
= P$. Conversely, given $P \in S^4$ with $J(P)=P$ we may write
$P=Q+Q'$ for some $Q$, $Q' \in \CP^3$ (as remarked already, this is
not a unique decomposition but given $Q$ we have $Q'=P-Q$, and
$J$ acts independently of this decomposition).  The result is now
obvious since $J$ is nondegenerate, so $Q'= J(Q'')$ for some
$Q''$, but we must have $Q''=Q$ since $J$ is an involution.
\endproof

As promised, there is a fibration of $\CP^3$ over $S^4$ given at the
coordinate algebra level by an inclusion $\C[S^4] \hookrightarrow
\C^-[\CP^3]$.  In terms of the $\C^-[\CP^3]$ coordinates
$a,b,c,x,y,z$ used in Proposition~3.8 and the $\C[S^4]$ coordinates
$z,w,a$ of Proposition~3.3 we have the following (we note the
overlap in notation between these propositions and rely on the
context for clarity).
\begin{prop}
There is an algebra inclusion
\[\eta : \C[S^4] \hookrightarrow \C^-[\CP^3]\]
given by
\[ \eta (a) = a+b, \quad \eta(z)=y+v^*,\quad \eta(w)=w-z^*.\]
\end{prop}

\proof That this is an algebra map is a matter of rewriting the
previous proposition in our explicit coordinates, for example that
$\eta(z) = y+v^*=y+J(y)$. The sole relation to investigate is the
image of the sphere relation $zz^*+ww^*=a(1-a)$. Applying $\eta$ to
the left hand side, we obtain
\[ \eta(zz^*+ww^*) = yy^* + yv + v^*y^* + v^*v + ww^* - wz - z^*w^* + z^*z. \]
Now using the relations of Proposition~3.5 we compute that
\[ ayv = yav = yx^*z = x^*yz = awz,\]
where we have relied upon the commutativity of the algebra.
Similarly one computes that $byv = bwz$, $cyv = cwz$, $dyv = dwz$,
so that adding these four relations yields that $yv=wz$ in
$\C^-[\CP^3]$.  Then finally using the relations in Proposition~3.8
we see that
\[\eta(zz^*+ww^*) = Y + V+ W+ Z = (a+b)(c+d) = (a+b)(1-(a+b)) =
\eta( a(1-a)),\] as required.  \endproof

We now look at the typical fibre $\CP^1$ of the fibration $\CP^3
\rightarrow S^4$, but now in the coordinate algebra picture.

\begin{prop}
The quotient of the algebra $\C^-[\CP^3]$ obtained by setting
$\eta(a)$, $\eta(z)$, $\eta(w)$ to be constant numerical values is
isomorphic to the coordinate algebra of a $\CP^1$.
\end{prop}

\proof If we suppose that we are in the patch where $a\ne 0$ in
$\C^-[\CP^3]$ then we can view $x$, $y$, $z$ as the variables and
$X=ab$, $Y=ac$, $Z=ad$ as the relations. The generators $u$, $v$,
$w$ are defined by the equations in Proposition~3.8 and the rest are
redundant.

Now suppose that $a+b=A$, a fixed real number, and $y+v^*=B$,
$w-z^*=C$, fixed complex numbers, such that
\[BB^*+CC^*=A(1-A)\]
(an element of $S^4$). Then we have just one equation
\[X=a(A-a)=(A/2)^2-s^2\]
if we set $s=a- A/2$. This is a $\CP^1$ of radius $A/2$ in place of
the usual radius $1/2$.

The equation $Y=ac$ is viewed as a definition of $c$.  The equation
$Z=ad$ is then equivalent to $Y+Z=a(1-A)$. We'll see that this is
automatic and that $y$, $z$ are uniquely determined by $x$, $a$ and
our fixed parameters $A$,$B$,$C$ so are not in fact free variables.

Indeed, $av=x^*z$ and $aw=x^*y$ determine $v$ and $w$ as mentioned
above, so our quotient is
\[aB=ay + z^*x,\quad aC=x^*y-az^*,\]
which implies that $a^2(BB^*+CC^*)=(a^2+X)(Y+Z)=aA(Y+Z)$, so
$Y+Z=a(1-A)$ necessarily holds if $A$, $B$, $C$ lie in $S^4$.

We also combine the equations to find $ax^*B=a^2C+z^*aA$ and
$aCx=aAy-a^2B$ so that at least if $A\ne 0$ we have $z$,$y$
determined. (In fact one has $By^*-z^*C=a(1-A)$ from the above so if
$z$ is determined then so is $y$ if $B$ is not zero etc).  Thus
$\C[\CP^1]$ is viewed inside $\C[\CP^3]$ in this patch as
\begin{equation} \label{CP1 fibre} \left( \begin{array}{cc} \begin{array}{cc} a & x \\ x^* & A-a
\end{array} & *  \\ * & * \end{array} \qquad \right), \qquad x^*x =
a(A-a),\end{equation} where the unspecified entries are determined
as above using the relations in terms of $x$ and $a$.

Similar analysis holds in the other coordinate patches, although we
shall not check this here as this is a well-known classical result.
In other patches we would see the various copies of $\C[\CP^1]$
appearing elsewhere in the above matrix. \endproof

This situation now provides us with yet another way to view the
instanton bundle.  Let $\mathcal{M}$ be a finite rank projective
$\C[\hCM]$-module.  Then $J$ induces a module map $J:\mathcal{M}
\rightarrow \mathcal{M}$ whose fixed point submodule is a finite
rank projective $\C[S^4]$-module.  In particular, if we take
$\mathcal{M}$ to be the $\C[\hCM]$-module given by the defining
tautological projector (\ref{twistorB}), then as explained above as
well as in that section, the fixed point submodule is precisely the
tautological bundle $\E = \C[S^4]^4e$ of Proposition~3.2 which
defines the instanton bundle over $S^4$.

Now the map $\eta: \C[S^4] \rightarrow \C^-[\CP^3]$ induces the
`push-out' of the $\C[S^4]$-module $\E$ along $\eta$ to obtain an
`auxiliary' $\C^-[\CP^3]$-module $\tilde \E$, given by viewing the
projector $e \in M_4(\C[S^4])$ as a projector $\tilde e \in
M_4(\C^-[\CP^3])$, so that $\tilde \E := \C^-[\CP^3]^4 \tilde e$,
giving a bundle over twistor space.  Explicitly, we have
\begin{eqnarray*} \tilde e &=& \begin{pmatrix} \eta(a) & 0 & \eta(z) & \eta(-w^*) \\
0 & \eta(a) & \eta(w) & \eta(z^*) \\ \eta(z^*) & \eta(w^*) &
\eta(1-a) & 0 \\ \eta(-w) & \eta(z) & 0 & \eta(1-a) \end{pmatrix} \\
&=& \begin{pmatrix} a+b & 0 & y+v^* & z-w^* \\
0 & a+b & w-z^* & y^*+v \\ y^*+v & w^*-z & 1-(a+b) & 0
\\ z^*-w & y+v^* & 0 & 1-(a+b) \end{pmatrix} \in M_4(\CC^-[\CP^3]).
\end{eqnarray*}
Moreover, if one sets $a+b=A$, $y+v^*=B$, $w-z^*=C$ for fixed real
$A$ and complex $B$, $C$ as in Proposition~3.13, then we have
\[ \tilde e = \begin{pmatrix} A & 0 & B & -C^* \\
0 & A & C & B^* \\ B^* & C^* & 1-A & 0
\\ -C & B & 0 & 1-A \end{pmatrix},\]
a constant projector of rank two.  Then viewing the fibre
$\C[\CP^1]$ as a subset of $\C[\CP^3]$ as in (\ref{CP1 fibre}), it
is easily seen that $\C[\CP^1]^4 \tilde e$ is a free
$\C[\CP^1]$-module of rank two. This is just the coordinate algebra
version of saying that for all $x = (A,B,C) \in S^4$ the instanton
bundle pulled back from $S^4$ to $\CP^3$ is trivial upon restriction
to each fibre $\hat{x} = \CP^1$, and we may thus see the instanton
bundle $\E$ over $\C[S^4]$ as coming from the bundle $\tilde \E$
over twistor space.  This is an easy example of the
\textit{Penrose-Ward transform}, which we shall discuss in more
detail later.

\section{The Quantum Conformal Group}
\label{section quantum conformal group} The advantage of writing
space-time and twistor space as homogeneous spaces in the language
of coordinate functions is that we are now free to apply the
standard theory of quantisation by a cocycle twist.

To this end, we recall that if $H$ is a Hopf algebra with coproduct
$\Delta: H \rightarrow H \otimes H$, counit $\ep: H \rightarrow \CC$
and antipode $S: H \rightarrow H$, then a two-cocycle $F$ on $H$ means
 $F: H \otimes H \rightarrow \CC$ which is convolution invertible
and unital (i.e. a 2-cochain) in the sense
\[ F(h\o,g\o)F^{-1}(h\t,g\t)=F^{-1}(h\o,g\o)F(h\t,g\t)=\eps(h)\eps(g)\]
(for some map $F^{-1}$) and obeys $\partial F = 1$ in the sense
\[ F(g\o,f\o)F(h\o,g\o f\t)F^{-1}(h\t g\thr,f\thr)F^{-1}(h\thr,g\fo)=\eps(f)\eps(h)\eps(g).\]
We have used Sweedler notation $\Delta(h) = h_{(1)} \otimes h_{(2)}$
and suppressed the summation. In this case there is a `cotwisted'
Hopf algebra $H_F$ with the same coalgebra structure and counit as
$H$ but with  modified product $\bu$ and antipode $S_F$ \cite{Ma:book}
\begin{equation}
\label{twisted product}
h \bu g = F(h_{(1)}\otimes g_{(1)})\, h_{(2)}
g_{(2)}\, F^{-1}(h_{(3)}\otimes g_{(3)})\end{equation}
\[ S_F(h)=U(h\o)Sh\t U^{-1}(h\t),\quad U(h)=F(h\o,Sh\t)\]
for $h,g \in H$, where we use the product and antipode of $H$ on the right hand
sides and $U^{-1}(h\o)U(h\t)=\eps(h)=U(h\o)U^{-1}(h\t)$ defines the inverse functional. If $H$ is a coquastriangular Hopf algebra then so is $H_F$. In
particular, if $H$ is commutative then $H_F$ is cotriangular with
'universal R-matrix' and induced (symmetric) braiding given by
\[ R(h,g)=F(g\o,h\o)F^{-1}(h\t,g\t),\quad \Psi_{V,W}(v\tens w)=R(w\o,v\o) w\t\tens v\t \]
for any two left comodules $V,W$. We use the Sweedler notation for
the left coactions as well. In the cotriangular case one has
$\Psi^2={\rm id}$, so every object on the category of
$H_F$-comodules inherits nontrivial statistics in which
transposition is replaced by this non-standard transposition.

The nice property of this construction is that the category of
$H$-comodules is actually equivalent to that of $H_F$-comodules, so
there is an invertible functor which `functorially quantises' any
construction in the first category (any $H$-covariant construction)
to give an $H_F$-covariant one. So not only is the classical Hopf
algebra $H$ quantised but also any $H$-covariant construction as
well. This is a particularly easy example of the `braid statistics
approach' to quantisation, whereby deformation is achieved by
deforming the category of vector spaces to a braided one
\cite{Ma:book}.

In particular, if $A$ is a left $H$-comodule algebra, we
automatically obtain a left $H_F$-comodule algebra $A_F$ which as a
vector space is the same as $A$, but has the modified product
\begin{equation}
\label{twisted comodule product}
a \bu b = F(a\o, b\o) a\t b\t ,
\end{equation}
for $a,b \in A$, where we have again used the Sweedler notation
$\Delta_L (a) = a\o \otimes a\t$ for the coaction $\Delta_L : A
\rightarrow H\tens A$. The same applies to any other covariant
algebra. For example if $\Omega(A)$ is an $H$-covariant differential
calculus (see later) then this functorially quantises as
$\Omega(A_F):=\Omega(A)_F$ by this same construction.

Finally, if $H'\to H$ is a homomorphism of Hopf algebras then any
cocycle $F$ on $H$ pulls back to one on $H'$ and as a result one has
a homomorphism $H'_F\to H_F$. In what follows we take $H=\CC[\CC^4]$
(the translation group of $\CC^4$) and $H'$ variously the coordinate
algebras of $\tilde K,\tilde H, \GL_4$.

In particular, since the group $\GL_4$ acts on the quadric $\TCM$,
we have (as in section \ref{conformal space time}) a coaction
$\Delta_L$ of the coordinate ring $\CC[\GL_4]$ on $\CC[\TCM]$, and
we shall first deform this picture. In order to do this we note
first that the conformal transformations of $\hCM$ break down into
compositions of translations, rotations, dilations and inversions.
Written with respect to the aforementioned double null coordinates,
$\GL_4$ decomposes into $2 \times 2$ blocks
\begin{equation}
\label{generators of conf group}
\left( \begin{array}{cc}
\gamma & \tau \\ \sigma & \tilde{\gamma} \end{array} \right)
\end{equation}
with overall non-zero determinant, where the entries of $\tau$
constitute the translations and the entries of $\sigma$ contain the
inversions. The diagonal blocks $\gamma \times \tilde \gamma$
constitute the space-time rotations as well as the dilations.
Writing $M_2:=M_2(\CC)$, $\GL_4$ decomposes as the subset of nonzero
determinant
\[ \GL_4 \subset \CC^4 \rtimes (M_2 \times M_2) \ltimes \CC^4 \]
where the outer factors denote $\sigma, \tau$ and $\gamma \times
\tilde \gamma \in M_2 \times M_2$.
In practice it is convenient to work in a `patch' $\GL_4^-$ where $\gamma$ is assumed
invertible. Then by factorising the matrix we deduce that
\[ \det\begin{pmatrix}\gamma & \tau\\ \sigma & \tilde \gamma\end{pmatrix}=\det(\gamma)\det (\tilde\gamma-\sigma\gamma^{-1}\tau)\]
which is actually a part of a universal formula for determinants of
matrices with entries in a noncommutative algebra (here the algebra
is $M_2$ and we compose with the determinant map on this algebra).
We see that as a set, $\GL_4^-$ is $\CC^4\times \GL_2\times
\GL_2\times \CC^4$, where the two copies of $\GL_2$ refer to
$\gamma$ and $\tilde\gamma-\sigma\gamma^{-1}\tau$. There is of
course another patch $\GL_4^+$ where we similarly assume
$\tilde\gamma$ invertible.

In terms of coordinate functions for $\C[\GL_4]$ we therefore have
four matrix generators $\tau,\sigma,\gamma,\tilde\gamma$ organised
as above. These together have a matrix form of coproduct
\[ \Delta \left( \begin{array}{cc} \gamma & \tau \\ \sigma & \tilde{\gamma} \end{array}
\right) = \left( \begin{array}{cc} \gamma & \tau \\ \sigma &
\tilde{\gamma} \end{array} \right) \otimes \left(
\begin{array}{cc} \gamma & \tau \\ \sigma & \tilde{\gamma}
\end{array} \right). \]
In the classical case the generators commute and an invertible
element $D$ obeying $D=\det a$ is adjoined. For $\CC[\GL_4^-]$ we
instead adjoin inverses to $d=\det(\gamma)$ and $\tilde
d=\det(\tilde\gamma-\sigma\gamma^{-1}\tau)$.

We focus next on the translation sector $H=\CC[\CC^4]$ generated by
some $t^{A}_{\pA}$, where $\pA\in\{3,4\}$ and $A\in\{1,2\}$
to line up with our conventions for $\GL_4$. These generators have a
standard additive coproduct. We let $\del^{\pA}_A$ be the Lie
algebra of translation generators dual to this, so
\[ \la \del^{\pA}_{A},t^B_{\pB} \ra=\delta^{\pA}_{\pB}\delta^B_A\]
which extends to the action on products of the $t^A_{\pA}$ by
differentiation and evaluation at zero (hence the notation). In this
notation the we use cocycle
\[ F(h,g) = \la\textup{exp} ({\imath\over
2} \theta^{AB}_{\pA \pB} \, \p^{\pA}_A \otimes \p^{\pB}_B),h\tens g
\ra
.\]
Cotwisting here does not change $H$ itself,  $H=H_F$, because its
coproduct is cocommutative (the group $\C^4$ is Abelian) but it
twists $A=\CC[\CC^4]$ as a comodule algebra into the Moyal plane.
This is by now well-known both in the module form and the above
comodule form. We now pull this cocycle back to $\CC[\GL_4]$, where
it takes the same form as above on the generators $\tau^A_{\pA}$
(which project onto $t^A_{\pA}$). The pairing extends as zero on the
other generators. One can view the $\del^{\pA}_A$ in the Lie algebra
of $\GL_4$ as the nilpotent  $4 \times 4$ matrices with entry $1$ in
the $A, \pA$ position for some $A = 1,2$, $\pA = 3,4$ and zeros
elsewhere, extending the above picture. Either way, one computes
\begin{eqnarray*}
F^{\mu \alpha}_{\nu \beta} = F(a^\mu_\nu, a^\alpha_\beta)
&=& \la \textup{exp} ({\imath\over
2} \theta^{AB}_{\pA \pB} \p^{\pA}_A \otimes \p^{\pB}_B),a^\mu_\nu\tens a^\alpha_\beta \ra \\
&=& \delta^{\mu}_\nu \delta^\alpha_\beta + {\imath\over
2} \theta^{AB}_{\pA \pB} \delta^\mu_A \delta^\alpha_B \delta_\nu^{\pA} \delta_\beta^{\pB} \\
&=& \delta^\mu_\nu \delta^\alpha_\beta  + {\imath\over 2}
\theta^{\mu \alpha}_{\nu \beta},
\end{eqnarray*}
where it is understood that $\theta^{\mu \alpha}_{\nu \beta}$ is
zero when $\{ \mu, \alpha \} \ne \{ 1,2 \}, \{ \nu, \beta \} \ne \{
3,4 \}$.  We also compute
\[ U(a^\mu_\nu)=F(a^\mu_a,Sa^a_\nu)= \la \textup{exp} (-{\imath\over
2} \theta^{AB}_{\pA \pB} \p^{\pA}_A \otimes \p^{\pB}_B),a^\mu_\alpha\tens a^\alpha_\nu \ra=
 \delta^\mu_\nu   - {\imath\over 2}
\theta^{\mu a}_{a \beta}=\delta^\mu_\nu.\]

Then following equations (\ref{twisted product}) the deformed
coordinate algebra $\CC_F [\GL_4]$ has undeformed antipode on the
generators and deformed product
\[ a^\mu_\nu \bu a^\alpha_\beta = F^{\mu \alpha}_{mn} a^m_p a^n_q
F^{-1}{}^{pq}_{\nu \beta}\]%
where $a^\mu_\nu \in \CC[\GL_4]$ are the generators of the classical
algebra. The commutation relations can be written in R-matrix form (as for any
matrix coquasitriangular Hopf algebra) as
\[ R^{\mu\nu}_{\alpha\beta}a^\alpha_\gamma\bullet a^\beta_\delta=a^\nu_\beta\bullet a^\mu_\alpha R^{\alpha\beta}_{\gamma\delta}, \quad R^{\mu\nu}_{\alpha\beta}=F^{\nu\mu}_{\delta\gamma}F^{-1}{}^{\gamma\delta}_{\alpha\beta}\]
where in our particular case
\[ R^{\mu\alpha}_{\nu\beta}=\delta^\mu_\nu\delta^\alpha_\beta-\imath\theta^-{}^{\mu\alpha}_{\nu\beta},\quad \theta^-{}^{\mu\alpha}_{\nu\beta}={1\over 2}(\theta^{\mu\alpha}_{\nu\beta}-\theta^{\alpha\mu}_{\beta\nu})\]
has the same form but now with only the antisymmetric part of $\theta$ in the sense shown.

We give the resulting relations explicitly in the
$\gamma,\tilde\gamma,\sigma,\tau$ block form  (\ref{generators of
conf group}). These are in fact all $2\times 2$ matrix relations
with indices $A,A'$ etc. as explained but when no confusion can
arise we write the indices in an apparently $\GL_4$ form. For
example, in writing $\gamma^\mu_\nu$ it is implicit that $\mu, \nu
\in \{1,2\}$, whereas for $\sigma^\mu_\nu$ it is understood that
$\mu \in \{ 3,4\}$ and $\nu \in \{1,2\}$.

\begin{thm}
\label{quantum conformal group}
The quantum group coordinate algebra $\CC_F[\GL_4]$ has deformed product
\begin{eqnarray*}
&\tau^\mu_{\nup} \bu \tau^\alpha_{\betap} = \tau^\mu_{\nup}
\tau^\alpha_{\betap} ~+~ {\imath\over 2} \theta^{\mu \alpha}_{\pc
\pd} ~\tilde \gamma^{\pc}_{\alphap} \tilde \gamma^{\pd}_{\betap} ~-~
{\imath\over 2} \gamma^\mu_c \gamma^\alpha_d ~\theta^{cd}_{\nup
\betap} ~+~
{1\over 4} \theta^{\mu \alpha}_{ab} ~\sigma^a_c \sigma^b_d ~\theta^{cd}_{\nup \betap}, \\
&\gamma^\mu_\nu \bu \tau^\alpha_{\betap} = \gamma^\mu_\nu
\tau^\alpha_{\betap} ~+~ {\imath\over 2}\theta^{\mu \alpha}_{c
d}~ \sigma^{c}_\nu \tilde \gamma^{d}_{\betap}, \qquad
\tau^\mu_{\nup} \bu \gamma^\alpha_\beta ~~=~~ \tau^\mu_{\nup}
\gamma^\alpha_\beta ~+~
{\imath\over 2}\theta^{\mu \alpha}_{cd} ~\tilde \gamma^{c}_{\nup} \sigma^{d}_\beta, \\
&\tilde \gamma^{\mup}_{\nup} \bu \tau^\alpha_{\betap} = \tilde
\gamma^{\mup}_{\nup} \tau^\alpha_{\betap} ~-~ {\imath\over 2}
\sigma^{\mup}_c \gamma^\alpha_d ~\theta^{cd}_{\nup \betap}, \qquad
\tau^\mu_{\nup} \bu \tilde \gamma^{\alphap}_{\betap} ~~=~~
\tau^\mu_{\nup} \tilde \gamma^{\alpha}_{\betap} ~-~ {\imath\over
2} \gamma^\mu_c \sigma^{\alphap}_d ~\theta^{cd}_{\nup \betap}, \\
&\gamma^\mu_\nu \bu \gamma^\alpha_\beta = \gamma^\mu_\nu
\gamma^\alpha_\beta ~+~ {\imath\over 2} \theta^{\mu \alpha}_{c
d} ~\sigma^{c}_\nu \sigma^{d}_\beta, \qquad
\tilde \gamma^{\mup}_{\nup} \bu \tilde \gamma^{\alphap}_{\betap}
~~=~~ \tilde \gamma^{\mup}_{\nup} \tilde \gamma^{\alphap}_{\betap}
~-~ {\imath\over 2} \sigma^{\mup}_c \sigma^{\alphap}_d
~\theta^{cd}_{\nup \betap}
\end{eqnarray*}
with the remaining relations, antipode and coproduct undeformed on the generators.  The quantum group is generated by matrices $\gamma,\tilde\gamma,\tau,\sigma$ of generators with commutation relations
\[ [\gamma^\mu_\nu,\gamma^\alpha_\beta]_\bullet=\imath\theta^-{}^{\mu\alpha}_{cd}\sigma^c_\nu\bullet \sigma^d_\beta,\quad [\tilde\gamma^\mu_\nu,\tilde\gamma^\alpha_\beta]_\bullet=-\imath \sigma^\alpha_d\bullet \sigma^\mu_c \theta^{-}{}^{cd}_{\nu\beta}\]
\[ [\gamma^\mu_\nu,\tau^\alpha_\beta]_\bullet=\imath\theta^-{}^{\mu\alpha}_{cd}\sigma^c_\nu\bullet \tilde\gamma^d_\beta,\quad  [\tilde\gamma^\mu_\nu,\tau^\alpha_\beta]_\bullet =-\imath \gamma^\alpha_d\bullet \sigma^\mu_c \theta^{-}{}^{cd}_{\nu\beta}\]
\[ [\tau^\mu_\nu,\tau^\alpha_\beta]_\bullet=\imath\theta^-{}^{\mu\alpha}_{cd}\tilde\gamma^c_\nu\bullet \tilde\gamma^d_\beta-\imath \gamma^\alpha_d\bullet \gamma^\mu_c \theta^{-}{}^{cd}_{\nu\beta}\]
and  a certain determinant inverted.
\end{thm}
\proof Finishing the computations above with the explicit form of $F$ we have
\[ a^\mu_\nu\bullet a^\alpha_\beta=a^\mu_\nu a^\alpha_\beta +{\imath\over 2}\theta^{\mu\alpha}_{cd}a^c_\nu a^d_\beta-{\imath\over 2}a^\mu_c a^\alpha_d\theta^{cd}_{\nu\beta}+{1\over 4}\theta^{\mu\alpha}_{ab}a^a_c a^b_d \theta^{cd}_{\nu\beta}.\]
Noting that $\theta^{\mu\alpha}_{\nu\beta}=0$ unless $\mu,\alpha\in\{1,2\}$ and $\nu,\beta\in\{3,4\}$ we can write these for the $2\times 2$ blocks as shown. For the commutation relations we have similarly
\[ [a^\mu_\nu,a^\alpha_\beta]_\bullet=\imath\theta^-{}^{\mu\alpha}_{cd}a^c_\nu\bullet a^d_\beta-\imath a^\alpha_d\bullet a^\mu_c \theta^{-}{}^{cd}_{\nu\beta}\]
which we similarly decompose as stated. Note that different terms
here drop out due to the range of the indices for nonzero
$\theta_-$, which are same as for $\theta$. There is in principle a
formula also for the determinant written in terms of the $\bu$
product. It can be obtained via braided `antisymmetric tensors' from
the $R$-matrix and will necessarily be product of $2\times 2$
determinants in the `patches' where $\gamma$ or $\tilde\gamma$ are
invertible in the noncommutative algebra. \endproof

One may proceed to compute these more explicitly, for example
\[ [\gamma^\mu_\nu,\gamma^\alpha_\beta]_\bullet=\imath\theta^-{}^{\mu\alpha}_{34}\sigma^3_\nu \sigma^4_\beta
+\imath\theta^-{}^{\mu\alpha}_{43}\sigma^3_\beta \sigma^4_\nu+\imath\theta^-{}^{\mu\alpha}_{33}\sigma^3_\nu \sigma^3_\beta+\imath\theta^-{}^{\mu\alpha}_{44}\sigma^4_\nu \sigma^4_\beta\]
and so forth.

We similarly calculate the resulting products on
the coordinate algebras of the deformed homogeneous spaces.  Indeed,
using equation (\ref{twisted comodule product}), we have the
following results.

\begin{prop}
\label{twisted space-time algebra} The covariantly twisted algebra
$\CC_F[\TCM]$ has the deformed product
\[ x^{\mu \nu} \bu x^{\alpha \beta}=x^{\mu\nu}x^{\alpha\beta}+{\imath\over 2}(\theta^{\mu\beta}_{ad}x^{a\nu}x^{\alpha d}+\theta^{\mu\alpha}_{ac}x^{a\nu}x^{c\beta}+\theta^{\nu\beta}_{bd}x^{\mu b}x^{\alpha d}+\theta^{\nu\alpha}_{bc}x^{\mu b}x^{c\beta})\]
\[ -{1\over 4}\left( \theta^{\nu\alpha}_{bc}\theta^{\mu\beta}_{ad} + \theta^{\nu\beta}_{bd}\theta^{\mu\alpha}_{ac} \right)x^{ab}x^{cd}\]
and is isomorphic to the subalgebra
\[ \CC_F[\GL_4]^{\CC_F[\tilde H]} \]
where $F$ is pulled back to $\CC[\tilde H]$. Products of generators with $t=x^{34}$ are undeformed.
\end{prop}
\proof  The isomorphism $\CC_F[\TCM] \cong
\CC_F[\GL_4]^{\CC_F[\tilde H]}$ is a consequence of the
functoriality of the cocycle twist.  The deformed product is simply
a matter of calculating the twisted product on $\CC_F[\TCM]$. The
coaction $\Delta_L (x^{\mu \nu}) = a^\mu_a a^\nu_b\tens x^{ab} $
combined with the formula (\ref{twisted comodule product}) yields
\[ x^{\mu \nu} \bu x^{\alpha \beta} = F(a^\mu_a a^\nu_b,a^\alpha_c
a^\beta_d) x^{ab} x^{cd} =F^{\mu\alpha}_{mn}F^{m\beta}_{ap}F^{\nu
n}_{qc}F^{qp}_{bd}x^{ab} x^{cd}, \]
using that $F$ in our particular case is multiplicative (a Hopf
algebra bicharacter on $\CC[\CC^4]$ and hence when pulled back to
$\CC[GL_4]$). Alternatively, one may compute it directly from the
original definition as an exponentiated operator, going out to
$\del\del\tens\del\del$ terms before evaluating at zero in $\CC^4$.
Either way we have the result stated when  we recall that
$\theta^{\mu \alpha}_{bd}$ is understood to be zero unless $\{\mu,
\alpha\} = \{1,2\}$, $\{b,d\} = \{3,4\}$.  We note also the
commutation relations
\[ x^{\alpha\beta}\bullet x^{\mu\nu}=  R^{\mu\alpha}_{mn}R^{m\beta}_{ap}R^{\nu n}_{qc}R^{qp}_{bd}x^{ab} \bullet x^{cd}  \]
following from $b\bullet a=R(a\o,b\o)a\t\bullet b\t$ computed in the
same was as above but now with $R$ in place of $F$. Since
$R^{-1}=R_{21}$ these relations may be written in a `reflection'
form on regarding $x$ as a matrix. Finally, since $\theta$ is zero
when $\{\mu, \alpha\} \ne \{1,2\}$ we see that the $\bullet$ product
of the generator $t = x^{34}$ with any other generator is
undeformed, which also implies that $t$ is central in the deformed
algebra. \endproof

Examining the resulting relations
associated to the twisted product more closely, one finds
\[
[z,\tz]_\bu ~=~ [-x^{23},x^{14}]_\bu = -{\imath\over
2}\theta^{21}_{43} x^{43}x^{34} + {\imath\over 2} \theta^{12}_{34}
x^{34}x^{43} =\imath \theta^-{}^{21}_{43} t^2,
\]
\[
[w,\tw]_\bu ~=~ [-x^{13},x^{24}]_\bu = -{\imath\over 2}
\theta^{12}_{43} x^{43}x^{34} + {\imath\over 2} \theta^{21}_{34}
x^{34}x^{43} =\imath\theta^-{}^{12}_{43}t^2,\]
with the remaining commutators amongst these affine Minkowski space
generators undeformed. Since products with $t$ are undeformed, the
above can  be viewed as  the commutation relations among the affine
generators $t,w,\tilde w,z,\tilde z$ of $\C_F[\TCM]$. The
commutation relations for the $s$ generator are
\begin{eqnarray*}
[s,z]_\bu &=& [x^{12},-x^{23}]_\bu = -{\imath\over
2}\theta^{12}_{ac} x^{a2}x^{c3} + {\imath \over 2}
\theta^{21}_{ac}x^{a3}x^{c2} \\ &=& -{\imath\over 2}\theta^{12}_{34}
x^{32}x^{43} -{\imath\over 2}\theta^{12}_{44}x^{42}x^{43} + {\imath
\over 2} \theta^{21}_{43}x^{43}x^{32} + {\imath \over 2}
\theta^{21}_{44}x^{43}x^{42} \\ &=& \imath\theta^-{}^{12}_{34} zt +
\imath\theta^-{}^{21}_{44} \tw t, \end{eqnarray*}
for example, as well as
\[ [s,\tz]_\bu = \imath \theta^-{}^{12}_{33} wt +  \imath\theta^-{}^{21}_{43} \tz
t,\]
\[ [s,w]_\bu = \imath\theta^-{}^{12}_{43} wt +  \imath\theta^-{}^{21}_{44} \tz
t,\]
\[ [s,\tw]_\bu = \imath\theta^-{}^{12}_{33} zt +  \imath\theta^-{}^{21}_{34} \tw
t. \]
Again, we may equally well use the $\bu$ product on the right hand
side of each equation.

\begin{prop}
\label{twisted twistor algebra} The twisted algebra $\CC_F[\tilde
T]$ has deformed product
\[ Z^\mu \bu Z^\nu = Z^\mu Z^\nu + {\imath\over
2} \theta^{\mu \nu}_{ab} Z^a Z^b. \]
It is isomorphic to the subalgebra
\[ \CC_F[\GL_4]^{\CC_F[\tilde K]} \]
where $F$ is pulled back to $\CC[\tilde K]$. Products of generators with $Z^3,Z^4$ are undeformed.
\end{prop}
\noindent \proof The isomorphism $\CC_F[\tilde
T]\cong\CC_F[\GL_4]^{\CC_F[\tilde K]}$ is again a consequence of the
theory of cocycle twisting.  An application of equation
(\ref{twisted comodule product}) gives the new product
\[ Z^\nu\bullet Z^\mu=F^{\nu\mu}_{ab}Z^a Z^b =F^{\nu\mu}_{ab} F^{-1}{}^{ba}_{cd} Z^c\bullet Z^d=R^{\mu\nu}_{cd}Z^c\bullet Z^d \]
and we compute the first of these explicitly. The remaining
relations tell us that this is a `braided vector space' associated
to the $R$-matrix, see \cite[Ch. 10]{Ma:book}. As before, the form of $\theta$ implies that products with  $Z^3,Z^4$ are undeformed. Hence these are central.   \endproof

We conclude that the only nontrivial commutation relation is the
$Z^1$-$Z^2$ one, which we compute explicitly as
\[ [Z^1,Z^2]_\bu =-\imath\theta^-{}^{21}_{ab} Z^a\bullet Z^b =\imath(\theta^-{}^{12}_{34}+\theta^-{}^{12}_{43})Z^3Z^4+\imath\theta^-{}^{12}_{33}Z^3Z^3+\imath\theta^-{}^{12}_{44}Z^4Z^4.\]
where we could as well use the $\bullet$ on the right.

We remark that although our deformation of the conformal group is
different from and more general than that previously obtained in
\cite{kch:ncts} (which insists on only first order terms in the
deformation parameter), it is of note that the deformed commutation
relations associated to twistor space and conformal space-time are
in agreement with those proposed in recent literature
\cite{ns:ins,kko:nitt,kch:ncts,sb:thesis} {\em provided} one supposes that the four parameters
\[ \theta^-{}^{12}_{33}=\theta^-{}^{12}_{44}=\theta^-{}^{11}_{34}=\theta^-{}^{22}_{34}=0.\]
This says that $\theta^-{}^{AB}_{\pA\pB}$ as a $4\times 4$ matrix with rows $A\pA$ and columns $B\pB$  (the usual presentation) has the
form
\[ \theta^-=\begin{pmatrix}0&0&0&\theta_1\\ 0&0&-\theta_2& 0\\ 0&\theta_2&0&0\\ -\theta_1&0&0&0\end{pmatrix},\quad \theta_1=\theta^-{}^{12}_{34}, \quad \theta_2=\theta^-{}^{21}_{34}.\]
The quantum group deformation we propose  then agrees with
\cite{kch:ncts} on the generators $\gamma^\mu_\nu,
\sigma^\alpha_\beta$ and this in fact is the reason that the
space-time and twistor algebras then agree, since  their  generators
may be viewed as living in the subalgebra generated by the first two
columns of $a$ \textit{via} the isomorphisms given in propositions
\ref{twisted space-time algebra} and \ref{twisted twistor algebra}.

Finally we give the commutation relations in the twisted coordinate
algebra $\C_F[\tilde \F]$ of the correspondence space, which may be
computed either by viewing it as a twisted comodule algebra for
$\C_F [\GL_4]$ or by identification with the appropriate subalgebra
of $\C_F [\GL_4]$ and calculating there.  Either way, one obtains
\begin{equation} \label{correspondence space relations}
[s,Z^1]_\bu = -\imath \theta^-{}^{21}_{33} w Z^3 - \imath
\theta^-{}^{21}_{34} w Z^4 + \imath \theta^-{}^{21}_{43} \tz Z^3 +
\imath \theta^-{}^{21}_{44} \tz Z^4 + \imath \theta^-{}^{11}_{34} t
Z^2,\end{equation}
\[ [s,Z^2]_\bu = \imath \theta^-{}^{12}_{33} z Z^3 + \imath \theta^-{}^{12}_{34}
z Z^4 - \imath \theta^-{}^{12}_{43} \tw Z^3 - \imath
\theta^-{}^{12}_{44} \tw Z^4 + \imath \theta^-{}^{22}_{34} t Z^1,\]
\[ [z,Z^1]_\bu = \imath \theta^-{}^{21}_{43} tZ^3 + \imath \theta^-{}^{21}_{44} tZ^4, \qquad [z,Z^2]_\bu = \imath \theta^-{}^{22}_{43} tZ^3,\]
\[ [\tz,Z^1]_\bu = \imath \theta^-{}^{11}_{34} tZ^4, \qquad [\tz,Z^2]_\bu = \imath \theta^-{}^{12}_{33} tZ^3 + \imath \theta^-{}^{12}_{34}
tZ^4,\]
\[ [w,Z^1]_\bu = \imath \theta^-{}^{11}_{43} tZ^3, \qquad [w,Z^2]_\bu = \imath \theta^-{}^{12}_{43} tZ^3 + \imath \theta^-{}^{12}_{44} tZ^4, \]
\[ [\tw,Z^1]_\bu = \imath \theta^-{}^{21}_{33} tZ^3 + \imath \theta^-{}^{21}_{34} tZ^4, \qquad [\tw,Z^2]_\bu = \imath \theta^-{}^{22}_{34}
tZ^4,\]
where we may equally write $\bu$ on the right hand side of each
relation.  The generators $t,Z^3,Z^4$ are of course central.  The
relations (\ref{alpha planes}) twist by replacing the old product by
$\bu$.  In terms of the old product they become
\begin{eqnarray}
\label{twisted alpha planes}
&& \tz Z^3 + w Z^4 -t Z^1=0,\\
\nonumber && \tw Z^3 + z Z^4 - t Z^2=0,\\
\nonumber && s Z^3+w Z^2+ -z Z^1 + {\imath \over
2}((\theta^{12}_{43}-
\theta^{21}_{43})tZ^3 + \theta^-{}^{12}_{44}tZ^4)=0, \\
\nonumber && s Z^4-\tz Z^2+\tw Z^1 + {\imath \over 2}
(\theta^-{}^{21}_{33} tZ^3 + (\theta^{21}_{34} -
\theta^{12}_{34})tZ^4)=0.
\end{eqnarray}

\section{Quantum Differential Calculi on $\CC_F[\GL_4]$, $\CC_F[\TCM]$ and $\CC_F[\tilde T]$}
\label{theta quantum calculi}

We recall that a \textit{differential calculus} of an algebra $A$
consists of an $A$-$A$-bimodule $\Omega^1 A$ and a map $\D :A
\rightarrow \Omega^1 A$ obeying the Leibniz rule such that $\Omega^1
A$ is spanned by elements of the form $a \D b$.  Every unital
algebra has a universal calculus $\Omega^1_{un}=\ker\mu$ where $\mu$
is the product map of $A$. The differential is $\D_{un}(a)=1\tens
a-a\tens 1$. Any other calculus is a quotient of $\Omega_{un}^1$ by
a sub-bimodule $N_A$.

When $A$ is a Hopf algebra, it coacts on itself by left and right translation
\textit{via} the coproduct $\Delta$: we say a calculus on $A$ is
\textit{left covariant} if this coaction extends to a left coaction
$\Delta_L: \Omega^1 A \rightarrow A \otimes \Omega^1 A$ such that
$\D$ is an intertwiner and $\Delta_L$ is a bimodule map, that is
\[\Delta_L (\D a) = ({\rm id} \otimes \D) \Delta (a), \]
\[ a \cdot \Delta_L (\omega) = \Delta_L (a \cdot \omega), \quad \Delta_L
(\omega) \cdot b = \Delta_L (\omega \cdot b) \]
for all $a,b \in A, \omega \in \Omega^1 A$, where $A$ acts on
$A\tens\Omega^1A$ in the tensor product representation.  Equivalently, $\Delta_L(a\cdot\omega)=(\Delta a)\cdot\Delta_L(\omega)$ etc., with the second product as an $A\tens A$-module. We
then say that a one-form $\omega \in \Omega^1 A$ is \textit{left
invariant} if it is invariant under left translation by $\Delta_L$.
Of course, similar definitions may be made with `left' replaced by
`right'. It is bicovariant if both definitions hold and the left and
right coactions commute. We similarly have the notion of the
calculus on an $H$-comodule algebra $A$ being $H$-covariant, namely
that the coaction extends to the calculus such that it commutes with
$\D$ and is multiplicative with respect to the bimodule product.

We now quantise the differential structures on our spaces and groups
by the same covariant twist method. For groups the important thing
to know is that the classical exterior algebra of differential forms
$\Omega(\GL_4)$ (in our case) is a super-Hopf algebra where the
coproduct on degree zero elements is that of $\CC[\GL_4]$ while on
degree one it is $\Delta_L+\Delta_R$ for the classical coactions
induced by left and right translation (so $\Delta_L \D a=a\o\tens\D
a\t$ etc.) We view $F$ as a cocycle on this super-Hopf algebra  by
extending it as zero, and make a cotwist in the super-algebra
version of the cotwist of $\CC[\GL_4]$. Then $\Omega(\CC_F[\GL_4])$
has the bimodule and wedge products
\[ a^\mu_\nu\bullet \D a^\alpha_\beta=F^{\mu \alpha}_{mn} a^m_p\D a^n_q F^{-1}{}^{pq}_{\nu\beta},\quad
\D  a^\mu_\nu\bullet a^\alpha_\beta=F^{\mu \alpha}_{mn} (\D a^m_p)
a^n_q F^{-1}{}^{pq}_{\nu\beta},\]
\[ \D a^\mu_\nu\bullet \D a^\alpha_\beta=F^{\mu \alpha}_{mn}\D a^m_p\wedge \D a^n_q F^{-1}{}^{pq}_{\nu\beta},\]
while $\D$ itself is not deformed.  The commutation relations are
\[  R^{\mu\alpha}_{ab}a^a_\nu\bullet \D a^b_\beta= \D a^\alpha_b\bu a^\mu_b R{}^{ab}_{\nu\beta},\quad
R^{\mu\alpha}_{ab}\D a^a_\nu\bullet \D a^b_\beta=- \D a^\alpha_b\bu \D a^\mu_b R{}^{ab}_{\nu\beta}.\]
In terms of the decomposition
(\ref{generators of conf group}) the deformed products come out as
\[ \gamma^\mu_\nu \bu \D \gamma^\alpha_\beta = \gamma^\mu_\nu \D
\gamma^\alpha_\beta+ {\imath\over 2} \theta^{\mu \alpha}_{ab}
\sigma^a_\nu \D \sigma^b_\beta, \quad
\D \gamma^\mu_\nu \bu \gamma^\alpha_\beta = (\D \gamma^\mu_\nu)
\gamma^\alpha_\beta+ {\imath\over 2} \theta^{\mu \alpha}_{ab} (\D
\sigma^a_\nu) \sigma^b_\beta,\]
\[ \D \gamma^\mu_\nu \bu \D \gamma^\alpha_\beta = \D \gamma^\mu_\nu \wedge \D
\gamma^\alpha_\beta+ {\imath\over 2} \theta^{\mu \alpha}_{ab} \D
\sigma^a_\nu \wedge \D \sigma^b_\beta, \]
\[ \tilde \gamma^\mu_\nu \bu \D \tilde \gamma^\alpha_\beta = \tilde \gamma^\mu_\nu \D
\tilde \gamma^\alpha_\beta- {\imath\over 2}\sigma^\alpha_c \D
\sigma^\alpha_d \theta^{cd}_{\nu \beta}, \quad
\D \tilde \gamma^\mu_\nu \bu \tilde \gamma^\alpha_\beta = (\D \tilde
\gamma^\mu_\nu) \tilde \gamma^\alpha_\beta- {\imath\over 2} (\D
\sigma^\alpha_c) \sigma^\alpha_d \theta^{cd}_{\nu \beta},\]
\[ \D \tilde \gamma^\mu_\nu \bu \D \tilde \gamma^\alpha_\beta = \D \tilde \gamma^\mu_\nu \wedge \D
\tilde \gamma^\alpha_\beta- {\imath\over 2}\D \sigma^\mu_c \wedge \D
\sigma^\alpha_d \theta^{cd}_{\nu \beta}, \]
\[ \gamma^\mu_\nu \bu \D \tau^\alpha_\beta = \gamma^\mu_\nu \D
\tau^\alpha_\beta + {\imath\over 2} \theta^{\mu \alpha}_{ab}
\sigma^a_\nu \D \tilde \gamma^b_\beta, \quad
\D \gamma^\mu_\nu \bu \tau^\alpha_\beta = (\D \gamma^\mu_\nu)
\tau^\alpha_\beta + {\imath\over 2} \theta^{\mu \alpha}_{ab} (\D
\sigma^a_\nu) \tilde \gamma^b_\beta, \]
\[ \D \gamma^\mu_\nu \bu \D \tau^\alpha_\beta = \D \gamma^\mu_\nu
\wedge \D \tau^\alpha_\beta + {\imath\over 2} \theta^{\mu
\alpha}_{ab} \D \sigma^a_\nu \wedge \D \tilde \gamma^b_\beta, \]
\[\tau^\mu_\nu \bu \D \gamma^\alpha_\beta = \tau^\mu_\nu \D
\gamma^\alpha_\beta + {\imath\over 2}\theta^{\mu \alpha}_{ab} \tilde
\gamma^a_\nu \D \sigma^b_\beta, \quad
\D \tau^\mu_\nu \bu \gamma^\alpha_\beta = (\D \tau^\mu_\nu)
(\gamma^\alpha_\beta) + {\imath\over 2}\theta^{\mu \alpha}_{ab} (\D
\tilde \gamma^a_\nu) \sigma^b_\beta, \]
\[ \D \tau^\mu_\nu \bu \D \gamma^\alpha_\beta = \D \tau^\mu_\nu \wedge
\D \gamma^\alpha_\beta + {\imath\over 2}\theta^{\mu \alpha}_{ab} \D
\tilde \gamma^a_\nu \wedge \D \sigma^b_\beta, \]
\[ \tilde \gamma^\mu_\nu \bu \D \tau^\alpha_\beta = \tilde \gamma^\mu_\nu \D \tau^\alpha_\beta
- {\imath\over 2} \sigma^\mu_c \gamma^\alpha_d \theta^{cd}_{\nu
\beta}, \quad
\D \tilde \gamma^\mu_\nu \bu \tau^\alpha_\beta = (\D \tilde
\gamma^\mu_\nu) \tau^\alpha_\beta - {\imath\over 2} (\D
\sigma^\mu_c) \gamma^\alpha_d \theta^{cd}_{\nu \beta}, \]
\[ \D \tilde \gamma^\mu_\nu \bu \D \tau^\alpha_\beta = \D \tilde \gamma^\mu_\nu \wedge \D \tau^\alpha_\beta
- {\imath\over 2} \D \sigma^\mu_c \wedge \gamma^\alpha_d
\theta^{cd}_{\nu \beta},\]
\[ \tau^\mu_\nu \bu \D \tilde \gamma^\alpha_\beta = \tau^\mu_\nu \D \tilde
\gamma^\alpha_\beta - {\imath\over 2} \gamma^\mu_c \D
\sigma^d_\alpha \theta^{cd}_{\nu \beta}, \quad
\D \tau^\mu_\nu \bu \tilde \gamma^\alpha_\beta = (\D \tau^\mu_\nu)
\tilde \gamma^\alpha_\beta - {\imath\over 2} (\D \gamma^\mu_c)
\sigma^d_\alpha \theta^{cd}_{\nu \beta}, \]
\[ \D \tau^\mu_\nu \bu \D \tilde \gamma^\alpha_\beta = \D \tau^\mu_\nu \wedge \D \tilde
\gamma^\alpha_\beta - {\imath\over 2} \D \gamma^\mu_c \wedge \D
\sigma^d_\alpha \theta^{cd}_{\nu \beta}, \]
\[ \tau^\mu_\nu \bu \D \tau^\alpha_\beta = \tau^\mu_\nu
\D \tau^\alpha_\beta + {\imath\over 2} \theta^{\mu \alpha}_{ab}
\tilde \gamma^a_\nu \D \tilde \gamma^b_\beta - {\imath\over 2}
\gamma^\mu_c \D \gamma^\alpha_d \theta^{cd}_{\nu \beta} + {1\over 4}
\theta^{\mu \alpha}_{ab} \sigma^a_c \D \sigma^b_d \theta^{cd}_{\nu
\beta}, \]
\[ \D \tau^\mu_\nu \bu \tau^\alpha_\beta = (\D \tau^\mu_\nu)
\tau^\alpha_\beta + {\imath\over 2} \theta^{\mu \alpha}_{ab} (\D
\tilde \gamma^a_\nu) \tilde \gamma^b_\beta - {\imath\over 2} (\D
\gamma^\mu_c) \gamma^\alpha_d \theta^{cd}_{\nu \beta} + {1\over 4}
\theta^{\mu \alpha}_{ab} (\D \sigma^a_c) \sigma^b_d \theta^{cd}_{\nu
\beta}, \]
\[ \D \tau^\mu_\nu \bu \D \tau^\alpha_\beta = \D \tau^\mu_\nu \wedge
\D \tau^\alpha_\beta + {\imath\over 2} \theta^{\mu \alpha}_{ab} \D
\tilde \gamma^a_\nu \wedge \D \tilde \gamma^b_\beta - {\imath\over
2} \D \gamma^\mu_c \wedge \D \gamma^\alpha_d \theta^{cd}_{\nu \beta}
+ {1\over 4} \theta^{\mu \alpha}_{ab} \D \sigma^a_c \wedge \D
\sigma^b_d \theta^{cd}_{\nu \beta}, \]
with remaining relations undeformed.  As before we adopt the
convention that in each set of equations, the indices $\alpha,
\beta, \mu, \nu$ lie in the appropriate ranges for each $2 \times 2$
block. One may also calculate the explicit commutation relations in closed form; they
will be similar to the above but with $\theta^-$ in place of $\theta$.

Similarly, since the classical differential structures on
$\TCM,\tilde T$ are covariant under $\GL_4$, we have coactions on
their classical exterior algebras induced from the coactions on the
spaces themselves, such that $\D$ is equivariant. We can hence
covariantly twist these in the same way as the algebras themselves.
Thus $\Omega(\CC_F[\tilde T])$ has structure
\begin{equation} \label{twistor calculus}
Z^\nu\bullet \D Z^\mu=F^{\nu\mu}_{ab}Z^a \D Z^b,\quad \D
Z^\nu\bullet Z^\mu=F^{\nu\mu}_{ab}(\D Z^a )Z^b \end{equation}
\[ \D Z^\nu\bullet \D Z^\mu=F^{\nu\mu}_{ab}\D Z^a\wedge \D Z^b\]
The commutation relations are similarly
\[ Z^\nu\bullet \D Z^\mu=R^{\mu\nu}_{ab}\D Z^a\bu Z^b,\quad \D Z^\nu\bullet \D Z^\mu=- R^{\mu\nu}_{ab}\D Z^a\bu \D Z^b.\]
These formulae are essentially as for the coordinate algebra, but now with $\D$ inserted, and are (in some form) standard for braided linear spaces define by an $R$-matrix.  More explicitly,
\[ Z^\mu\bullet\D Z^\nu=Z^\mu\D Z^\nu+{\imath\over 2}\theta^{\mu\nu}_{ab}Z^a\D Z^b,\quad
\D Z^\mu\bullet Z^\nu=(\D Z^\mu) Z^\nu+{\imath\over 2}\theta^{\mu\nu}_{ab}\D Z^a Z^b\]
\[  \D Z^\mu\bullet\D Z^\nu=\D Z^\mu\wedge \D Z^\nu+{\imath\over 2}\theta^{\mu\nu}_{ab}\D Z^a\wedge \D Z^b\]
so that the $Z^3,Z^4,\D Z^3, \D Z^4$ products are undeformed. In terms of commutation relations
\[   [Z^\mu, \D Z^\nu]_\bu=\imath\theta^-{}^{\mu\nu}_{ab}Z^a\D Z^b,\quad  \{\D Z^\mu, \D Z^\nu\}_\bu=\imath\theta^-{}^{\mu\nu}_{ab}\D Z^a\wedge \D Z^b\]
where the right hand sides are for $a,b\in\{3,4\}$ and could be
written with the bullet product equally well. The only nonclassical
commutation relations here are those with $\mu,\nu\in \{1,2\}$.

Similarly, $\Omega(\CC_F[\TCM])$ has structure
\[ x^{\mu \nu} \bu \D x^{\alpha \beta} =F^{\mu\alpha}_{mn}F^{m\beta}_{ap}F^{\nu n}_{qc}F^{qy}_{bd}x^{ab} \D x^{cd},\quad \D x^{\mu \nu} \bu x^{\alpha \beta}  =F^{\mu\alpha}_{mn}F^{m\beta}_{ap}F^{\nu n}_{qc}F^{qy}_{bd}(\D x^{ab}) x^{cd}\]
\[ \D x^{\mu \nu} \bu \D x^{\alpha \beta}  =F^{\mu\alpha}_{mn}F^{m\beta}_{ap}F^{\nu n}_{qc}F^{qy}_{bd}\D x^{ab}\wedge \D x^{cd}.\]
On the affine Minkowski generators and $t$ we find explicitly:
\begin{equation} \label{affine minkowski calculus} z \bu \D \tz = z \D \tz + {\imath\over 2} \theta^{21}_{43}t \D
t, \quad
\D z \bu \tz = (\D z) \tz + {\imath\over 2} \theta^{21}_{43}(\D
t)t,\end{equation}
\[ \tz \bu \D z = \tz \D z + {\imath\over 2} \theta^{12}_{34}t \D t,
\quad
\D \tz \bu z = (\D \tz) z + {\imath\over 2} \theta^{12}_{34}(\D
t)t,\]
\[ \D z \bu \D \tz = \D z \wedge \D \tz + {\imath\over 2} \theta^{21}_{43}\D t \wedge \D
t = \D z \wedge \D \tz,\]
\[ w \bu \D \tw = w \D \tw + {\imath\over 2} \theta^{12}_{43} t \D
t, \quad
\D w \bu \tw = (\D w)\tw + {\imath\over 2} \theta^{12}_{43} (\D
t)t,\]
\[ \tw \bu \D w = \tw \D w + {\imath\over 2} \theta^{21}_{34} t \D
t, \quad
\D \tw \bu w = (\D \tw)w + {\imath\over 2} \theta^{21}_{34} (\D
t)t,\]
\[ \D w \bu \D \tw = \D w \wedge \D \tw + {\imath\over 2} \theta^{12}_{43} \D t \wedge \D
t = \D w \wedge \D \tw,\]
with other relations amongst these generators undeformed, as are the
relations involving $\D t$, whence we may equally use the $\bu$
product in terms which involve $t,\D t$. The relations in the
calculus involving $s, \D s$ are more complicated. We write just the
final commutation relations for these:
\[ [s,\D z]_\bu = -{\imath\over 2}\theta^{12}_{34} z\D t +
{\imath\over 2}\theta{}^{21}_{43} t \D z -{\imath \over
2}\theta^{12}_{44}\tw \D t + {\imath\over 2}\theta^{21}_{44}t\D
\tw;\]
\[ [\D s, z]_\bu = -{\imath\over 2}\theta^{12}_{34}(\D z) t +
{\imath\over 2}\theta{}^{21}_{43}(\D t) z -{\imath \over
2}\theta^{12}_{44}(\D \tw) t + {\imath\over 2}\theta^{21}_{44}(\D t)
\tw;\]
\[ [\D s,\D z]_\bu = -{\imath\over 2}\theta^{12}_{34}\D z\wedge \D t +
{\imath\over 2}\theta{}^{21}_{43}\D t\wedge \D z -{\imath \over
2}\theta^{12}_{44}\D \tw\wedge \D t + {\imath\over
2}\theta^{21}_{44}\D t\wedge\D \tw;\]
\[ [s,\D \tz]_\bu = -{\imath\over 2} \theta^{21}_{33} w\D t + {\imath\over 2}\theta^{12}_{33} t \D w + {\imath\over 2}\theta^{12}_{43} \tz \D t - {\imath\over 2}\theta^{12}_{34} t \D \tz;\]
\[ [\D s,\tz]_\bu = -{\imath\over 2} \theta^{21}_{33} (\D w) t + {\imath\over 2}\theta^{12}_{33}(\D t) w + {\imath\over 2}\theta^{12}_{43}(\D \tz) t - {\imath\over 2}\theta^{12}_{34} (\D t) \tz;\]
\[ [\D s,\D \tz]_\bu = -{\imath\over 2} \theta^{21}_{33} \D w\wedge \D t + {\imath\over 2}\theta^{12}_{33} \D t \wedge\D w + {\imath\over 2}\theta^{12}_{43} \D \tz\wedge \D t - {\imath\over 2}\theta^{12}_{34} \D t\wedge \D \tz;\]
\[ [s,\D w]_\bu = -{\imath\over 2} \theta^{21}_{34} w\D t + {\imath\over 2}\theta^{12}_{43} t \D w + {\imath\over 2}\theta^{21}_{44} \tz \D t - {\imath\over 2}\theta^{12}_{44} t \D \tz;\]
\[ [\D s,w]_\bu = -{\imath\over 2} \theta^{21}_{34} (\D w) t + {\imath\over 2}\theta^{12}_{43}(\D t) w + {\imath\over 2}\theta^{12}_{44}(\D \tz) t - {\imath\over 2}\theta^{12}_{44} (\D t) \tz;\]
\[ [\D s,\D w]_\bu = -{\imath\over 2} \theta^{21}_{34} \D w\wedge \D t + {\imath\over 2}\theta^{12}_{43} \D t \wedge\D w + {\imath\over 2}\theta^{21}_{44} \D \tz \wedge\D t - {\imath\over 2}\theta^{12}_{44} \D t \wedge\D \tz;\]
\[ [s,\D \tw]_\bu = {\imath\over 2} \theta^{12}_{33} z\D t - {\imath\over 2}\theta^{21}_{33} t \D z - {\imath\over 2}\theta^{12}_{43} \tw \D t + {\imath\over 2}\theta^{21}_{34} t \D \tw;\]
\[ [\D s, \tw]_\bu = {\imath\over 2} \theta^{12}_{33} (\D z) t - {\imath\over 2}\theta^{21}_{33} (\D t) z - {\imath\over 2}\theta^{12}_{43} (\D \tw) t + {\imath\over 2}\theta^{21}_{34} (\D t) \tw;\]
\[ [\D s,\D \tw]_\bu = {\imath\over 2} \theta^{12}_{33} \D z\wedge \D t - {\imath\over 2}\theta^{21}_{33} \D t\wedge \D z - {\imath\over 2}\theta^{12}_{43} \D \tw\wedge \D t + {\imath\over 2}\theta^{21}_{34} \D t \wedge \D \tw.\]

The calculus of the correspondence space algebra $\C_F[\tilde \F]$
is generated by $\D Z^\mu$, $\D s$, $\D t$, $\D z$, $\D \tz$, $\D
w$, $\D \tw$ with twisted relations given by (\ref{correspondence
space relations}) with $\D$ inserted where appropriate, as well as
the relations given by differentiating (\ref{twisted alpha planes})
using the Leibniz rule. Explicitly we have
\begin{equation}
\label{twisted correspondence space calculus} [z,\D Z^1]_\bu =
\imath \theta^-{}^{21}_{43} t\D Z^3 + \imath \theta^-{}^{21}_{44} t
\D Z^4, \qquad [z,\D Z^2]_\bu = \imath \theta^-{}^{22}_{43} t \D
Z^3,\end{equation}
\[ [\D z,Z^1]_\bu = \imath \theta^-{}^{21}_{43} (\D t)Z^3 + \imath
\theta^-{}^{21}_{44} (\D t)Z^4, \qquad [\D z,Z^2]_\bu = \imath
\theta^-{}^{22}_{43} (\D t)Z^3\]
\[[\D z,\D Z^1]_\bu = \imath \theta^-{}^{21}_{43} \D t\wedge \D Z^3 + \imath
\theta^-{}^{21}_{44} \D t\wedge \D Z^4, \qquad [\D z, \D Z^2]_\bu =
\imath \theta^-{}^{22}_{43} \D t\wedge \D Z^3,\]
for example, where we may equally use the $\bu$ product on the
right-hand sides.  Moreover,
\begin{eqnarray}
\nonumber && (\D \tz )Z^3 + (\D w) Z^4 -(\D t) Z^1 +  \tz \D Z^3 + w \D Z^4 -t \D Z^1=0,\\
\nonumber && (\D \tw) Z^3 + (\D z) Z^4 - (\D t) Z^2 + \tw \D Z^3 + z \D Z^4 - t \D Z^2=0,\\
\nonumber && (\D s) Z^3+(\D w) Z^2+ -(\D z) Z^1 + {\imath \over
2}((\theta^{12}_{43}- \theta^{21}_{43})(\D t)Z^3 +
\theta^-{}^{12}_{44}(\D t)Z^4) \\
\nonumber && \qquad + s \D Z^3+w \D Z^2+ -z \D Z^1 + {\imath \over
2}((\theta^{12}_{43}-
\theta^{21}_{43})t\D Z^3 + \theta^-{}^{12}_{44}t\D Z^4)=0, \\
\nonumber && (\D s) Z^4-(\D \tz) Z^2+(\D \tw) Z^1 + {\imath \over 2}
(\theta^-{}^{21}_{33} (\D t)Z^3 + (\theta^{21}_{34} -
\theta^{12}_{34})(\D t)Z^4) \\
\nonumber && \qquad + s \D Z^4-\tz \D Z^2+\tw \D Z^1 + {\imath \over
2} (\theta^-{}^{21}_{33} t\D Z^3 + (\theta^{21}_{34} -
\theta^{12}_{34})t\D Z^4)=0.
\end{eqnarray}
Of course, there relations may be written much more compactly using
the $\bu$ product (as explained, the twisted relations are the same
as the classical case, save for replacing the old product by the
twisted one).  The other relations are obtained similarly, hence we
refrain from writing them out explicitly, although we stress that
$\D Z^3, \D Z^4, \D t$ are central in the calculus.

\section{Quantisation by cotwists in the $\SU_n$  $*$-algebra version}
\label{Quantisation by cotwists in the *-algebra version}

Our second setting is to work with twistor space and space-time as
real manifolds in a $*$-algebra context. To this end we gave a
`projector' description of our spaces $\CP^3,\hCM,\F$ as well as
their realisation as $\SU_4$ homogeneous spaces. We also show that
this setting too quantises nicely by a cochain twist. This approach
is directly compatible with $C^*$-algebra methods, although we shall
not perform the $C^*$-algebra completions here.  We shall however
simultaneously quantise all real manifolds defined by $n\times
n$-matrices of generators with $\SU_n$-covariant conditions, a class
which (as we have seen) includes all partial flag varieties based on
$\C^n$, with $n=4$ the case relevant for the paper.

The main new ingredient is that in the general theory of cotwisting,
we should add that $H$ is a Hopf $*$-algebra in the sense that
$\Delta$ is a $*$-algebra map and $(S\circ *)^2={\rm id}$, and that
the cocycle is real in the sense \cite{Ma:book}
\[ \overline{F(h,g)}=F((S^2g)^*,(S^2h)^*).\]
In this case the twisted Hopf algebra acquires a new $*$-structure
\[ h^{*_F}=\sum \overline{V^{-1}(S^{-1}h\o)}(h\t)^*\overline{V(S^{-1}h\thr)},\quad V(h)=U^{-1}(h\o)U(S^{-1}h\t).\]
(We note the small correction to the formula stated in
\cite{Ma:book}). Also, if $A$ is a left comodule algebra and a
$*$-algebra, we require the coaction $\Delta_L$ to be a $*$-algebra
map. Then $A_F$ has a new $*$-structure
\[ a^{*_F}=\overline{V^{-1}(S^{-1}a\o)}(a\t)^*.\]

The twisting subgroup in our application will be $\S(\U(1)\times
\U(1)\times \U(1)\times \U(1))$ since this is contained in all
relevant subgroups of $\SU_4$, or rather the larger group $\SU_n$
with twisting subgroup appearing in the Hopf $*$-algebra picture as
\[ H=\C[\S(\U(1)^n)]=\C[t_\mu,t_\mu^{-1};\ \mu=1,\cdots, n]/\la t_1\cdots t_n=1\ra\]
\[ \quad t_\mu^*=t_\mu^{-1},\quad \Delta t_\mu=t_\mu\tens t_\mu\quad \eps t_\mu=1,\quad St_\mu=t_\mu^{-1}.\]
This is also the group algebra of the Abelian group
$\Z^n/\Z(1,1,\cdots,1)$ (where the vector here has $n$ entries, all
equal to $1$). We define a basis of $H$ by
\[ t^{\vec a}=t_1^{a_1}\cdots t_n^{a_n},\quad \vec a\in \Z^n/\Z(1,1,\cdots,1).\]
Note that we could of course eliminate one of the $\U(1)$ factors
and identify the group with $\U(1)^{n-1}$ and the dual group with
$\Z^{n-1}$, and this would be entirely equivalent in what follows
but not canonical. We prefer to keep manifest the natural inclusion
on the diagonal of $\SU_n$. This inclusion appears now as the
$*$-Hopf algebra surjection
\[ \pi:\C[\SU_n]\to H,\quad \pi(a^\mu_\nu)=\delta^\mu_\nu t_\mu.\]

Next, we define a cocycle $F:H\tens H\to \C$ by
\[ F(t^{\vec a},t^{\vec b})=e^{\imath \vec a\cdot\theta\cdot \vec b}\]
where $\theta\in M_n(\C)$ is any matrix for which every row and
every column adds up to zero (so that $(1,1,\cdots, 1)$ is in the
null space from either side).  Such a matrix is fully determined by
an arbitrary choice of (say) lower $(n-1)\times (n-1)$ diagonal
block. Thus the data here is an arbitrary $(n-1)\times (n-1)$ matrix
just as we would have if we had eliminated $t_1$ in the first place.
The reality condition and the functional $V$ work out as
\[\theta^\dagger=-\theta,\quad U(t^{\vec a})=e^{-\imath \vec a\theta\vec a},\quad V(t^{\vec a})=1.\]
The latter means that the $*$-structures do not deform.

We now pull back this cocycle under $\pi$ to a cocycle on
$\C[\SU_n]$ with matrix
\[ \mathbb{F} = F(a^\mu_\nu, a^\alpha_\beta) = \delta^\mu_\nu \delta ^\alpha_\beta F(t_\nu,t_\beta) = \delta^\mu_\nu \delta ^\alpha_\beta e^{\imath \theta_{\nu \beta}}.\]
We then use these new $F$-matrices in place of those in Sections~4,5
since the general formulae in terms of $F$-matrices are identical,
being determined by the coactions and coproducts, with the
additional twist of the $*$-operation computed in the same way. Thus
we find easily that $ \C_F[\SU_n]$ has the deformed product and
antipode (and undeformed $*$-structure):
\[ a^\mu_\nu \bullet a^\alpha_\beta = e^{\imath(\theta_{\mu\alpha} -
\theta_{\nu\beta})}a^\mu_\nu a^\alpha_\beta =
e^{\imath(\theta_{\mu\alpha}-\theta_{\nu\beta}-\theta_{\alpha\mu}+\theta_{\beta\nu})}
a^\alpha_\beta\bullet a^\mu_\nu, \quad
S_Fa^\mu_\nu=e^{-\imath(\theta_{\mu\mu}-\theta_{\nu\nu})}Sa^\mu_\nu\]
which means that $\C_F[\SU_n]$ has a compact form if we use new
generators
\[ \hat a^\mu_\nu=e^{{\imath\over 2}(\theta_{\mu\mu}-\theta_{\nu\nu})}a^\mu_\nu\]
in the sense
\[\Delta \hat a^\mu_\nu=\hat a^\mu_\alpha\tens \hat a^\alpha_\nu,
\quad S_F\hat a^\mu_\alpha\bullet \hat
a^\alpha_\nu=\delta^\mu_\nu=\hat a^\mu_\alpha\bullet  S_F\hat
a^\alpha_\nu,\quad \hat a^\mu_\nu{}^*=S_F \hat a^\nu_\mu.\]
They have the same form of commutation relations as the $a^\mu_\nu$
with respect to the $\bullet$ product. Note that for a $C^*$-algebra
treatment we will certainly want the commutation relations in `Weyl
form' with a purely phase factor and hence $\theta$ to be
real-valued, hence antisymmetric and hence with zeros on the
diagonal. So in this case natural case there will be no difference
between the $\hat a$ and the $a$ generators.

Similarly we find by comodule cotwist that
\[ Z^\mu\bullet Z^\nu=e^{\imath\theta_{\mu\nu}}Z^\mu Z^\nu,
\quad Z^\mu{}^*\bullet Z^\nu=e^{-\imath\theta_{\mu\nu}}Z^\mu{}^*
Z^\nu,\quad Z^\mu{}\bullet Z^\nu{}^*=e^{-\imath\theta_{\mu\nu}}Z^\mu
Z^\nu{}^*,\]
or directly from the unitary transformation of any projectors
\[ P^\mu_\nu\bullet P^\alpha_\beta=F(a^\mu_a S a^c_\nu,a^\alpha_b a^d_\beta)P^a_c P^b_d
=F(t_\mu t_\nu^{-1},t_\alpha t_\beta^{-1})P^\mu_\nu
P^\alpha_\beta=e^{\imath(\theta_{\mu
\alpha}-\theta_{\nu\alpha}+\theta_{\nu\beta}-\theta_{\mu\beta})}P^\mu_\nu
P^\alpha_\beta.\]
The commutation relations for the entries of $Z, P$ respectively
have the same form as the deformation relations but with
$\theta_{\mu\nu}$ replaced by $2\theta^-_{\mu\nu}\equiv
\theta_{\mu\nu}-\theta_{\nu\mu}$.

The $P^\mu_\nu$ are no longer projectors with respect to the $\bullet$ product but the new generators
\[ \hat P^\mu_\nu=e^{-\imath\theta_{\mu\nu}+{\imath\over 2}(\theta_{\mu\mu}+\theta_{\nu\nu})}P^\mu_\nu\]
are. They enjoy the same commutation relations as the $P^\mu_\nu$ with respect to the bullet product.  Moreover
\[ \Tr ~\hat P=\Tr ~P,\quad \hat P^\mu_\nu{}^*=\hat P^\nu_\mu.\]
Thus we see for example that  $\C_F[\CP^3]$ has quantised
commutation relations for the matrix entries of generators $\hat Q$
with further relations $\Tr ~\hat Q=1$ and $*$-structure $\hat
Q^\dagger=\hat Q$, i.e. exactly the same form for the
matrix-$\bullet$ relations as in the classical case.  Applying these
computations but now with $\Tr~ P=2$ for the matrix generator $P$,
we obtain $\C_F[\hCM]$ in exactly the same way but with  $\Tr~\hat
P=2$, so that in the projector picture we cover both cases at the
same time but with different values for the trace.

For $\C_F[\F]$ we have projectors $\hat Q$ of trace $1$ and $\hat P$
of trace $2$. Their products are deformed in the same manner as the
$P$-$P$ relations, leading to
\[ \hat P^\mu_\nu\bullet \hat Q^\alpha_\beta =e^{2\imath(\theta^-_{\mu \alpha}-\theta^-_{\nu\alpha}+\theta^-_{\nu\beta}-\theta^-_{\mu\beta})}\hat Q^\alpha_\beta\bullet \hat P^\mu_\nu\]
for the quantised commutation relations between entries. Moreover,
$\hat P\bullet \hat Q=\hat Q=\hat Q\bullet\hat P$ by a similar
computation as for $\hat P^2=\hat P$. In particular,  we see that
$\hat Q\in M_4(\C_F[\CP^3])$ and $\hat P\in M_4(\C_F[\hCM])$ are
projectors which define tautological quantum vector bundles over
these quantum spaces and their pull-backs to $\C_F[\F]$.

Rather than proving all these facts for each algebra, let us prove
them for the quantisation of any real manifold $X\subset M_n(\C)^r$
defined as the set of $r$-tuples of matrices $P_1,\cdots P_r$
obeying relations defined by the operations of: (a) matrix product;
(b) trace; (c) the $(\ )^\dagger$ operation of Hermitian
conjugation. We say that $X$ is defined by `matrix relations'.
Clearly any such $X$ has on it an action of $\SU_n$ acting by
conjugation. We define $\C[X]$ to be the (possibly $*$-) algebra
defined by treating the matrix entries $P_i{}^\mu_\nu$ as polynomial
generators, the matrix relations as relations  in the algebra, and
$P_i^\dagger$ (when specified) as a definition of
$P_i{}^\nu_\mu{}^*$. The coaction of $\C[\SU_n]$ is
\[ \Delta_L P_i{}^\mu_\nu=a^\mu_a S a^c_\nu \tens P_i{}^a_c.\]
We have already seen several examples of such coordinate algebras with matrix relations.

\begin{prop}
Let $\C[X]$ be a  $*$-algebra defined by `matrix relations' among
matrices of generators $P_i$, $i=1,\cdots,r$. Its quantisation
$\C_F[X]$ by cocycle cotwist using the cocycle above is the free
associative algebra with matrices of generators $\hat P_i$ modulo
the commutation relations
\[ \hat P_i{}^\mu_\nu\bullet \hat P_j{}^\alpha_\beta =e^{2\imath(\theta^-_{\mu \alpha}-\theta^-_{\nu\alpha}+\theta^-_{\nu\beta}-\theta^-_{\mu\beta})}\hat P_j{}^\alpha_\beta\bullet \hat P_i{}^\mu_\nu\]
and the matrix relations of $\C[X]$ with $P_i$ replaced by $\hat P_i$.
\end{prop}
\proof All the $P_i$ have the same coaction, hence for the deformed
product for any $P,Q\in \{P_1,\cdots,P_r\}$ we have
\[ P^\mu_\nu\bullet Q^\alpha_\beta=e^{\imath(\theta_{\mu \alpha}-\theta_{\nu\alpha}+\theta_{\nu\beta}-\theta_{\mu\beta})}P^\mu_\nu Q^\alpha_\beta\]
by the same computation as for the $P$-$P$ relations above. This
implies the commutation relations stated for the entries of $P_i$
and hence of the $\hat P_i$. Again motivated by the example we
define
\[ \hat P_i{}^\mu_\nu=e^{-\imath\theta_{\mu\nu}+{\imath\over 2}(\theta_{\mu\mu}+\theta_{\nu\nu})}P_i{}^\mu_\nu\]
and verify that
\[\hat P\bullet\hat Q= \widehat{(PQ)},\quad (\hat P)^\dagger=\widehat{(P^\dagger)},\quad \Tr~\hat P=\Tr ~P\]
for any $P,Q$ taken from our collection. Thus
\begin{eqnarray*}(\hat P\bullet\hat Q)^\mu_\nu &=&\hat P^\mu_\alpha\bullet \hat P^\alpha_\nu=e^{-\imath\theta_{\mu\alpha}+{\imath\over 2}(\theta_{\mu\mu}+\theta_{\alpha\alpha})-\imath\theta_{\alpha\nu}+{\imath\over 2}(\theta_{\alpha\alpha}+\theta_{\nu\nu})}P^\mu_\alpha\bullet P^\alpha_\nu\\
&=&e^{-\imath\theta_{\mu\alpha}+{\imath\over 2}(\theta_{\mu\mu}+\theta_{\alpha\alpha})-\imath\theta_{\alpha\nu}+{\imath\over 2}(\theta_{\alpha\alpha}+\theta_{\nu\nu})}e^{\imath(\theta_{\mu \alpha}-\theta_{\alpha\alpha}+\theta_{\alpha\nu}-\theta_{\mu\nu})}P^\mu_\alpha Q^\alpha_\nu \\
&=&e^{-\imath\theta_{\mu\nu}+{\imath\over 2}(\theta_{\mu\mu}+\theta_{\nu\nu})}P^\mu_\alpha Q^\alpha_\nu=\widehat{(PQ)}^\mu_\nu\\
\hat
P^\nu_\mu{}^*&=&e^{\imath\overline{\theta_{\nu\mu}}-{\imath\over
2}(\overline{\theta_{\nu\nu}}+\overline{\theta_{\mu\mu}})}P^\nu_\mu{}^*=e^{-\imath\theta_{\mu\nu}+{\imath\over
2}(\theta_{\mu\mu}+\theta_{\nu\nu})}P^\dagger{}^\mu_\nu=\widehat{(P^\dagger)}^\mu_\nu.\end{eqnarray*}
The proof for the trace is immediate from the definition. We also
note that $\hat\delta^\mu_\nu=\delta^\mu_\nu$ as the quantisation of
the constant identity projector (the identity for matrix
multiplication).
\endproof

For example, we now have the quantisation of all flag varieties
$\C_F[\F_{k_1,\cdots,k_r}(\C^n)]$ with projectors $\hat P_i$ having
this new form of commutation relations for their matrix entries, but
with matrix-$\bullet$ products having the same form as in the
classical case given in Section~3. Again, the $\hat P_i$ define $r$
tautological projectors with values in the quantum algebra and hence
$r$ tautological classes in the noncommutative $K$-theory, strictly
quantising the commutative situation.

Let us also make some immediate observations from the relations in
the proposition. We see that diagonal elements $\hat P_i{}^\mu_\mu$
(no sum) are central. We also see that $\hat P_i{}^\mu_\nu$ and
$\hat P_i{}^\nu_\mu=(\hat P_i{}^\mu_\nu)^*$ commute (so all matrix
entry generators are normal in the *-algebra sense). On the other
hand nontrivial commutation relations arise when three of the four
indices are different and that if we take the adjoint of a generator
on both sides of a commutation relation, we should also invert the
commutation factor. This means that elements of the form $(\hat
P_i{}^\mu_\nu)^*\hat P_i{}^\mu_\nu$ (no summation) are always
central.

These observations  mean that $\C_F[\CP^1]$ is necessarily undeformed
in the new generators. For a nontrivial deformation the smallest
example is then $\C_F[\CP^2]$. Writing its matrix generator as
\[ \hat Q=\begin{pmatrix} a& x & y\\ x^* & b & z\\ y^* & z^* & c\end{pmatrix}\]
one has $a,b,c$ are self-adjoint and central with $a+b+c=1$ and
\[ xy=e^{\imath\theta}yx,\quad yz=e^{\imath\theta}zy,\quad xz=e^{\imath\theta}zx,\quad \theta=2(\theta^-_{12}+\theta^-_{23}+\theta^-_{31})\]
and the projection relations exactly as stated in Proposition~3.4
(whose statement and proof assumed only that $a,b,c,x^*x,y^*y,z^*z$
are central; no other commutation relations were actually needed).
Also note that since $a,b,c$ are central it is natural to set them
to constants even in the quantum case. In the quotient
$a=b=c={1\over 3}$ we can define $U=3x$, $W=3y$, $V=3z$. Then we
have the algebra
\[ UV=W,\quad U^*=U^{-1},\quad V^*=V^{-1},\quad W^*=W^{-1}\]
\[ UW=e^{\imath\theta}WU,\quad WV=e^{\imath\theta}VW,\quad UV=e^{\imath\theta}VU.\]
Actually this is just the usual noncommutative torus
$\C_\theta[S^1\times S^1]$ with $W$ defined by the above relations
and no additional constraints. Similarly in general, for any actual
values with $b,c>0$ and $b+c<1$ we will have the same result but
with different rescaling factors for $U,V,W$, i.e. again
noncommutative tori as quantum versions of a family of inclusions
$S^1\times S^1\subset\CP^2$. We can consider this family as a
quantum analogue of $\C^*\times\C^*\subset\CP^2$. This conforms to
our expectation of  $\C_F[\CP^2]$ as a `quantum toric variety'.
Moreover, by arguments analogous to the classical case in
Section~3.2 we can view the localisation $\C_F[\CP^2][a^{-1}]$ as a
punctured quantum $S^4$ with generators $x,y,a,a^{-1}$ and the
relations
\[ xy=e^{\imath\theta}yx,\quad x^*x+y^*y=a(1-a),\]
with $a$ central.

One can also check that the cotriangular Hopf $*$-algebra
$\C_F[\SU_n]$ coacts on $\C_F[X]$ now as the quantum version of our
classical coactions, as is required by the general theory, namely
that these quantised spaces may be realised as quantum homogeneous
spaces if one wishes. Again from general theory, these spaces are
$\Psi$-commutative with respect to the induced involutive braiding
built from $\theta^-$ and appearing in the commutation relations for
the matrix entries.

Finally, the differential calculi are twisted by the same methods as
in Section~5.  The formulae are similar to the deformation of the
coordinate algebras, with the insertion of $\D$ just as before. Thus
$\Omega(\C_F[\SU_n])$ has structure
\[ a^\mu_\nu\bu \D a^\alpha_\beta=e^{\imath(\theta_{\mu\alpha}-\theta_{\nu\beta})}a^\mu_\nu\D a^\alpha_\beta,\quad \D a^\mu_\nu\bu a^\alpha_\beta=e^{\imath(\theta_{\mu\alpha}-\theta_{\nu\beta})}\D a^\mu_\nu a^\alpha_\beta,\]
\[ \D a^\mu_\nu\bu \D a^\alpha_\beta=e^{\imath(\theta_{\mu\alpha}-\theta_{\nu\beta})}\D a^\mu_\nu\wedge \D a^\alpha_\beta\]
and commutation relations
\[\hat a^\mu_\nu\bu \D \hat a^\alpha_\beta=e^{2\imath(\theta^-_{\mu\alpha}-\theta^-_{\nu\beta})}\hat a^\mu_\nu\bu \D \hat a^\alpha_\beta,\quad \D \hat a^\mu_\nu\bu \D \hat a^\alpha_\beta=-e^{2\imath(\theta^-_{\mu\alpha}-\theta^-_{\nu\beta})}\D\hat a^\mu_\nu\bullet \D \hat a^\alpha_\beta.\]
Similarly, for the quantisation $\C_F[X]$ above of an algebra with
matrix relations, $\Omega(\C_F[X])$ has  commutation relations
\[  \hat P_i{}^\mu_\nu\bullet \D \hat P_j{}^\alpha_\beta =e^{2\imath(\theta^-_{\mu \alpha}-\theta^-_{\nu\alpha}+\theta^-_{\nu\beta}-\theta^-_{\mu\beta})}\D \hat P_j{}^\alpha_\beta\bullet \hat P_i{}^\mu_\nu,\]
\[ \D \hat P_i{}^\mu_\nu\bullet \D \hat P_j{}^\alpha_\beta =- e^{2\imath(\theta^-_{\mu \alpha}-\theta^-_{\nu\alpha}+\theta^-_{\nu\beta}-\theta^-_{\mu\beta})}\D \hat P_j{}^\alpha_\beta\bullet \D \hat P_i{}^\mu_\nu.\]

\goodbreak

\subsection{Quantum $\C_F[\hCM]$, $\C_F[S^4]$ and the quantum instanton}

In this section we specialise the above general theory to $\C[\hCM]$
and its cocycle twist quantisation. We also find a natural
one-parameter family of the $\theta_{\mu\nu}$-parameters for which
$\C_F[\hCM]$ has a $*$-algebra quotient $\C_F[S^4]$. We find that
this recovers the $S^4_\theta$ previously introduced by Connes and
Landi \cite{cl:id} and that the quantum tautological bundle (as a
projective module) on $\C_F[\hCM]$ pulls back in this case to a
bundle with Grassmann connection equal to the noncommutative
instanton found in \cite{CS}.  This is very different from the
approach in \cite{CS}.

We start with some notations.  Since here $n=4$, $\theta_{\mu\nu}$
is a $4\times 4$ matrix with all rows and columns summing to zero.
For convenience we limit ourselves to the case where $\theta$ is
real and hence (see above) antisymmetric (only the antisymmetric
part enters into the commutation relations so this is no real loss).
As a result it is equivalent to giving a $3\times 3$ real
antisymmetric matrix, i.e. it has within it only three independent
parameters.

\begin{lemma}
\[ \theta_A=\theta_{12}+\theta_{23}+\theta_{31},\quad \theta_B=\theta_{23}+\theta_{34}+\theta_{42},\quad \theta=\theta_{13}-\theta_{23}+\theta_{24}-\theta_{14}\]
determine any antisymmetric theta completely.
\end{lemma}
\begin{proof} We write out the three independent equations $\sum_j \theta_{ij}=0$ for $j=1,2,3$  as
\[ \theta_{12}+\theta_{13}+\theta_{14}=0,\quad -\theta_{12}+\theta_{23}+\theta_{24}=0,\quad -\theta_{13}-\theta_{23}+\theta_{34}=0\]
and sum the first two, and sum the last two to give:
\[ \theta_{13}+\theta_{23}+\theta_{14}+\theta_{24}=0,\quad \theta_{24}+\theta_{34}-\theta_{12}-\theta_{13}=0.\]
Adding the first to $\theta$ equations tells us that $\theta=2(\theta_{13}+\theta_{24})$, while adding the second
to $\theta_A-\theta_B$ tells us that $\theta_A-\theta_B=2(\theta_{24}-\theta_{13})$. Hence, knowing $(\theta,\theta_A-\theta_B)$ is equivalent to knowing $\theta_{13},\theta_{24}$. Finally,
\[ \theta_A+\theta_B=\theta_{12}+2\theta_{23}+\theta_{34}-\theta_{13}-\theta_{24}=\theta_{12}+3\theta_{23}-\theta_{24}=4(\theta_{12}-\theta_{24})\]
using the third of our original three equations to identify $\theta_{23}$ and then the second of our original three equations to replace it. Hence knowing $\theta,\theta_A-\theta_B$, we see that knowing $\theta_A+\theta_B$ is equivalent to knowing $\theta_{12}$. The remaining $\theta_{14},\theta_{23},\theta_{34}$ are determined from our original three equations. This completes the proof, which also provides the explicit formulae:
\[ \theta_{24}={1\over 4}(\theta+\theta_A-\theta_B),\quad \theta_{13}={1\over 4}(\theta-\theta_A+\theta_B),\quad  \theta_{12}={1\over 4}\theta+{1\over 2}\theta_A\]
\[ \theta_{14}=-{1\over 2}\theta-{1\over 4}\theta_A-{1\over 4}\theta_B,\quad \theta_{23}={1\over 4}\theta_A+{1\over 4}\theta_B,\quad\theta_{34}={1\over 4}\theta+{1\over 2}\theta_B\]
\end{proof}

We are now ready to compute the commutation relations between the entries of
\[ \hat P=\begin{pmatrix}\hat A&\hat B\\ \hat B^\dagger & \hat D \end{pmatrix},\quad \hat B=\begin{pmatrix}z & \tilde w\\ w & \tilde z\end{pmatrix}\]
as in (\ref{twistorB}), except that now the matrix entry generators
are for the quantum algebra $\C_F[\hCM]$ (we omit their hats). From
the general remarks after Proposition~6.1 we know that the
generators along the diagonal, i.e. $a,\alpha_3,\delta_3$, are
central (and self-conjugate under $*$). Also from general remarks we
know that {\rm all} matrix entry generators $\hat P$ are normal
(they commute with their own conjugate under $*$).  Moreover, it is
easy to see that if $xy=\lambda yx$ is a commutation relation
between any two matrix entries then so is $xy^*=\bar\lambda y^*x$,
again due to the form of the factors in Proposition~6.1.

\begin{prop} The nontrivial commutation relations of $\C_F[\hCM]$ are
\[ \alpha\delta=e^{2\imath\theta}\delta\alpha,\]
\[ \alpha z=e^{2\imath\theta_A}z\alpha,\quad \alpha \tilde w=e^{2\imath(\theta+\theta_A)} \tilde w\alpha,\quad \alpha w=e^{2\imath\theta_A} w\alpha,\quad \alpha \tilde z=e^{2\imath(\theta+\theta_A)}\tilde z\alpha,\]
\[ \delta z=e^{-2\imath(\theta+\theta_B)}z\delta,\quad \delta w=e^{-2\imath\theta_B} w\delta,\quad  \delta \tilde w=e^{-2\imath(\theta+\theta_B)}\tilde w\delta,\quad \delta \tilde z=e^{-2\imath\theta_B}\tilde z\delta,\]
\[ z \tilde w=e^{2\imath (\theta+\theta_B)}\tilde w z,\quad z w=e^{2\imath \theta_A }wz,\quad  z \tilde z=e^{2\imath (\theta+\theta_A+\theta_B)}\tilde z z\] \[ \tilde w w=e^{2\imath(\theta_A- \theta_B)}w \tilde w,\quad  \tilde w \tilde z=e^{2\imath (\theta+\theta_A)}\tilde z \tilde w,\quad w \tilde z=e^{2\imath \theta_B}\tilde z w,\]
and similar relations with inverse coefficient when a generator is replaced by its conjugate under $*$. The further (projector) relations of $\C_F[\hCM]$ are exactly the same as in stated in Corollary~3.3 except for the last two auxiliary relations:
\[  (\alpha_3+\delta_3)\begin{pmatrix}z\\ \tilde z\end{pmatrix}=\begin{pmatrix}-\alpha & -e^{-2\imath(\theta+\theta_B)} \delta^*\\ e^{2\imath\theta_B}\delta & \alpha^*\end{pmatrix}\begin{pmatrix}w\\ \tilde w\end{pmatrix},\]
\[ (\alpha_3-\delta_3)\begin{pmatrix}w\\ \tilde w\end{pmatrix}=\begin{pmatrix}-\alpha^* & -e^{-2\imath\theta_B}\delta^*\\ e^{2\imath(\theta+\theta_B)}\delta & \alpha\end{pmatrix}\begin{pmatrix}z\\ \tilde z\end{pmatrix}.\]
\end{prop}
\begin{proof} Here the product is the twisted $\bullet$ product which we do not denote explicitly.  We use the commutation relations in Proposition~6.1, computing the various instances of
\[ \theta_{ijkl}=\theta_{ik}-\theta_{jk}+\theta_{jl}-\theta_{il}\]
in terms of the combinations in Lemma~6.2. This gives
\[ \alpha \hat B^i{}_1=e^{2\imath\theta_A}\hat B^i{}_1\alpha,\quad \alpha \hat B^i{}_2=e^{2\imath(\theta+\theta_A)}\hat B^i{}_2\alpha,\quad i=1,2\]
\[  \delta \hat B^1{}_i=e^{-2\imath(\theta+\theta_B)}\hat B^1{}_i\delta,\quad \delta \hat B^2{}_i=e^{-2\imath\theta_B}\hat B^2{}_i\delta,\quad i=1,2\]
\[ \hat B^1{}_1 \hat B^1{}_2=e^{2\imath (\theta+\theta_B)}\hat B^1{}_2 \hat B^1{}_1,\quad \hat B^1{}_1 \hat B^2{}_1=e^{2\imath \theta_A }\hat B^2{}_1 \hat B^1{}_1,\]
\[ \hat B^1{}_1\hat B^2{}_2=e^{2\imath (\theta+\theta_A+\theta_B)}\hat B^2{}_2 \hat B^1{}_1,\quad \hat B^1{}_2 \hat B^2{}_1=e^{2\imath (\theta_A-\theta_B)}\hat B^2{}_1 \hat B^1{}_2,\]
\[ \hat B^1{}_2 \hat B^2{}_2=e^{2\imath (\theta+\theta_A)}\hat B^2{}_2 \hat B^1{}_2,\quad \hat B^2{}_1 \hat B^2{}_2=e^{2\imath \theta_B}\hat B^2{}_2 \hat B^2{}_1,\]
which we write out more explicitly as stated. As explained above,
the diagonal elements of $A,D$ are central and for general reasons
the conjugate relations are as stated. Finally, we explicitly
recompute the content of the noncommutative versions of
(\ref{hcma})-(\ref{hcmb}) to find the relations required for $\hat
P$ to be a projector (this is equivalent to computing the bullet
product from the classical relations).  The $a,\alpha,\alpha^*$ form
a commutative subalgebra, as do $a,\delta,\delta^*$, so the
calculations of $A(1-A)$ and $(1-D)D$ are not affected. We can
compute $BB^\dagger$ without any commutativity assumptions, and in
fact we stated all results from (\ref{hcma}) in Corollary~3.3
carefully so as to still be correct without such assumptions. Being
similarly careful for (\ref{hcmb}) gives the remaining two auxiliary
equations (without any commutativity assumptions) as
\[ (\alpha_3+\delta_3)\begin{pmatrix} z \\ \tilde z\end{pmatrix}=\begin{pmatrix} -\tilde w\delta^*-\alpha w\\ w\delta+\alpha^*\tilde w\end{pmatrix},\quad (\alpha_3-\delta_3)\begin{pmatrix} w \\ \tilde w\end{pmatrix}=\begin{pmatrix} \tilde z\delta^*+\alpha^* z\\ -z\delta-\alpha\tilde z\end{pmatrix}\]
which we write in `matrix' form using the above deformed commutation
relations.
\end{proof}

Note that the `Cartesian' decomposition $\hat B=t+\imath
x\cdot\sigma$ may also be computed but it involves $\sin$ and $\cos$
factors, whereas the `twistor' coordinates, where we work with $\hat
B^i{}_j$ directly as generators, have simple phase factors as above.
Next, we look at the possible cases where $\C_F[\hCM]$ has a
quotient analogous to $\C[S^4]$ in the classical case.  We saw in
the classical case that $\alpha=\delta=0$ and $t,x$ are Hermitian,
or equivalently that
\begin{equation}\label{starS4} z^*=\tilde z,\quad w^*=-\tilde w.\end{equation}
Now in the quantum case the $*$-operation on the entries of $\hat P$
are given by a multiple of the undeformed $*$-operation (as shown in
the proof of Proposition~6.1). Hence the analogous relations in
$\C_F[S^4]$, if it exists as a $*$-algebra quotient, will have the
same form as (\ref{starS4}) but with some twisting factors.

\begin{prop}
The twisting quantisation $\C_F[\hCM]$ is compatible with the
$*$-algebra quotient $\C[\hCM]\to \C[S^4]$ {\em if and only if}
$\theta_A=\theta_B=-{1\over 2}\theta$. In this case
\[  z^*=e^{{\imath\over 2}\theta} \tilde z,\quad w^*=-e^{-{\imath\over 2}\theta}\tilde w\]
\end{prop}
\begin{proof} If the twisting quantisation is compatible with the $*$-algebra quotient, we have
\[ \hat B^1{}_1=e^{-\imath\theta_{13}}B^1{}_1,\quad \hat B^1{}_2=e^{-\imath\theta_{14}}B^1{}_2,
\quad \hat B^2{}_1=e^{-\imath\theta_{23}}B^2{}_1,\quad \hat
B^2{}_2=e^{-\imath\theta_{24}}B^2{}_2\]
from which we deduce the required $*$-operations for the quotient.
For example, $\hat
B^1{}_1{}^*=e^{\imath\theta_{13}+\imath\theta_{24}}\hat B^2{}_2$ and
use the above lemma to identify the factor here as
$e^{\imath\theta/2}$. Therefore we obtain the formulae as stated for
the $*$-structure necessarily in the quotient.  Next, working out
$\C_F[\hCM]$ using Proposition~6.3 we have on the one hand
\[ z \tilde w^*=e^{-2\imath (\theta+\theta_B)}\tilde w{}^* z\]
and on the other hand
\[ z w=e^{2\imath \theta_A}w z.\]
For these to coincide as needed by any relation of the form
(\ref{starS4}) (independently of any deformation factors there) we
need $-(\theta+\theta_B)=\theta_A$. Similarly for compatibility of
the $z^*w$ relation with the $\tilde z w$ relation, we need
$\theta_A=\theta_B$. This determines $\theta_A,\theta_B$ as stated
for the required quotient to be a $*$-algebra quotient. These are
also sufficient as far as the commutation relations are concerned.
The precise form of $*$-structure stated allows one to verify the
other relations in the quotient as well.   \end{proof}

We see that while $\C_F[\hCM]$ has a 3-parameter deformation, there
is only a 1-parameter deformation that pulls back to $\C_F[S^4]$.
The latter has only $a,z,w,z^*,w^*$  as generators with relations
\[ [a,z]=[a,w]=0,\quad zw=e^{-\imath\theta}wz,\quad zw^*=e^{\imath\theta}w^*z, \quad
z^*z+w^*w=a(1-a), \]
which after a minor change of variables is exactly the $S^4_\theta$
in \cite{cl:id}. The `pull-back' of the projector $\hat P$ to
$\C_F[S^4]$ is
\[ \hat e=\begin{pmatrix} a & \hat B\\ \hat B^\dagger &1-a\end{pmatrix},\quad a^*=a,\quad \hat B=\begin{pmatrix}z & -e^{\imath\theta\over 2}w^*\\ w  & e^{-{\imath\theta\over 2}}z^*\end{pmatrix}\]
which up to the change of notations is the `defining projector' in
the Connes-Landi approach to $S^4$. Whereas it is obtained in
\cite{cl:id} from considerations of cyclic cohomology, we obtain it
by a straightforward twisting-quantisation. In view of
Proposition~3.2 we define the noncommutative 1-instanton to be the
Grassmann connection for the projector $\hat e$ on
$\CE=\C_F[S^4]^4\hat e$. This should not come as any surprise since
the whole point in \cite{cl:id} was to define the noncommutative
$S^4$ by a projector generating the K-theory as the 1-instanton
bundle does classically. However, we now obtain $\hat e$ not by this
requirement but by twisting-quantisation and as a `pull-back' of the
tautological bundle on $\C_F[\hCM]$.

Finally, our approach also canonically constructs
$\Omega(\C_F[\hCM])$ and one may check that this quotients in the
one-parameter case to $\Omega(\C_F[S^4])$, coinciding with the
calculus used in \cite{cl:id}. As explained above, the classical
(anti)commutation relations are modified by the same phase factors
as in the commutation relations above. One may then obtain explicit
formulae for the instanton connection and for the Grassmann
connection on $\C_F[\hCM]^4\hat P$.

\subsection{Quantum twistor space $\C_F[\CP^3]$}

In our $*$-algebra approach the classical algebra $\C[\CP^3]$ has a
matrix of generators $Q^\mu{}_\nu$ with exactly the same form as for
$\C[\hCM]$, with the only difference being now $\Tr ~Q=1$, which
significantly affects the content of the `projector' relations of
the $*$-algebra. However, the commutation relations in the quantum
case $\C_F[\CP^3]$ according to Proposition~6.1 have exactly the
same form as $\C_F[\hCM]$ if we use the same cocycle $F$. Hence the
commutation relations between different matrix entries in the
quantum case can be read off from  Proposition~6.3. We describe them
in the special 1-parameter case  found in the previous section where
$\theta_A=\theta_B=-{1\over 2}\theta$.

Using the same notations as in Section~3.2 but now with potentially
noncommutative generators, $\C_F[\CP^3]$ has a matrix of generators
\[ \hat Q=\begin{pmatrix} a & x & y & z\\ x^* & b & w & v \\ y^* & w^* & c & u \\ z^* & v^* & u^* & d\end{pmatrix},\quad a^*=a,\ b^*=b,\ c^*=c,\ d^*=d,\ a+b+c+d=1\]
(we omit hats on the generators). We will use the same shorthand
$X=x^*x$ etc as before. As we know on general grounds above, any
twisting quantisation $\C_F[\CP^3]$ has all entries of the quantum
matrix $\hat Q$ normal (commuting with their adjoints), and the
diagonal elements and quantum versions of $X,Y,Z,U,V,W$ central.
Moreover, all proofs and statements in Section~3.2 were given whilst
being careful {\em not} to assume that $x,y,z,u,v,w$ mutually
commute, only that these elements are normal and $X,Y,Z,U,V,W$ central. Hence the relations
stated there are also exactly the projector relations for this
algebra:

\begin{prop} For the 1-parameter family of cocycles $\theta_A=\theta_B=-{1\over 2}\theta$ the
quantisations $\C_F[\CP^3]$ and $\C^-_F[\CP^3]$ have exactly the
projection relations as in Propositions~3.5 and~3.7 but now with the
commutation relations
\[ xz=e^{\imath\theta}zx,\quad yx=e^{\imath\theta}xy,\quad yz=e^{\imath\theta}zy,\]
and the auxiliary commutation relations
\[ uv=e^{\imath\theta}vu,\quad uw=e^{\imath\theta}wu,\quad vw=e^{\imath\theta}wv,\]
\[ x(u,v,w)=(e^{2\imath\theta}u,e^{\imath\theta}v,e^{-\imath\theta}w)x,\]
\[ y(u,v,w)=(e^{\imath\theta}u,v,e^{-\imath\theta}w)y,\quad z(u,v,w)=(e^{\imath\theta}u,e^{\imath\theta}v,w)z,\]
and similar relations with inverse factor if any generator in a
relation is replaced by its adjoint under $*$. \end{prop}
\begin{proof}  As explained, the commutation relations are the same as for the
entries of $\hat P$ in $\C_F[\hCM]$ with a different notation of the
matrix entries.  We read them off and specialise to the 1-parameter
case of interest.  The `auxiliary' set are deduced from those among
the $x$, $y$, $z$ if $a\ne 0$, since in this case $u$, $v$, $w$ are
given in terms of these and their adjoints.  \end{proof}

If we localise by inverting $a$, then by analogous arguments to the
classical case, the resulting `patch' of $\C^-_F[\CP^3]$ becomes a
quantum punctured $S^6$ with complex generators $x,y,z$, invertible
central self-adjoint generator $a$ and commutation relations as
above, and the relation $x^*x+y^*y+z^*z=a(1-a)$. Also by the same
arguments as in the classical case, if we set $a,b,c,d$ to actual
fixed numbers (which still makes sense since they are central) then
\[ \C^-_F[\CP^3]|_{{b,c,d>0\atop b+c+d=1}}=\C_\theta[S^1\times S^1\times S^1] \]
where the right hand side has relations as above for three circles
but now with commutation relations between the $x,y,z$ circle
generators as stated in the proposition. For each set of values of
$b,c,d$ we have a quantum analogue of $S^1\times S^1\times
S^1\subset \CP^3$ and if we leave them undetermined then in some
sense a quantum version of $\C^*\times\C^*\times\C^*\subset \CP^3$,
i.e. $\C^-_F[\CP^3]$ is in this sense a `quantum toric variety'. We
have seen the same pattern of results already for $\C_F[\CP^2]$ and
$\C[\CP^1]$.

Having established the quantum versions of $\C[S^4]$ and twistor
space $\C^-[\CP^3]$, we now investigate the quantum version of the
fibration $\CP^3 \rightarrow S^4$.  In terms of coordinate algebras
one has an antilinear involution $J: \C_F^-[\CP^3] \rightarrow
\C_F^-[\CP^3]$ analogous to Lemma~3.9. The form of $J$, however, has to be modified by some phase
factors to fit the commutation relations of Proposition~6.5, and is
now given by
\[ J(y) = e^{\imath \theta}v^*, \quad J(y*)=e^{-\imath \theta}v,
\quad J(v) = e^{-\imath \theta}y^*, \quad J(v^*)=e^{\imath
\theta}y,\]
\[ J(w)=z^*, \quad J(w^*)=-z, \quad J(z) = -w^*, \quad J(z^*)=-w,\]
\[ J(x)=-x, \quad J(u)=-u, \quad J(a)=b, \quad J(b)=a.\]
The map $J$ then extends to $\C_F[\hCM]$ and by arguments analogous
to those given in the previous section and in Section~3.2, the fixed
point subalgebra under $J$ is once again precisely $\C_F[S^4]$. We arrive in this way at the
analogous main conclusion, which we verify directly:

\begin{prop}
There is an algebra inclusion
\[\eta : \C_F[S^4] \hookrightarrow \C_F^-[\CP^3]\]
given by
\[ \eta (a) = a+b, \quad \eta(z)=e^{\imath \theta}y+v^*,\quad \eta(w)=w-z^*.\]
\end{prop}

\proof Once again the main relation to investigate is the image of
the sphere relation $zz^*+ww^*=a(1-a)$.  Applying $\eta$ to the left
hand side, we obtain
\[ \eta(zz^*+ww^*) = yy^* + e^{\imath \theta}yv + e^{-\imath \theta}v^*y^* + v^*v + ww^* - wz - z^*w^* + z^*z. \]
We now compute that
\[ ayv = yav = yx^*z = e^{-\imath \theta}x^*yz = e^{-\imath \theta}awz,\]
where the first equality uses centrality of $a$, the second uses the
projector relation $av=x^*z$ and the third uses Proposition~6.5.
Similarly one obtains that $byv = e^{-\imath \theta} bwz$, $cyv =
e^{-\imath \theta}cwz$, $dyv = e^{-\imath \theta}dwz$, so that
adding these four relations now reveals that $yv=e^{-\imath
\theta}wz$ in $\C_F^-[\CP^3]$. Finally using the relations in
Proposition~3.8 (which are still valid as the projector relations in
our noncommutative case) we see that
\[\eta(zz^*+ww^*) = Y + V+ W+ Z = (a+b)(c+d) = (a+b)(1-(a+b)) =
\eta( a(1-a)).\]
To verify the preservation of the algebra structure of $\C_F[S^4]$
under $\eta$ we also have to check the commutation relations, of
which the nontrivial one is $zw=e^{-\imath\theta}wz$. Indeed,
$\eta(zw)=(e^{\imath
\theta}y+v^*)(w-z^*)=e^{-\imath\theta}(w-z^*)(e^{\imath
\theta}y+v^*)=\eta(e^{-\imath\theta}wz)$ using the commutation
relations in Proposition~6.5.
\endproof

Just as in Section~3.2 we may compute the `push-out' of the quantum
instanton bundle along $\eta$, given by viewing the tautological
projector $\hat e$ as an element $\tilde e \in M_4(\C^-_F[\CP^3])$.
Explicitly, we have (following the method of Section~3.2)
\[ \tilde e = \begin{pmatrix}  a+b & \hat M \\ \hat M^\dagger & 1-(a+b) \end{pmatrix} \in
M_4(\C^-_F[\CP^3]), \quad \hat M= \begin{pmatrix} y+v^* & -e^{\imath
\theta \over 2}(w^*-z) \\ w-z^* & e^{-\imath \theta \over 2}(y^*+v)
\end{pmatrix}, \]
and the auxiliary bundle over twistor space is then $\tilde \E
=\C^-_F[\CP^3]^4 \tilde e$. In this way the quantum instanton may be
thought of as coming from a bundle over quantum twistor space, just
as in the classical case.

\section{The Penrose-Ward Transform}
\label{penrose-ward tranform}

The main application of the double fibration (\ref{double
fibration}) is to study relationships between vector bundles over
twistor space and space-time, and between their associated geometric
data. We begin by discussing how differential forms on these spaces
are related, before considering more general vector bundles.
Although we restrict our attention to the special case of
(\ref{double fibration}), our remarks are not specific to this
example and one should keep in mind the picture of a pair of
fibrations of homogeneous spaces
\begin{eqnarray}
\label{general double fibration}
\setlength{\unitlength}{1mm}
\begin{picture}(50,15)(0,0)
 \put(25,14){\makebox(0,0){$G/R$}}
\put(9,0){\makebox(0,0){$G/H$}} \put(41,0){\makebox(0,0){$G/K,$}}
\put(23,11){\vector(-1,-1){9}} \put(27,11){\vector(1,-1){9}}
\put(15,7){\makebox(0,0){$p$}} \put(35,7){\makebox(0,0){$q$}}
\end{picture}
\end{eqnarray}
where $G$ is, say, a complex Lie group with parabolic subgroups
$H,K,R$ such that $R = H \cap K$.  In general, the differential
calculi of the coordinate algebras of these spaces need not be
compatible in any sensible way, however when there is some form of
compatibility (namely a certain transversality condition on the
fibrations at the level of the correspondence space $G/R$) then the
double fibration has some nice properties. Indeed, as we shall see,
such a compatibility between calculi allows one to `transform'
geometric data from $G/H$ to $G/K$ and \textit{vice versa}.

This is the motivation behind the {\em Penrose-Ward transform},
which we shall describe in this section.  The idea is that certain
classes of vector bundles over a given subset of $G/H$ correspond to
vector bundles over an appropriate subset of $G/K$ equipped with a
connection possessing anti-self-dual curvature, the correspondence
given by pull-back along $p$ followed by direct image along $q$.  In
fact, as discussed, we have already seen this transform in action,
albeit in the simplified case where $G/H = \CP^3$ and the given
subset of $G/K = \hCM$ is $S^4$, so that the double fibration
collapses to a single fibration.  In Section~3.2 we gave a
coordinate algebra description of this fibration, with the analogous
quantum version computed in Section~6.2.  In this section we outline
how the transform works in a different situation, namely between the
affine piece $\CM\subset \hCM$ of space-time and the corresponding
patch of twistor space (so that we are utilising the full geometry
of the double fibration): we also outline how it quantises under the
cotwist given in Section~4.

\subsection{Localised coordinate algebras}
We note that classically the driving force behind the construction
that we shall describe is the theory of holomorphic functions and
holomorphic sections of vector bundles.  Of course, if one works
globally then one cannot expect to be able to say very much at all
(in general, one cannot expect to find enough functions).  Hence
here we resort to working \textit{locally}, on some open set $U$ of
space-time $\CM^\#$ having nice topological
properties\footnote{Namely that the intersection of each
$\alpha$-plane $\hat{Z}$ with $U$ is either empty or connected and
simply connected, so that if we write $W := q^{-1}(U)$ and $Y:=p(W)$
then the fibres of the maps $p: W\rightarrow Y$ and $q:  W
\rightarrow U$ are connected, so that operations such as pull-back
and direct image are well-behaved.}. As remarked earlier, of
particular interest are the principal open sets $U_f$, for which one
can explicitly write down the coordinate algebra $\CO_{\hCM}(U_f)$
in our language and expect the construction to go through.  In what
follows, however, we choose to be explicit and restrict our
attention to the case of
\[ \CO_{\hCM}(U_t) = \CC[\TCM](t^{-1})_0,\]
the `affine' piece of space-time, for which we have the inclusion of
algebras $\CC[\CM] \hookrightarrow \CC[\TCM](t^{-1})_0$ given by
\[ x_1 \rightarrow t^{-1}z, \quad x_2 \rightarrow t^{-1}\tz, \quad x_3 \rightarrow
t^{-1}w, \quad x_4 \rightarrow t^{-1}\tw.\]
In fact we also choose to write $\ts := t^{-1} s$, so that the
quadric relation in $\CC[\TCM]$ becomes $\ts = x_1x_2 - x_3x_4$.
Thus in the algebra $\CC[\TCM](t^{-1})_0$ the generator $\ts$ is
redundant and we may as well work with $\CC[\CM]$.

Since the points in $\TCM$ for which $t=0$ correspond to the points
in $\tilde T$ where $Z^3=Z^4=0$, the homogeneous twistor space of
affine space-time $\CC[\CM]$ has coordinate algebra $\C[\tilde T]$
with the extra conditions that $Z^3,Z^4$ are not both zero. The
twistor space $T_t$ of $\CM$ is thus covered by two coordinate
patches, where $Z^3 \neq 0$ and where $Z^4 \neq 0$. These patches
have coordinate algebras
\[ \C[T_{Z^3}]:= \C[\tilde T]((Z^3)^{-1})_0, \qquad \C[T_{Z^4}] := \C[\tilde
T]((Z^4)^{-1})_0\]
respectively.  Note that when both $Z^3$ and $Z^4$ are non-zero the
two algebras are isomorphic (even in the twisted case, since
$Z^3,Z^4$ remain central in the algebra), with `transition
functions' $Z^\mu(Z^3)^{-1} \mapsto Z^\mu(Z^4)^{-1}$.  This
isomorphism simply says that both coordinate patches look locally
like $\C^3$, in agreement that our expectation that twistor space is
a complex 3-manifold.

In passing from $\hCM$ to $\CM$ we delete the `region at infinity'
where $t=0$, and similarly we obtain the corresponding twistor space
$T_t$ by deleting the region where $Z^3=Z^4=0$.  At the homogeneous
level this region has coordinate algebra
\[ \C[\tilde T]/\la Z^3=Z^4=0 \ra \cong \C[Z^1,Z^2],\]
describing a line $\CP^1$ at the projective level.  In the twisted
framework, the generators $Z^3,Z^4$ are central so the quotient
still makes sense and we have
\[ \C_F[\tilde T]/\la Z^3=Z^4=0 \ra \cong \C_F[Z^1,Z^2],\]
describing a noncommutative $\CP^1$.

At the level of the homogeneous correspondence space we also have
$t^{-1}$ adjoined and $Z^3, Z^4$ not both zero, and we write
$\CC[\tilde \F_t]: = \C[\tilde \F](t^{-1})_0$ for the corresponding
coordinate algebra. We note that in terms of the $\CC[\CM]$
generators, the relations (\ref{alpha planes}) now read
\begin{equation} \label{trivialisation of taut bundle} Z^1 = x_2 Z^3 + x_3 Z^4, \quad Z^2 = x_4 Z^3 + x_1
Z^4.\end{equation}
At the projective level, the correspondence space is also covered by
two coordinate patches.  Thus when $Z^3 \ne 0$ and when $Z^4 \ne 0$
we respectively mean
\[ \C[\F_{Z^3}]:= \C[\tilde \F_t]((Z^3)^{-1})_0, \qquad \C[\F_{Z^4}]:= \C[\tilde \F_t]((Z^4)^{-1})_0.\]
Again when $Z^3$ and $Z^4$ are both non-zero these algebras are seen
to be isomorphic under appropriate transition functions.  Thus at
the homogeneous level, the coordinate algebra $\C[\tilde \F_t]$
describes a local trivialisation of the correspondence space in the
form $\tilde \F_t \cong \C^2 \times \CM$.  The two coordinate
patches $\F_{Z^3}$ and $\F_{Z^4}$ together give a trivialisation of
the projective correspondence space in the form $\F_t = \CP^1 \times
\CM$.

Regarding the differential calculus of $\CC[\TCM](t^{-1})_0$, it is
easy to see that since $\D(t^{-1}t) = 0$, from the Leibniz rule we
have $\D(t^{-1}) = - t^{-2}\D t$.  Since $\D t$ is central even in
the twisted calculus, adjoining this extra generator $t^{-1}$ causes
no problems, and it remains to check that the calculus is
well-defined in the degree zero subalgebra $\CC[\TCM](t^{-1})_0$.
Indeed, we see that for example
\[ \D x_1 = \D (t^{-1} z) = t^{-1}\D z - (t^{-1}z)(t^{-1}\D t),\]
which is again of overall degree zero.  Similar statements hold
regarding the differential calculi of $\CC[T_{Z^3}]$ and
$\CC[T_{Z^4}]$ as well as those of $\CC[\F_{Z^3}]$ and
$\CC[\F_{Z^4}]$.

\subsection{Differential aspects of the double fibration}
We now consider the pull-back and direct image of one-forms on our
algebras. Indeed, we shall examine how the differential calculi
occurring in the fibration (\ref{double fibration}) are related and
derive the promised `transversality condition' required in order to
transfer data from one side of the fibration to the other.  For now,
we consider only the classical (i.e. untwisted) situation.

Initially we work at the homogeneous level, with $\C[\tilde T]$ and
$\C[\tilde \F_t]$, later passing to local coordinates by adjoining
an inverse for either $Z^3$ or $Z^4$, as described above. We define
\begin{equation} \label{relative forms}
\Omega^1_p : = \Omega^1 \CC[\tilde \F_t] / p^* \Omega^1 \CC[\tilde
T] \end{equation}
to be the set of \textit{relative} one-forms (the one-forms which
are dual to those vectors which are tangent to the fibres of $p$).
Note that $\Omega^1_p$ is just the sub-bimodule of $\Omega^1
\C[\tilde \F_t]$ spanned by $\D \tilde s$, $\D x_1$, $\D x_2$, $\D
x_3$, $\D x_4$, so that
\begin{equation} \label{relative decomposition} \Omega^1 \CC[\tilde \F_t] = p^* \Omega^1 \CC[\tilde T] \oplus
\Omega^1_p.\end{equation}
There is of course an associated projection
\[ \pi_p: \Omega^1 \CC[\tilde \F_t] \rightarrow \Omega^1_p\]
and hence an associated relative exterior derivative $\D_p:
\CC[\tilde \F_t] \rightarrow \Omega^1_p$ given by composition of
$\D$ with this projection,
\begin{equation} \label{relative connection} \D_p = \pi_p \circ \D: \CC[\tilde \F_t] \rightarrow \Omega^1_p.\end{equation}
Similarly, we define the relative two-forms by
\[ \Omega^2_p := \Omega^2 \CC[\tilde \F_t] / (p^* \Omega^1 \C[\tilde T] \wedge \Omega^1 \C[\tilde \F_t]) \]
so that $\D_p$ extends to a map
\[ \D_p:\Omega^1_p \rightarrow \Omega^2_p\]
by composing $\D:\Omega^1 \C[\tilde \F_t] \rightarrow \Omega^2 \C[\tilde \F_t]$ with the projection $\Omega^1 \C[\tilde \F_t]  \rightarrow \Omega^2_p$. We see that $\D_p$ obeys the relative Leibniz rule,
\begin{equation} \label{partial leibniz} \D_p (fg) = (\D_p f) g + f(\D_p g), \quad f,g \in \CC[\tilde
\F_t].\end{equation}
It is clear by construction that the kernel of $\D_p$ consists
precisely of the functions in $\CC[\tilde \F_t]$ which are constant
on the fibres of $p$, whence we recover $\CC[\tilde T]$ by means of
the functions in $\CC[\tilde \F_t]$ which are covariantly constant
with respect to $\D_p$ (since functions in $\C[\tilde T]$ may be
identified with those functions in $\C[\tilde \F_t]$ which are
constant on the fibres of $p$).  Moreover, the derivative $\D_p$ is
\textit{relatively flat}, i.e. its curvature $\D_p^2$ is zero.

The next stage is to consider the direct image of the one-forms
$\Omega^1_p$ along $q$.  In the usual theory the direct image
\cite{ww:tgft} of a vector bundle $\pi: E' \rightarrow \F_t$ over
$\F_t$ along the fibration $q: \F_t \rightarrow \CM$ is by
definition the bundle $E:=q_*E' \rightarrow \CM$ whose fibre over $x
\in \CM$ is $H^0(q^{-1}(x), E')$, the space of global sections of
the restriction of $E'$ to $q^{-1}(x)$.  Of course, this definition
does not in general result in a well-defined vector bundle.  For
this, we assume that each $q^{-1}(x)$ is compact and connected so
that $H^0(q^{-1}(x), E')$ is finite dimensional, and that
$H^0(q^{-1}(x), E')$ is of constant dimension as $x \in \CM$ varies.
In our situation these criteria are clearly satisfied.

We recall that in the given local trivialisation, the correspondence
space looks like $\F_t = \CP^1 \times \CM$ with coordinate functions
$\zeta$, $x_1$, $x_2$, $x_3$, $x_4$, where $\zeta = Z^3(Z^4)^{-1}$
or $Z^4(Z^3)^{-1}$ (depending on which coordinate patch we are in).
Writing $\E'$ for the space of sections of the bundle $E'
\rightarrow \F_t$, the space of sections $\E$ of the direct image
bundle $E$ is in general obtained by computing $\E'$ as a
$\C[\CM]$-module.  So although any section $\xi \in \E'$ of $E'$ is
in general a function of the twistor coordinate $\zeta$ as well as
the space-time coordinates $x_j$, its restriction to any fibre
$q^{-1}(x)$ is by Liouville's theorem independent of $\zeta$ (i.e.
$\zeta$ is constant on $q^{-1}(x)$). Of course, when we restrict the
section $\xi$ to each $q^{-1}(x)$, its dependence on $\zeta$ varies
as $x \in \CM$ varies (in that although $\zeta$ is constant on each
$q^{-1}(x)$, it possibly takes different values as $x$ varies), and
this means that in taking the direct image we may write $\zeta$ as a
function of the space-time coordinates $x_j$, $j=1,\ldots 4$.  One
could also argue by by noting that $\CM$ is topologically trivial so
has no cohomology, whence one has for example
\[H^0(\CP^1 \times \CM,\C) = H^0(\CP^1,\C)\otimes H^0(\CM,\C) \cong
H^0(\CM,\C).\]

At the homogeneous level, although the fibres of the map $q: \tilde
\F_t \rightarrow \CM$ are not compact, the effect of the above
argument may be achieved by regarding the twistor coordinates $Z^3,
Z^4$ as functions of the space-time coordinates under the direct
image.  The upshot of this argument is that to compute the direct
image of a $\C[\tilde \F_t]$-module $\E'$, we `remove' the
dependence of $\E'$ on the twistor coordinates $Z^i$, hence
obtaining a $\C[\CM]$-module $\E$.  In what follows, this will be
our naive definition of direct image, by analogy with the classical
case.  Since this argument really belongs in the language of
cohomology, we satisfy ourselves with this definition for now,
deferring a more precise treatment to a sequel.

\begin{prop} \label{pull back direct image}
There is an isomorphism $q_*q^* \Omega^1 \CC[\CM] \cong \Omega^1
\CC[\CM]$.
\end{prop}
\proof  The generators $\D x_j$, $j=1,\ldots 4$ of $\Omega^1
\C[\CM]$ pull back to their counterparts $\D x_i, i=1,\ldots 4$ in
$\CC[\tilde \F_t]$ and these span $q^* \Omega^1 \CC[\CM]$ as a
$\C[\tilde \F_t]$-bimodule. Taking the direct image involves
computing $q^* \Omega^1 \CC[\CM]$ as a $\CC[\CM]$-bimodule, and as
such it is spanned by elements of the form $\D \tilde s$, $\D x_i$,
$Z^i \D \ts$, $Z^i \D x_j$, for $i,j = 1,\ldots 4$. As already
observed, the generator $\ts$ is essentially redundant, hence so are
the generators involving $\D \ts$. Moreover, the relations
(\ref{trivialisation of taut bundle}) allow us to write $Z^1$ and
$Z^2$ in terms of $Z^3, Z^4$, whence we are left with elements of
the form $\D x_j$, $Z^3 \D x_j$ and $Z^4 \D x_j$ for $j = 1,\ldots
4$.

As explained above, the direct image is given by writing the
$Z$-coordinates as functions of the space-time coordinates, and it
is clear that the resulting differential calculus $q_*q^* \Omega^1
\CC[\CM]$ of $\C[\CM]$ is isomorphic to the one we first thought of.
We remark that we have used the same symbol $\D$ to denote the
exterior derivative in the different calculi $\Omega^1 \C[\CM]$,
$q^* \Omega^1 \CC[\CM]$ and $q_*q^* \Omega^1 \CC[\CM]$. Although
they are in principle different, our notation causes no confusion
here.
\endproof

\begin{prop}
\label{direct images of forms} For $s=1,2$ the direct images $q_*
\Omega^s_p$ are given by:
\begin{itemize}
\item $q_* \Omega^1_p \cong \Omega^1 \CC[\CM]$;
\item $q_* \Omega^2_p \cong \Omega^2_+\CC[\CM]$,
\end{itemize}
where $\Omega^2_+\CC[\CM]$ denotes the space of two-forms in
$\Omega^2 \CC[\CM]$ which are self-dual with respect to the Hodge
$*$-operator defined by the metric $\eta = 2(\D x_1 \D x_2 - \D x_3
\D x_4)$.
\end{prop}
\proof We first consider the $\CC[\tilde \F_t]$-bimodule
$\Omega^1_p$, and write $q_* \Omega^1_p$ for the same vector space
considered as a $\C[\CM]$-bimodule.  Quotienting $\Omega^1 \C[\tilde
\F_t]$ by the one-forms pulled back from $\CC[\tilde T]$ means that
$\Omega^1_p$ is spanned as a $\CC[\tilde \F_t]$-bimodule by $\D \ts$
and $\D x_i$, $i=1,\ldots 4$.

The direct image is computed just as before, so $q_* \Omega^1_p$ is
as a vector space the same calculus $\Omega^1_p$ but now considered
as a $\C[\CM]$-bimodule.  In what follows we shall write $\D$ for
the exterior derivatives in the calculi $\Omega^1 \C[\CM]$ and
$\Omega^1 \C[\tilde \F_t]$ as they are the usual operators (the ones
we which we wrote down and quantised in section \ref{theta quantum
calculi}). However, in calculating the direct image $q_* \Omega^1_p$
of the calculus $\Omega^1_p$, we must introduce different notation
for the image of the operator $\D_p$ under $q_*$. To this end, we
write $q_* \D_p := \uD$, so that as a $\C[\CM]$-bimodule the
calculus $q_*\Omega^1_p$ is spanned by elements of the form $\uD
\ts$ and $\uD x_j$, $i,j=1,\ldots 4$. As already observed, the
generator $\ts$ is essentially redundant, hence so is the generator
$\uD \ts$.

The identity (\ref{partial leibniz}) becomes a Leibniz rule for
$\uD$ upon taking the direct image $q_*$, whence $(q_* \Omega^1_p,
\uD)$ is a first order differential calculus of $\CC[\CM]$.  We must
investigate its relationship with the calculus $\Omega^1 \C[\CM]$.

Differentiating the relations (\ref{trivialisation of taut bundle})
and quotienting by generators $\D Z^i$ yields the relations
\begin{equation} \label{relative calculus relations1} Z^3 \D_p x_2 + Z^4 \D_p x_3 = 0, \quad Z^3 \D_p x_4 + Z^4 \D_p x_1
=0 \end{equation}
in the relative calculus $\Omega^1_p$.  Thus as a $\C[\tilde
\F_t]$-bimodule, $\Omega^1_p$ has rank two, since the basis elements
$\D_p x_j$ are not independent.  However, in the direct image $q_*
\Omega^1_p$, which is just $\Omega^1_p$ considered as a
$\C[\CM]$-bimodule, these basis elements (now written $\uD x_j$) are
independent.  Thus it is clear that the calculus $(q_*
\Omega^1_p,\uD)$ is isomorphic to $(\Omega^1 \C[\CM], \D)$ in the
sense that as bimodules they are isomorphic, and that this
isomorphism is an intertwiner for the derivatives $\uD$ and $\D$.

Using the relations (\ref{relative calculus relations1}) it is easy
to see that in the direct image bimodule $q_* \Omega^2_p :=
\Lambda^2 (q_* \Omega^1_p)$ we have
\[ \uD x_1 \wedge \uD x_4 = \uD x_2 \wedge \uD x_3 = \uD x_1 \wedge \uD
x_2 + \uD x_3 \wedge \uD x_4 = 0,\]
which we recognise (since we are in double null coordinates) as the
anti-self-dual two-forms, whence it is the self-dual two-forms which
survive under the direct image. It is evident that these arguments
are valid upon passing to either coordinate patch at the projective
level, i.e. upon adjoining either $(Z^3)^{-1}$ or $(Z^4)^{-1}$ and
taking the degree zero part of the resulting calculus.
\endproof

It is clear that the composition of maps
\[ q^* \Omega^1 \CC[\CM] \rightarrow \Omega^1 \CC[\tilde \F_t] \rightarrow
\Omega^1_p\]
determines the isomorphism
\[ \Omega^1 \CC[\CM] \rightarrow q_* \Omega^1_p\]
by taking direct images\footnote{We note that for the general double
fibration (\ref{general double fibration}), a similar analysis will
go through provided one has the transversality condition that $q_*
\Omega^1_p \cong \Omega^1 \C[G/H]$.} (this is also true in each of
the coordinate patches). The sequence
\begin{equation}
\label{correspondence space sequence} \CC[\tilde \F_t] \rightarrow
\Omega^1_p \rightarrow \Omega^2_p,
\end{equation}
where the two maps are just $\D_p$, becomes the sequence
\[\CC[\CM] \rightarrow \Omega^1 \CC[\CM] \rightarrow
\Omega^2_+\CC[\CM]\]
upon taking direct images (the maps become $\uD$).  The condition
that $\D_p^2 = 0$ is then equivalent to the statement that the
curvature $\uD^2$ is annihilated by the map
\[ \Omega^2 \CC[\CM] \rightarrow q_* \Omega^2_p = \Omega^2_+,\]
i.e. that the curvature $\uD^2$ is anti-self-dual. So although the
connection $\D_p$ is flat, its image under $q_*$ is not. Although
the derivatives $\D$ and $\uD$ agree as maps $\C[\CM] \rightarrow
\Omega^1 \C[\CM]$, they do not agree beyond the one-forms.  The
image of $\D: \Omega^1 \C[\CM] \rightarrow \Omega^2 \CC[\CM]$
consists of all holomorphic two-forms, whereas $\uD$ maps $\Omega^1
\C[\CM] \cong q_* \Omega^1_p$ onto the anti-self-dual two-forms. At
the level of the correspondence space, the reason for this is that
in the sequence (\ref{correspondence space sequence}), the calculus
$\Omega^1_p$ comes equipped with relations (\ref{relative calculus
relations1}), whereas the pull-back $q^*\Omega^1 \C[\CM]$ has no
such relations.

We are now in a position to investigate how this construction
behaves under the twisting discussed in Sections \ref{section
quantum conformal group} and \ref{theta quantum calculi}.  Of
course, we need only check that the various steps of the procedure
remain valid under the quantisation functor. To this end, our first
observation is that although the relations in the first order
differential calculus of $\C[\CM]$ are deformed, equations
(\ref{affine minkowski calculus}) show that the two-forms are
undeformed, as is the metric $\upsilon$, hence the Hodge
$*$-operator is also undeformed in this case. It follows that the
notion of anti-self-duality of two-forms is the same as in the
classical case.

Furthermore, the definition of relative one-forms (\ref{relative
forms}) still makes sense since the decomposition (\ref{relative
decomposition}) is clearly unaffected by the twisting.  Since we are
working with the affine piece of (now noncommutative) space-time,
the relevant relations in the correspondence space algebra
$\C_F[\tilde \F_t]$ are given by (\ref{trivialisation of taut
bundle}), which are unchanged under twisting since $t$, $Z^3$, $Z^4$
remain central in the algebra.

Finally we observe that the proofs of Propositions \ref{direct
images of forms} and \ref{pull back direct image} go through
unchanged. The key steps use the fact that when the generators
$Z^3$, $Z^4$ and $t$ are invertible, one may adjoin their inverses
to the coordinate algebras and differential calculi and take the
degree zero parts. Since these generators remain central under
twisting, this argument remains valid and we have the following
twisted consequence of Proposition \ref{direct images of forms}.

\begin{prop}
Let $\D_p : \C_F [\tilde \F_t] \rightarrow \Omega^1_p$ be the
differential operator defined by the composition of maps $\D_p = \pi
\circ \D$,
\[ \D_p: \C_F [\tilde \F_t] \rightarrow \Omega^1 \CC_F[\tilde F_t]
\rightarrow \Omega^1_p.\]
Then the direct image $\uD: \C_F[\CM] \rightarrow q_*\Omega^1_p
\cong \Omega^1 \C_F [\CM]$ is a differential operator whose
curvature $\uD^2$ takes values in the anti-self-dual two-forms
$\Omega^2_- \C_F[\CM]$.
\end{prop}

We remark that this is no coincidence:  the transform between
one-forms on $T_t$ and the operator $\uD$ on the corresponding
affine patch $\CM$ of space-time goes through to this noncommutative
picture precisely because of the choice of cocycle made in Section
\ref{section quantum conformal group}.  Indeed, we reiterate that
any construction which is covariant under a chosen symmetry group
will also be covariant after applying the quantisation functor.  In
this case, the symmetry group is the subgroup generated by conformal
translations (see Section \ref{section quantum conformal group}),
which clearly acts covariantly on affine space-time $\CM$.  By the
very nature of the twistor double fibration, this translation group
also acts covariantly on the corresponding subsets $\tilde T$ and
$\tilde \F_t$ of (homogeneous) twistor space and the correspondence
space respectively, and it is therefore no surprise that the
transform outlined above works in the quantum case as well.

As discussed, it is true classically that one can expect such a
transform between subsets of $\hCM$ and the corresponding subsets of
twistor space $T$ provided the required topological properties (such
as connectedness and simple-connectedness of the fibres) are met. We
now see, however, that the same is not necessarily true in the
quantum case.  For a given open subset $U$ of $\hCM$, we expect the
transform between $U$ and its twistor counterpart
$\hat{U}=p(q^{-1}(U)) \subset T$ to carry over to the quantum case
provided the twisting group of symmetries is chosen in a way so as
to preserve $U$.

\subsection{Outline of the Penrose-Ward transform for vector bundles}
The previous section described how the one-forms on twistor space
give rise to a differential operator on forms over space-time having
anti-self-dual curvature.  The main feature of this relationship is
that bundle data on twistor space correspond to differential data on
space-time.  The idea of the full Penrose-Ward transform is to
generalise this construction from differential forms to sections of
more general vector bundles.

We begin with a finite rank $\CC[\tilde T]$-module\footnote{In this
section we shall for convenience work with left modules over the
algebras in question.} $\tilde \E$ describing the holomorphic
sections of a (trivial) holomorphic vector bundle $\tilde E$ over
homogeneous twistor space $\tilde T$.  The pull-back $p^*\E$ is the
$\C[\tilde \F_t]$-module
\begin{equation}
\label{pull back decomposition} \tilde \E':= p^*\tilde \E =
\CC[\tilde \F_t] \otimes_{\CC[\tilde T]} \tilde \E.
\end{equation}
The key observation is then that there is a relative connection
$\n_p$ on $p^*\tilde \E$ defined by
\[ \n_p  = \D_p \otimes 1\]
with respect to the decomposition (\ref{pull back decomposition}).
Again there is a relative Leibniz rule
\begin{equation} \label{partial Leibniz connection} \n_p(f \xi) =f(\n_p \xi) + (\D_p f) \otimes \xi
\end{equation}
for $f \in \C[\tilde \F_t]$ and $\xi \in p^*\tilde \E$.  Moreover,
$\n_p$ extends to $\C[\tilde \F_t]$-valued $k$-forms by defining
\[ \Omega^k_p \tilde \E' := \Omega^k_p
\otimes_{\CC[\tilde \F_t]} \tilde \E' \]
for $k \geq 0$ and extending $\n_p = \D_p \otimes 1$ with respect to
this decomposition.  It is clear that the curvature satisfies
$\n_p^2 = 0$ and we say that $\n_p$ is \textit{relatively flat}.

Conversely, if the fibres of $p$ are connected and simply connected
(as they clearly are in our situation; we assume this in the general
case of (\ref{general double fibration})), then if $\tilde \E'$ is a
finite rank $\C[\tilde \F_t]$-module admitting a flat relative
connection $\n'$ (that is, a complex-linear map $\n':\tilde \E'
\rightarrow \Omega^1_p \tens_{\C[\tilde \F_t]} \tilde \E'$
satisfying (\ref{partial Leibniz connection}) and $(\n')^2=0$), we
may recover a finite rank $\C[\tilde T]$-module $\tilde \E$ by means
of the covariantly constant sections, namely
\[ \tilde \E : = \{ \xi \in \tilde \E' ~|~ \n' \xi = 0 \}. \]
This argument gives rise to the following result.

\begin{prop}
\label{holomorphic modules}
There is a one-to-one correspondence between finite rank $\CC[\tilde
T]$-modules $\tilde \E$ and finite rank $\CC[\tilde \F_t]$-modules
$\tilde \E'$ admitting a flat relative connection $\n_p$,
\[ \n_p(f\xi) =f(\n_p \xi) + (\D_p f) \otimes \xi, \qquad
\n_p^2=0\]
for all $f \in \CC[\tilde \F_t]$ and $\xi \in \tilde \E'$.
\end{prop}

We remark that here we do not see the non-trivial structure of the
bundles involved (all of our modules describing vector bundles are
free) as in our local picture all bundles are trivial.  However,
given these local formulae it will be possible at a later stage to
patch together what happens at the global level.

The Penrose-Ward transform arises by considering what happens to a
relative connection $\n_p$ under direct image along $q$.  We here
impose the additional assumption that $\tilde \E'$ is also the
pull-back of a bundle on space-time, so $\tilde \E' = p^* \tilde \E
= q^* \E$ for some finite rank $\C[\CM]$-module $\E$.  This is
equivalent to assuming that the bundle $\tilde E'$ is trivial upon
restriction to each of the fibres of the map $q: \tilde \F_t
\rightarrow \CM$.  The direct image is computed exactly as described
in the previous section.

Thus the direct image of $\tilde \E'$ is $q_* \tilde \E' = q_*q^*\E
\cong \E$ and, just as in Proposition \ref{direct images of forms},
equation (\ref{partial Leibniz connection}) for $\n_p$ becomes a
Leibniz rule for $\un : = q_* \n$, whence $\n_p$ maps onto a genuine
connection on $\E$. The sequence
\[\tilde \E' \rightarrow \Omega^1_p \tens_{\C[\tilde \F_t]} \tilde \E' \rightarrow
\Omega^2_p \tens_{\C[\tilde \F_t]} \tilde \E',\]
where the two maps are $\n_p$, becomes the sequence
\[\E \rightarrow \Omega^1 \CC[\CM] \tens_{\C[\CM]} \E \rightarrow
\Omega^2_+ \tens_{\C[\CM]} \E, \]
where the maps here are $\un$.  Moreover, just as in the previous
section it follows that under direct image the condition that
$\n_p^2 = 0$ is equivalent to the condition that the curvature
$\un^2$ is annihilated by the mapping
\[ \Omega^2 \C[\CM] \otimes_{\C[\CM]} \E \rightarrow q_* \Omega^2_p \otimes_{\C[\CM]} \E = \Omega^2_+ \otimes_{\C[\CM]} \E,\]
so that $\un$ has anti-self-dual curvature.

\subsection{Tautological bundle on $\CM$ and its Ward transform} We remark that in the previous section we began with a bundle over
homogeneous twistor space $\tilde T$, whereas it is usual to work
with bundles over the projective version $T$.  As such, we
implicitly assume that in doing so we obtain from $\tilde \E$
corresponding $\C[T_{Z^3}]$- and $\C[T_{Z^4}]$-modules which are
compatible in the patch where $Z^3$ and $Z^4$ are both non-zero, as
was the case for the coordinate algebras $\C[T_{Z^3}]$ and
$\C[T_{Z^4}]$, $\C[\F_{Z^3}]$ and $\C[\F_{Z^4}]$ \textit{via} the
transition functions $Z^\mu (Z^3)^{-1} \mapsto Z^\mu (Z^4)^{-1}$. We
similarly assume this for the corresponding calculi on these
coordinate patches.

In order to neatly capture these issues of coordinate patching, the
general construction really belongs in the language of cohomology:
as explained, the details will be addressed elsewhere.  For the time
begin we give an illustration of the transform in the coordinate
algebra framework, as well as an indication of what happens under
twisting, through the tautological example introduced in Section
\ref{section taut bundle on cmhash}.

We recall the identification of conformal space-time and the
correspondence space as flag varieties $\hCM = \F_2(\C^4)$ and $\F =
\F_{1,2}(\C^4)$ respectively, and the resulting fibration
\[q: \F_{1,2}(\C^4) \rightarrow \F_2(\C^4),\]
where the fibre over a point $x \in \F_2(\C^4)$ is the set of all
one-dimensional subspaces of $\C^4$ contained in the two-plane $x
\subset \C^4$, so is topologically a projective line $\CP^1$. We
also have a fibration at the homogeneous level,
\[ \tilde \F \rightarrow \F_2(\C^4),\]
where this time the fibre over $x \in \F_2(\C^4)$ is the set of all
vectors which lie in the two-plane $x$.  As explained, there are no
non-constant global holomorphic sections of this bundle, so as
before we pass to the affine piece of space-time $\CM$ in order to
avoid this trivial case.

We identify $\C^4$ with its dual and take the basis
$(Z^1,Z^2,Z^3,Z^4)$.  In the patch $\tilde \F_t$ the relations in
$\C[\tilde \F_t]$ are
\begin{equation} \label{taut trivialisation again} Z^1 = x_2 Z^3 + x_3 Z^4, \quad Z^2 = x_4 Z^3 + x_1
Z^4,\end{equation}
which may be seen as giving a trivialisation of the tautological
bundle over $\hCM$ in this coordinate patch.  In this trivialisation
the space $\E$ of sections of the bundle is just the free module
over $\C[\CM]$ of rank two, spanned by $Z^3$ and $Z^4$, i.e. $\E
\cong \C[\CM] \otimes \C^2$. We equip this module with the
anti-self-dual connection $\uD \otimes 1$ constructed in
Section~7.2.

It is now easy to see that the pull-back $\tilde \E'$ of $\E$ along
$q$ is just the free $\C[\tilde \F_t]$-module of rank two.  By
construction, the connection $\uD$ on $\E$ pulls back to the
relative connection $\D_p \otimes 1$ on $\tilde \E' = \C[\tilde
\F_t] \otimes \C^2$.  As discussed earlier the corresponding
$\C[\tilde T]$-module $\tilde \E$ is obtained as the kernel of the
partial connection $\n_p$, which is precisely the free rank two
$\C[\tilde T]$-module $\tilde \E = \C[\tilde T] \otimes \C^2$.  It
is also clear that these modules satisfy the condition that $\tilde
\E' = p^* \tilde \E = q^* \E$ (the corresponding vector bundles are
trivial in this case, whence they are automatically trivial when
restricted to each fibre of $p$ and of $q$).

Lastly, it is obvious that the bundle over $\tilde T$ described by
the module $\tilde \E$ descends to a (trivial) bundle over the
twistor space $T_t$ of $\CM$.  On the coordinate patches where $Z^3$
and $Z^4$ are non-zero one obtains respectively a free rank two
$\C[T_{Z^3}]$-module $\tilde \E_{Z^3}$ by inverting $Z^3$ and a
$\C[T_{Z^4}]$-module $\tilde \E_{Z^4}$ by inverting $Z^4$, and since
the algebras $\C[T_{Z^3}]$ and $\C[T_{Z^4}]$ agree when $Z^3$, $Z^4$
are both non-zero, the same is true of the sections in the localised
modules $\tilde \E_{Z^3}$ and $\tilde \E_{Z^4}$. The construction in
this example remains valid under the twisting described in Section
\ref{quantum conformal group} due to our earlier observation that
the relations (\ref{taut trivialisation again}) are unchanged.

Explicitly then, the tautological bundle of $\hCM$ in the affine
patch $\CM$ is simply the free rank two $\C[\CM]$-module $\E =
\C[\CM]^2 = \C[\CM]\otimes \C^2$, which we equip with the
anti-self-dual connection $\uD \otimes 1$ constructed in Section~7.2
(note that one reserves the term {\em instanton} specifically for
anti-self-dual connections over $S^4$).  Its Penrose-Ward transform
is the trivial rank two holomorphic bundle over the twistor space
$T_t$ of $\CM$.  Since these spaces are topologically trivial, we
see that the Penrose-Ward transform is here very different in
flavour to the transform given in Section~3.2.  It is however clear
that this example of the transform quantises in exactly the same way
as the rank one case of Section~7.2.

\section{The ADHM Construction}

\subsection{The classical ADHM construction}
We begin this section with a brief summary of the ADHM construction
for connections with anti-self-dual curvature on vector bundles over
Minkowski space $\CM$ \cite{adhm:ci,mw:isdtt}, with a view to
dualising and then twisting the construction.

A monad over $\tilde T = \mathbb{C}^4$ is a sequence of linear maps
\begin{center}
\setlength{\unitlength}{1mm}
\begin{picture}(50,15)(0,0)
\put(10,5){\makebox(0,0){A}}
\put(25,5){\makebox(0,0){B}}
\put(40,5){\makebox(0,0){C,}}
\put(17.5,8){\makebox(0,0){$\rho_{_{Z}}$}}
\put(32.5,8){\makebox(0,0){$\tau_{_{Z}}$}}
\put(13,5){\vector(1,0){9}}
\put(28,5){\vector(1,0){9}}
\end{picture}
\end{center}
between complex vector spaces $A,B,C$ of dimensions $k,2k+n, k$
respectively, such that for all $Z \in \mathbb{C}^4$, $\tau_{_{Z}}
\rho_{_{Z}} : A \longrightarrow C$ is zero and for all $Z \in \C^4$,
$\rho_{_{Z}}$ is injective and $\tau_{_{Z}}$ is surjective.
Moreover, we insist that $\rho_{_{Z}}$, $\tau_{_{Z}}$ each depend
linearly on $Z \in \tilde T$.  The spaces $A,B,C$ should be thought
of as typical fibres of trivial vector bundles over $\CP^3$ of ranks
$k, 2k+n, k$, respectively.

A monad determines a rank-$n$ holomorphic vector bundle on $T =
\mathbb{C}\mathbb{P}^3$ whose fibre at $[Z]$ is $\textup{Ker} \;
\tau_{_{Z}} / \textup{Im} \; \rho_{_{Z}}$ (where $[Z]$ denotes the
projective equivalence class of $Z \in \mathbb{C}^4$). Moreover, any
holomorphic vector bundle on $\mathbb{C}\mathbb{P}^3$ trivial on
each projective line comes from such a monad, unique up to the
action of $\textup{GL}(A) \times \textup{GL}(B) \times
\textup{GL}(C)$. For a proof of this we refer to \cite{oss:vb},
although we note that the condition $\tau_{_{Z}} \rho_{_{Z}} = 0$
implies $\textup{Im} \; \rho_{_{Z}} \subset \textup{Ker} \;
\tau_{_{Z}}$ for all $z \in \mathbb{C}^4$ so the cohomology makes
sense, and the fact that $\rho_{_{Z}}, \tau_{_{Z}}$ have maximal
rank at every $Z$ implies that each fibre has dimension $n$.  The
idea behind the ADHM construction is to use the same monad data to
construct a rank $n$ vector bundle over $\hCM$ with second Chern
class $c_2 = k$ (in the physics literature this is usually called
the \textit{topological charge} of the bundle).

For each $W,Z \in \tilde T$ we write $x = W \wedge Z$ for the
corresponding element $x$ of (homogeneous) conformal space-time
$\TCM \subset \Lambda^2 \tilde T$. Then define
$$E_x  = \textup{Ker} \; \tau_{_{Z}} \cap \textup{Ker} \; \tau_{_{W}},
\qquad F_x = (\rho_{_{Z}} A) \cap (\rho_{_{W}} A), \qquad \Delta_x =
\tau_{_{Z}} \rho_{_{W}}.$$

\begin{prop} \cite{mw:isdtt}
The vector spaces $E_x, F_x$ and the map $\Delta_x$ depend on $x$,
rather than on $Z,W$ individually.
\end{prop}
\proof We first consider $\Delta_x$ and suppose that $x = Z \wedge W
= Z \wedge W'$ for some $W'=W+\lambda Z, \lambda \in \C$.  Then we
have
\[ \tau_{_{Z}} \rho_{_{W'}} - \tau_{_{Z}} \rho_{_{W}} = \tau_{_{Z}}
\rho _{_{W'-W}} = \tau_{_{Z}} \rho_{_{\lambda Z}} =\lambda
\tau_{_{Z}} \rho_{_{Z}} = 0.\]
Now writing $Z' = Z + \lambda W'$ we see that
\[ \tau_{_{Z'}} \rho_{_{W'}} - \tau_{_{Z}} \rho_{_{W}} = (\tau_{_{Z}} + \lambda \tau _{_{W'}})\rho_{_{W'}} - \tau_{_{Z}} \rho_{_{W}} = \lambda \tau_{_{Z}}
\rho _{_{W'}} - \tau_{_{Z}} \rho_{_{W}} = 0\]
by the first calculation, proving the claim for $\Delta_x$.

Next we consider $b \in \textup{Ker} \; \tau_{_{Z}} \cap
\textup{Ker} \; \tau_{_{W}}$ and suppose $x = Z \wedge W = Z \wedge
W'$ where $W' = W + \lambda Z$.  Then
\[ \tau_{_{W'}}b = (\tau_{_{W}} + \lambda \tau_{_{Z}})b =0,\]
so that $\textup{Ker} \; \tau_{_{Z}} \cap \textup{Ker} \;
\tau_{_{W}} = \textup{Ker} \; \tau_{_{Z}} \cap \textup{Ker} \;
\tau_{_{W'}}$.  Moreover, if $Z' = Z + \lambda W'$ we see that
\[ \tau_{_{Z'}}b = (\tau_{_{Z}} + \lambda \tau_{_{W'}})b =0,\]
so that $\textup{Ker} \; \tau_{_{Z}} \cap \textup{Ker} \;
\tau_{_{W}} = \textup{Ker} \; \tau_{_{Z'}} \cap \textup{Ker} \;
\tau_{_{W'}}$, establishing the second claim.  The third follows
similarly.
\endproof

As in \cite{mw:isdtt}, we write $U$ for the set of $x \in \hCM$ on
which $\Delta_x$ is invertible.

\begin{prop}
For all $x \in U$ we have the decomposition,
\begin{equation}
\label{trivial bundle decomposition}
B = E_x \oplus \textup{Im} \; \rho_{_{Z}} \oplus \textup{Im} \;
\rho_{_{W}}.
\end{equation}
In particular, for all $x \in \hCM$ we have $F_x = 0$.
\end{prop}
\proof  For $x \in U$, define
$$P_x = 1 - \rho_{_{W}} \Delta^{-1}_x \tau_{_{Z}} + \rho_{_{Z}}
\Delta^{-1}_x \tau_{_{W}} \; : B \rightarrow B.$$
Now $P_x$ is linear in $x$ and is dependent only on $x$ and not on
$Z,W$ individually.  It is easily shown that $P_x^2 = P_x$ and $P_x
B = E_x$, whence $P_x$ is the projection onto $E_x$. Moreover,
$$-\rho_{_{W}} \Delta^{-1}_x \tau_{_{Z}}, \qquad \rho_{_{Z}} \Delta^{-1}_x \tau_{_{W}},$$
are the projections onto the second and third summands respectively.
Hence we have proven the claim provided we can show that $F_x = 0$,
since the sum $\textup{Im} \; \rho_{_{Z}} \oplus \textup{Im} \;
\rho_{_{W}}$ is then direct.

Suppose that $F_x = (\rho_{_{Z}} A) \cap (\rho_{_{W}} A)$ is not
zero for some $x \in \CM^\#$, so there exists a non-zero $b \in
\rho_{_{W'}} A$ for all $W' \in \tilde T$ such that $x = Z \wedge W
= Z \wedge W'$.   Then in particular $a:= \rho_{_{W'}}^{-1}b$ is
non-zero and defines a holomorphic section of a vector bundle over
the two-dimensional subspace of $\tilde T$ spanned by all such $W'$
(whose typical fibre is just $A$), and hence (at the projective
level) a non-zero holomorphic section of the bundle $\CO(-1) \otimes
A$ over $\hat{x} = \CP^1$, where $\CO(-1)$ denotes the tautological
line bundle over $\CP^1$.  It is however well-known that this bundle
has no non-zero global sections, whence we must in fact have $F_x =
0$. \endproof

This procedure has thus constructed a rank $n$ vector bundle over
$U$ whose fibre over $x \in U$ is $E_x$ (again noting that the
construction is independent of the scaling of $x \in \TCM$). The
bundle $E$ is obtained as a sub-bundle of the trivial bundle $U
\times B$: the projection $P_x$ identifies the fibre of $E$ at each
$x \in \CM^\#$ as well as defining a connection on $E$ by orthogonal
projection of the trivial connection on $U \times B$.

\subsection{ADHM in the $*$-algebra picture}
In this section we mention how the ADHM construction ought to
operate in our $\SU_4$ $*$-algebra framework.  In passing to the
affine description of our manifolds, we encode them as real rather
than complex manifolds, obtaining a global coordinate algebra
description. Thus we expect that the ADHM construction ought to go
through at some global level in our $*$-algebra picture.  For now we
suppress the underlying holomorphic structure of the bundles we
construct, with the complex structure to be added elsewhere.

Indeed, we observe that the key ingredient in the ADHM construction
is the decomposition (\ref{trivial bundle decomposition}),
\[ B = E_x \oplus \textup{Im} \; \rho_{_{Z}} \oplus \textup{Im} \;
\rho_{_{W}},\]
which identifies the required bundle over space-time as a sub-bundle
of the trivial bundle with fibre $B$.  We wish to give a version of
this decomposition labelled by points at the projective level, that
is in terms of points of $T=\CP^3$ and $\hCM = \F_2(\C^4)$, rather
than in terms of the homogeneous representatives $Z \in \tilde T$
and $x \in \TCM$ used above. We shall do this as before by
identifying points of $\CP^3$ with rank one projectors $Q$ in
$M_4(\C)$, similarly points of $\hCM$ with rank two projectors $P$.
Points of the correspondence space are identified with pairs of such
projectors $(Q,P)$ such that $QP=Q=PQ$.

Recall that in the previous description, given $x \in \Lambda^2
\tilde T$ and $Z \in \hat{x}$ there are many $W$ such that $x = Z
\wedge W$.  In the alternative description, given a projection $P
\in \hCM$, the corresponding picture is that there are $Q, Q' \in T$
with $QP=Q=PQ$ and $Q'P=Q'=PQ'$ such that $\Im \, P = \Im \, Q
\oplus \Im\, Q'$. Given $P \in \hCM$, $Q \in T$ with $(Q,P) \in \F$,
there is a canonical choice for $Q'$, namely $Q':=P-Q$. Indeed, $P$
identifies a two-dimensional subspace of $\C^4$ and $Q$ picks out a
one-dimensional subspace of this plane.  The projector $Q'=P-Q$
picks out a line in $\C^4$ in the orthogonal complement of the line
determined by $Q$.

As in the monad description of the construction outlined earlier,
the idea is to begin with a trivial bundle of rank $2k+n$ with
typical fibre $B=\C^{2k+n}$ and to present sufficient information to
canonically identify a decomposition of $B$ in the form
(\ref{trivial bundle decomposition}).  In the monad description this
was done by assuming that $\tau_{_{Z}} \rho_{_{Z}} = 0$ and that the
maps $\tau, \rho$ were linearly dependent on $Z \in \tilde T$.  Of
course, this makes full use of the additive structure on $\tilde T$,
a property which we do not have in the projector version: here we
suggest an alternative approach.

In order to obtain such a decomposition it is necessary to determine
the reason for each assumption in the monad construction and to then
translate this assumption into the projector picture.  The first
observation is that the effect of the map $\rho$ is to identify a
$k$-dimensional subspace of $B$ for each point in twistor space (for
each $Q \in T$ the map $\rho_{_{Q}}$ is simply the associated
embedding of $A$ into $B$). We note that this may be achieved
directly by specifying for each $Q \in T$ a rank $k$ projection
$\rho_{_{Q}}:B \rightarrow B$.

The construction then requires us to decompose each point of
space-time in terms of a pair of twistors.  As observed above, given
$P \in \hCM$ and any $Q \in T$ such that $PQ=Q=QP$ we have $P=Q+Q'$,
where $Q' = P-Q$ (it is easy to check that $Q'$ is indeed another
projection of rank one). The claim is then that the corresponding
$k$-dimensional subspaces of $B$ have zero intersection and that
their direct sum is independent of the choice of decomposition of
$P$.  In the monad construction this was obtained using the assumed
linear dependence of $\rho$ on points $Z \in \tilde T$, although in
terms of projectors this is instead achieved by assuming that the
projections $\rho_{_{Q}}$,$\rho_{_{Q'}}$ are orthogonal whenever
$Q,Q'$ are orthogonal (clearly if $Q+Q'=P$ then $Q$,$Q'$ are
orthogonal), i.e.
\[\rho_{_{Q}}\rho_{_{Q'}}=\rho_{_{Q'}}\rho_{_{Q}} = 0\]
for all $Q,Q'$ such that $QQ'=Q'Q = 0$.  Moreover, we impose that
$\rho_{_{Q}} + \rho_{_{Q'}}$ depends only on $Q+Q'$, so that the
direct sum of the images of these projections depends only on the
sum of the projections.

Then for each $P \in \hCM$ we define a subspace $B_P$ of $B$ of
dimension $n$,
\[ B = B_P \oplus \Im \,\rho_{_{Q}} \oplus \Im \,\rho_{_{Q'}}\]
by constructing the projection
\[e_{_P} := 1 - \rho_{_{Q}} - \rho_{_{Q'}}\]
on $\C^4$, which is well-defined since the assumptions we have made
on the family $\rho_{_{Q}}$ imply that the projectors $e_{_P}$,
$\rho_{_{Q}}$ and $\rho_{_{Q'}}$ are pairwise orthogonal. This
constructs a rank $n$ vector bundle over $\hCM$ whose fibre over $P
\in \hCM$ is $\Im \, e_{_P} = B_P$.

\subsection{The tautological bundle on $\hCM$ and its corresponding monad}  We illustrate these ideas by
constructing a specific example in the $*$-algebra picture.  We
construct a `tautological monad' which appears extremely naturally
in the projector version and turns out to correspond to the
1-instanton bundle of Section~3.2.

We once again recall that compactified space-time may be identified
with the flag variety $\F_2(\C^4)$ of two-planes in $\C^4$ and this
space has its associated tautological bundle whose fibre at $x \in
\F_2(\C^4)$ is the two-plane in $\C^4$ which defines $x$ (it is
precisely this observation which gave rise to the projector
description of space-time in the first place).  Then we take
$B=\C^4$ in the ADHM construction: note that we expect to take
$n=2,k=1$, which agrees with the fact that $\dim B = 2k+n=4$.

Then for each point $P \in \hCM$ we are required to decompose it as
the sum of a pair $Q,Q'$ of rank one projectors (each representing a
twistor).  Here this is easy to do:  we simply choose any
one-dimensional subspace of the image of $P$ and take $Q$ to be the
rank one projector whose image is this line.  As discussed, the
canonical choice for $Q'$ is just $Q':=P-Q$.

In doing so, we have tautologically specified the one-dimensional
subspace of $B=\C^4$ associated to $Q \in T$ (recalling that $k=1$
here), simply defining $\rho_{_{Q}} :=Q$.  We now check that these
data satisfy the conditions outlined in the previous section. It is
tautologically clear that for all $Q_1,Q_2 \in T$ we have that
$\rho_{_{Q_1}}$ and $\rho_{_{Q_2}}$ are orthogonal projections if
and only if $Q_1$ and $Q_2$ are orthogonal.  Moreover, if we fix $P
\in \hCM$, $Q \in T$ such that $PQ=Q=QP$ and take $Q' = P-Q$ then
\[ \rho_{_{Q}} + \rho_{_{Q'}} = Q + Q' = P,\]
which of course depends only on $Q+Q'$ rather than on $Q$, $Q'$
individually.

Thus we construct the subspace $B_P$ as the complement of the direct
sum of the images of $\rho_{_{Q}}$ and $\rho_{_{Q'}}$.  As
explained, this is done by constructing the projection
\[e_{_P} = 1- \rho_{_{Q}} - \rho_{_{Q'}} = 1-P.\]
Thus as $P$ varies we get a rank two vector bundle over $\hCM$,
which is easily seen to be the complement in $\C^4$ of the
tautological vector bundle over $\hCM$.  We equip this bundle with
the Grassmann connection obtained by orthogonal projection of the
trivial connection on the trivial bundle $\hCM \times \C^4$, as
discussed earlier.

This gives a monad description of the tautological bundle over $\hCM
= \F_2(\C^4)$. As explained, the instanton bundle over $S^4$ is
obtained by restriction of this bundle to the two-planes $x \in
\F_2(\C^4)$ which are invariant under the map $J$ defined in
Section~3.2.  Hence we consider this construction only for $x \in
S^4$, and by Proposition~3.11 we have that $Q'=J(Q)$ in the above.
We now wish to give the monad version of the corresponding bundle
over twistor space, for which we need the map $\tau_{_{Q}}$.  We
recall that $\tau$ is meant to satisfy $\tau_{_{Q}} \rho_{_{Q}}=0$
for all $Q \in \CP^3$, and we use this property to construct $\tau$
by putting $\tau_{_{Q}}:= \rho_{_{J(Q)}}=J(Q)$, so that in this
tautological example we have
\[ \tau_{_{Q}} \rho_{_{Q}} = \rho_{_{J(Q)}} \rho_{_{Q}} =
J(Q)Q = 0.\]
The bundle over twistor space corresponding to the instanton then
appears as the vector bundle whose fibre over $Q \in \CP^3$ is the
cohomology
\[\tilde \E = \textup{Ker} \; \tau_{_{Q}} / \textup{Im} \; \rho_{_{Q}} = \textup{Ker} \; \rho_{_{J(Q)}} / \textup{Im} \; \rho_{_{Q}},\]
which in the case of the 1-instanton is the rank two bundle $\tilde
E = \textup{Ker} \; J(Q) / \textup{Im} \; Q$, and one may easily
check that this bundle over twistor space agrees with the one
computed in Section~3.2 using the Penrose-Ward transform.  Indeed,
the crucial property is that it is trivial upon restriction to $\hat
P = p(q^{-1}(P))$ for all $P \in S^4$, which is straightforward to
see through the observation that as $Q$ varies with $P$ fixed, $Q$
and $J(Q)$ always span the same plane (the one defined by $P$), and
the fibre of $\tilde E$ over all such $Q$ is precisely the
orthogonal complement to this plane.

\subsection*{Acknowledgements}
We would like to thank G. Landi and S.P. Smith for helpful comments
on versions of the manuscript during their respective visits to the
Noncommutative Geometry Programme at the Newton Institute.

\end{document}